%% file: main.tex
\documentclass[preprint,11pt]{article}

\usepackage{mathtools}
\usepackage{amsmath,amsfonts,amsthm,bm}
\usepackage{amssymb}
\usepackage[margin=1.0in]{geometry}
\usepackage{graphicx}
\usepackage{array}
\usepackage[normalem]{ulem}
\usepackage{cases}
\usepackage{comment}
\usepackage{tcolorbox}
\usepackage{float}
\usepackage{indentfirst}
\usepackage{lineno}
\usepackage{cases}
\usepackage[font=scriptsize,labelfont=bf]{caption} 
\usepackage[font=scriptsize,subrefformat=parens]{subcaption}
\usepackage{multirow}
\usepackage{tikz}
\usepackage{tikz-cd}
\usepackage[nottoc]{tocbibind}
\usetikzlibrary{matrix, calc, arrows}

\usepackage{lipsum}
\usepackage{fonts}
\usepackage{algorithm}
\usepackage{algpseudocode}
\usepackage{hyperref} \usepackage{cleveref}
\hypersetup{breaklinks,colorlinks,linkcolor=blue,citecolor=blue}
\usepackage{authblk}
\usepackage{placeins}

\usepackage{newtxtext,newtxmath}
\usepackage{fonts}

\newtheorem{remark}{Remark}

\newcommand{\p}{\uppartial}
\newcommand{\WS}{\mathcal W}

\newcommand{\Tr}{{\mathsf{T}}}

\newcommand{\diff}{\mathop{}\!\mathrm{d}}
\newcommand\figref{Figure~\ref}
\newcommand\tableref{Table~\ref}

\newcommand{\fc}{f_{\rm{c}}}
\newcommand{\fck}{f_{{\rm{c}},k} }
\newcommand{\Fc}{F_{\rm{c}}}
\newcommand{\bfc}[1][]{ {\bm{f}}^{ #1}_{\rm{c}} }
\newcommand{\bqc}[1][]{ {\bm{q}}^{ #1}_{\rm{c}} }

\newcommand{\fu}{f_{\rm{u}}}

\newcommand{\Fu}{F_{\rm{u}}}
\newcommand{\bfu}[1][]{\bm{f}^{#1}_{\rm{u}}}
\newcommand{\bqu}[1][]{ {\bm{q}}^{ #1}_{\rm{u}} }

\newcommand{\xl}{x_{\rm{L}}}
\newcommand{\xr}{x_{\rm{R}}}

\newcommand{\vmax}{v_{\rm{max}}}

\usepackage{tcolorbox}

\newtcbox{\rfbox}[1][gray!20]{on line,
	colback=#1, colframe=black, boxsep=0pt, boxrule=0pt, size=small, arc=2mm}

\newtcbox{\sfbox}[1][white]{on line,
	colback=#1, colframe=black, boxsep=0pt, boxrule=0pt, size=small}

\newcommand\reals{{{\rm l} \kern -.15em {\rm R} }}

\date{}
\title{A Collision-Based Hybrid Method for the BGK Equation}

\author[a]{Minwoo Shin}
\author[b]{Cory D. Hauck\thanks{Corresponding author at Oak Ridge National Laboratory: \url{hauckc@ornl.gov}}}
\author[c]{Ryan G. McClarren\thanks{Corresponding author at the University of Notre Dame: \url{rmcclarr@nd.edu}}}

\affil[a]{Yonsei University, School of Mathematics and Computing, Seoul, 03722, Republic of Korea}
\affil[b]{Computer Science and Mathematics Division, Oak Ridge National Laboratory, Oak Ridge, TN 37831}
\affil[c]{Department of Aerospace and Mechanical Engineering, University of Notre Dame, Fitzpatrick Hall, Notre Dame, IN 46556, USA}

\begin{document}
\allowdisplaybreaks
\maketitle
\begin{abstract}
We apply the collision-based hybrid introduced in \cite{hauck} to the Boltzmann equation with the BGK operator and a hyperbolic scaling.  An implicit treatment of the source term is used to handle stiffness associated with the BGK operator. Although it helps the numerical scheme become stable with a large time step size, it is still not obvious to achieve the desired order of accuracy due to the relationship between the size of the spatial cell and the mean free path. Without asymptotic preserving property, a very restricted grid size is required to resolve the mean free path, which is not practical. Our approaches are based on the noncollision-collision decomposition of the BGK equation. We introduce the arbitrary order of nodal discontinuous Galerkin (DG) discretization in space with a semi-implicit time-stepping method; we employ the backward Euler time integration for the uncollided equation and the 2nd order predictor-corrector scheme for the collided equation, i.e., both source terms in uncollided and collided equations are treated implicitly and only streaming term in the collided equation is solved explicitly. This improves the computational efficiency without the complexity of the numerical implementation. Numerical results are presented for various Knudsen numbers to present the effectiveness and accuracy of our hybrid method. Also, we compare the solutions of the hybrid and non-hybrid schemes.
\end{abstract}



\input{sec1_intro}
\input{sec2_formulation}
\input{sec3_implementation_details}
\input{sec4_results}
\input{sec5_conclusion}
\section*{Declaration of competing interest}
The authors declare that they have no known competing financial interests or personal relationships that could have appeared to influence the work reported in this paper.
\section*{Acknowledgements}
This research is based upon work supported by the U.S. Department of Energy, Office of Science, Office of Advanced Scientific Computing Research, as part of their Applied Mathematics Research Program. The work was performed at the Oak Ridge National Laboratory, which is managed by UT-Battelle, LLC under Contract No. De-AC05-00OR22725. The United States Government retains and the publisher, by accepting the article for publication, acknowledges that the United States Government retains a non-exclusive, paid-up, irrevocable, world-wide license to publish or reproduce the published form of this manuscript, or allow others to do so, for the United States Government purposes. The Department of Energy will provide public access to these results of federally sponsored research in accordance with the DOE Public Access Plan (http://energy.gov/downloads/doe-public-access-plan). In addition, this research is supported by the National Research Foundation of Korea (NRF) grant funded by the Korea government (MSIT) (No. RS-2023-00220762).
\FloatBarrier
\bibliographystyle{unsrt}
\bibliography{reference.bib}
\end{document}

%% file: sec1_intro.tex
\section{Introduction}

The Bhatnagar-Gross-Krook (BGK) equation is a well-known kinetic theory model used to simulate the non-equilibrium behavior of dilute gases. It is an approximation of the Boltzmann equation that replaces the Boltzmann collision operator, a nonlinear integral operator modeling binary collisions, by the BGK operator \cite{bhatnagar,perthame}, a simpler nonlinear relaxation model. This modification significantly reduces the computational cost of kinetic simulations of dilute gases while still recovering both the equilibrium and streaming behavior of the Boltzmann equation in collision-dominated and collisionless limits, respectively.  The BGK operator also preserves the conservation and entropy dissipation properties  of the Boltzmann operator, and recent work has extended these properties to the multispecies setting \cite{klingenberg2017consistent,HHKPW,HHM}.  

Approximation with the BGK equation does have two major drawbacks.  First, it does not capture all of the transport properties of the Boltzmann equation in collisional regimes.  Most notably, it fails to capture the correct Prandtl number (essentially the ratio of viscosity to thermal conductivity).  Hence it may not agree with the compressible Navier-Stokes equations that are derived from the Boltzmann equation in collisional regimes.  However, this problem can be remedied via a slight modification that leads to an ellipsoidal statistical operator \cite{holway1965kinetic,andries2000gaussian}, sometimes referred to as ES-BGK \cite[Section 3.6]{struchtrup2005macroscopic}.  The second drawback is that the structural properties of the BGK operator rely on the assumption that the collision frequency is independent of the microscopic velocity, which is known to be unphysical \cite{struchtrup2005macroscopic, LeeMore}.  Extensions to velocity-dependent frequencies have been explored in single species \cite{Struchtrup97} and multi-species \cite{HHKPW, haack2022numerical} settings.  While these models add more physical realism, they come with a significant increase in computational cost.

\subsection{Related work}
On the numerical side, the BGK equation has been the object of extensive study.  Aside from the cost of discretizing phase space, there are two primary challenges, both of which occur in collisional regimes.  The first challenge is the stiffness of the BGK operator, which suggests an implicit treatment in order to avoid excessively small time steps. (See, however, \cite{melis2019projective} for a recent alternative approach.)  The second challenge is to recover---at the numerical level---consistent numerical discretizations of the Euler and Navier-Stokes system for compressible gas flow that are derived from the BGK model via a Chapman-Enskog expansion \cite[Chapter 4]{struchtrup2005macroscopic}.   

In \cite{coron}, it was shown that a backward Euler time discretization of the BGK operator could be implemented in an explicit fashion.  This numerical ``trick" enables stable time stepping for large time steps across a wide range of collisional regimes, and for this reason, it has become the basis of many subsequent numerical schemes.  In \cite{yang}, various discretizations of the BGK and Shaskov model \cite{shakhov1968generalization} are considered and compared.%
\footnote{The Shaskov model is designed to capture the correct Prandtl number, but unlike BGK, it is not guaranteed to generate a positive solution.}   In \cite{bennoune}, a method based on the decomposition of the kinetic distribution into macroscopic (equilibrium) and microscopic (non-equilibrium) parts is proposed in order to recover a consistent Navier-Stokes limit when the Knudsen number is sufficiently small.  Simulations of the ES-BGK model for small Knudsen numbers have also been developed \cite{mieussens,andries,hu} in order to recover the correct  Navier-Stokes equations at the numerical level.  There has also been extensive research on implicit–explicit schemes for the BGK equations and ES-BGK equation; see for instance \cite{pieraccini1,pieraccini2,filbet2}. A gas-kinetic scheme for the BGK model was proposed and compared to the compressible Navier–Stokes equations in the collisional regime \cite{xu}.

Even when the collision operator is treated implicitly, an explicit treatment of the advection operator in the BGK equation can lead to restrictive time steps.  There are two common reasons:  long time scales (i.e. a diffusive-type scaling) that lead to incompressible equations in the infinite collisional limit \cite{bardos1991fluid} or problems for which the maximum velocity in the computational domain is significantly larger than the fluid sound speed.  One way to address such restrictions is to take a fully implicit approach, which is common for time-dependent kinetic equations in radiation transport contexts \cite[Section 1.2]{larsen2010advances} and has been proposed for electron transport problems in \cite{decaria2021analysis,laiu2020fast}.  In both these cases, diffusive limits are important. Fully implicit methods have also been proposed for collisional dilute gases  \cite{taitano2014moment, bobylev1995,mieussens2000discrete} and collisionless plasmas \cite{garrett2018fast}.  All of these approaches require sophisticated iterative strategies to manage the cost and memory requirement of the implicit update.   An alternative is to use semi-Lagrangian methods \cite{ding2021semi,cho2021conservative, groppi2014high}, which rely on a characteristic formulation.  There the main effort lies in the interpolation of point-wise values since characteristics of the advection operator in the BGK equation do not always intersect grid points.

\subsection{Hybrid approach:  benefits and drawbacks}

In the current paper, we focus on the BGK equation in a hyperbolic scaling and propose a method to address stiffness due to large phase space velocities.  Rather than use an IMEX or fully implicit time-stepping strategy, we instead leverage the hybrid approach introduced in \cite{hauck}.  While this approach was originally proposed for linear transport equations with fully implicit time integration \cite{hauck,crockatt,crockatt1}, we use it here to solve the BGK model in a semi-implicit fashion \cite{McClarren2021,MCCLARREN20087561}.  

The hybrid strategy is based on the method of {\em first-collision source} \cite{Alcouffe1985,Alcouffe2017jt}, in which the kinetic equation is decomposed into the sum of an uncollided equation and a collided equation.  The uncollided equation is a damped linear advection equation with no in-scattering and hence can be solved implicitly rather easily with high resolution and in parallel;  meanwhile, the collided equation is expected to have a solution that can be approximated with a low-fidelity velocity discretization.  
 We use the Euler fluid approximation as the velocity discretization for the collided equation, but unlike previous applications of similar hybrid methods, the advection is solved explicitly with a 2nd-order Runge-Kutta (predictor-corrector) time integrator.   These choices are flexible;  for example, a Navier-Stokes model can be used for the collided equation and a semi-Lagrangian approach can be taken for the uncollided equation.
 
 There are two main benefits of the hybrid decomposition. The first benefit is the simplicity of the numeric implementation.  Indeed, if a fluid code is readily available, then one need only add a relatively simple implicit solver for the uncollided equations.  Moreover, as mentioned above, the choice of discretization for each equation can be flexible.  The second benefit is a less restrictive CFL condition that is based on the maximum eigenvalues of the flux Jacobian in the fluid solver rather than the maximum velocity in the computational domain. As a consequence, larger time steps can be taken.
 
 The main drawback of the hybrid approach is that it induces both general time-stepping errors and consistency errors in the coarse approximation of the collided equation that depend on the collisionality of the problem.   Because the method includes a remapping procedure after each time step, the mixing of errors can be difficult to track, making even a formal analysis difficult.  What is certain, however, is that the consistency error from the hybrid approximation will eventually cause saturation in convergence as the space-time mesh is resolved.  In our experiments, we observe that this saturation happens only at very low error levels, and if saturation does become a concern, there are ways to correct it \cite{crockatt, crockatt1}.

The remainder of the paper is organized as follows. We begin with a derivation of the uncollided-collided decomposition of a BGK equation in Section \ref{sec:formulation}. We present the discretization details for the uncollided and collided equations in Section \ref{sec:discretization} and perform an asymptotic analysis. 
In Section \ref{sec:results}, we present numerical results on several benchmark problems to demonstrate the simplicity, robustness, efficiency, asymptotic properties of the hybrid approach. We also compare the hybrid method to a previously developed implicit-explicit method. Finally, in Section \ref{sec:conclusion}, we draw conclusions and make suggestions for future work.

%% file: sec2_formulation.tex
\section{The hybrid formulation}\label{sec:formulation}
\subsection{The BGK equation}
Given a positive integer $d$, let $X \subset \bbR^d$ be an open bounded domain with Lipschitz boundary.  Let $\Gamma = \uppartial X \times \bbR^d$ and  $\Gamma^{\pm} = \{ (\bm{x},\bm{v}) \in \Gamma \colon \pm \bm{n}(\bm{x}) \cdot \bm{v}  >0 \}$, where $\bm{n}(\bm{x})$ is the outward normal with respect to $X$, defined at a.e. $\bm{x} \in \uppartial X$.  The BGK model for the \textit{kinetic distribution function} $F$ with a source term $S$ is
\begin{subequations}
\label{eq:bgk_system}
\begin{numcases}{}
   \uppartial_t F(\bm{x},\bm{v},t)+\bm{v}\cdot\nabla_{\bm{x}}{F}=\frac{1}{\epsilon}\left(\cM_F(\bm{x},\bm{v},t)-F(\bm{x},\bm{v},t)\right)+S(\bm{x},\bm{v},t),
      & $(\bm{x},\bm{v}) \in X \times \bbR^d, ~t>0$, \label{eq:bgk_eqn}\\
     F(\bm{x},\bm{v},t) = G^-(\bm{x},\bm{v},t),
    & $(\bm{x},\bm{v}) \in \Gamma^-, ~t>0$, \label{eq:bgk_bndy}\\
    F(\bm{x},\bm{v},0) = G^{0}(\bm{x},\bm{v}),
    & $(\bm{x},\bm{v}) \in X \times \bbR^d$, \label{eq:bgk_init}
\end{numcases} 
\end{subequations}
where the inflow $G^-$, initial condition $G^0$, and source term $S$ are given non-negative functions.  The function $F$, which depends on position $\bm{x} \in X$, velocity $\bm{v}\in\bbR^{d}$, and time $t \geq 0$, 
gives the mass density of particles with respect to the phase space measure $\bm{dv} \bm{dx}$.
The parameter $\epsilon > 0$ is the Knudsen number, i.e., the ratio of the mean free path between collisions to the length scale of the domain; in this work, we assume $\epsilon$ is a fixed constant.  The function $\cM_F$ is given by
\begin{equation}\label{eq:Maxwellian}
    \cM_F=\cM_F(\bm{x},\bm{v},t)=\frac{\rho_F (\bm{x},t)}{(2\pi  \theta_F(\bm{x},t))^{d/2}}\exp\left(-{\frac{|\bm{v}-\bm{u}_F(\bm{x},t)|^2}{2\theta_F(\bm{x},t)}}\right),\end{equation}
where the \textit{primitive variables} $\rho_F(\bm{x},t)\geq 0$, $\bm{u}_F(\bm{x},t) \in \bbR^d$, and $\theta_F(\bm{x},t) \geq 0$ are the local mass density, mean velocity, and temperature, respectively; functions of the form in \eqref{eq:Maxwellian} are referred to as \textit{Maxwellians}.

\subsection{Moment equations}
The local mass density $\rho_F$, momentum density $\bm{m}_F$, and total energy density $E_F$, are given by 
\begin{equation}
    \bm{q}_F:=\begin{pmatrix}
    \rho_F\\
    \bm{m_F}\\
    E_F
    \end{pmatrix}
    =\int_{\bbR^{d}}\bm{e}F \diff \bm{v}=\int_{\bbR^{d}}
    \begin{pmatrix}
    {F}\\
    {\bm{v}F}\\
    {\frac{1}{2}|\bm{v}|^2 F}
    \end{pmatrix}\diff \bm{v},
\end{equation}
where  the components of $\bm{e}=\bm{e}(\bm{v}):=(1,\bm{v},\frac{1}{2}|\bm{v}|^2)^T$ are referred to as \textit{collision invariants}, since
\begin{equation}
\label{eq:collision_invariants}
   \int_{\bbR^{d}} (\cM_F-F) \bm{e} dv =0.
\end{equation}
We also denote by $\mathcal{E}$ the map which takes a set of moments to the associated Mawellian, i.e., $\cM_F = \mathcal{E}(\bm{q}_F)$.

The condition in \eqref{eq:collision_invariants}, which is intended to mimic properties of the Boltzmann collision operator (see e.g. \cite[Chapter 3]{cercignani2013mathematical}), uniquely determines the primitive variables $\rho_F$, $\bm{u}_F$, $\theta_F$ that parameterize the Maxwellian in \eqref{eq:Maxwellian}.  It also establishes the relationship between the moments $ \bm{m}_F$ and $E_F$ and the variables $\bm{u}_F$ and $\theta_F$.  Specifically,
\begin{equation}
    \bm{m}_F=\rho_F \bm{u}_F \quad \quad \text{and} \quad
    E_F=\frac{1}{2}\rho_F |\bm{u}_F|^2+ \frac{d}{2}\rho_F \theta_F.
\end{equation}

If $F$ is a solution to the BGK equation \eqref{eq:bgk_eqn}, its moments satisfy the following system of conservation laws:
\begin{subequations}
\label{eq:cons_laws}
\begin{align}
\frac{\uppartial \rho_F}{\uppartial t}+\nabla_x\cdot(\rho_F \bm{u}_F)&=0,\\
\frac{\uppartial (\rho_F \bm{u}_F)}{\uppartial t}+\nabla_x\cdot(\rho_F \bm{u}_F\otimes  \bm{u}_F+\bm{P}_F)&=0,\\
\frac{\uppartial E_F}{\uppartial t}+\nabla_x\cdot(E_F\bm{u}_F+\bm{P}_F\bm{u}_F+Q_F)&=0,
\end{align}
\end{subequations}
where 
\begin{equation}
    \bm{P}_F=\int_{\bbR^{d}} (\bm{v}-\bm{u}_F)\otimes(\bm{v}-\bm{u}_F) F\diff  \bm{v} 
    \quad \text{and} \quad
    Q_F=\int_{\bbR^{d}}{\frac{1}{2}(\bm{v}-\bm{u}_F)|\bm{v}-\bm{u}_F|^2} F \diff \bm{v}
\end{equation}
are the pressure tensor and heat flux, respectively.  The equations in \eqref{eq:cons_laws}, which are not closed, are derived by integrating \eqref{eq:bgk_eqn} against the vector $\bm{e}$ of collision invariants and then invoking \eqref{eq:collision_invariants} to make the right-hand side vanish.  To close \eqref{eq:cons_laws}, one must specify formulas for $\bm{P}_F$ and $Q_F$ in terms of the moments $\rho_F$, $\bm{m}_F$, and $E_F$ or, equivalently, the primitive variables $\rho_F$, $\bm{m}_F$, and $\theta_F$.  In the fluid limit, when $\epsilon \to 0$, the scaling in \eqref{eq:bgk_eqn} suggests that $F$ can be approximated by $\mathcal{M}_F$ in the formulas for $\bm{P}_F$ and $Q_F$ (see for example \cite[Section 2.8]{cercignani1988boltzmann}), in which case
\begin{equation}
\label{eq:E_Fuler_approx}
    \bm{P}_F \approx \int_{\bbR^{d}}{(\bm{v}-\bm{u}_F)\otimes(\bm{v}-\bm{u}_F)} \mathcal{M}_F \diff  \bm{v} = \rho_F \theta_F \bm{I}
    \quad \text{and} \quad 
    Q_F \approx \int_{\bbR^{d}}{\frac{1}{2}(\bm{v}-\bm{u}_F)|\bm{v}-\bm{u}_F|^2} \mathcal{M}_F \diff \bm{v} = 0,
\end{equation}
where $\bm{I}$ is the $d\times d$ identity matrix. Using these approximations for $\bm{P}_F$ and $Q_F$ to close \eqref{eq:cons_laws} recovers the Euler equations for a gamma-law gas with pressure law $p (\rho,\theta)=\rho \theta$. Because a generic gamma-law gas has equation of state $p=\frac{d}{2}(\gamma-1)\rho_F\theta_F$ (see .g. \cite[pg.~13]{toro2013riemann}), it follows that $\gamma=\frac{d+2}{d}$.

\subsection{The hybrid method}\label{sec:hybrid}
We adopt the collision-based hybrid approach introduced in \cite{hauck} and decompose the distribution function $F$ as the sum of an uncollided component $F_{\rm{u}}$ and an uncollided component $F_{\rm{c}}$. 
Let $0 = t_0 < t_1 < \dots $ be a sequence of points in a temporal grid.  For  $t \in (t_n,t_{n+1})$, let $F_{\rm{u}}$ and $F_{\rm{c}}$ solve
\begin{subequations}
	\label{eq:bgku_system}
	\begin{numcases}{}
		\uppartial_t \Fu(\bm{x},\bm{v},t)+\bm{v}\cdot\nabla_{\bm{x}}{\Fu} + \frac{1}{\epsilon} \Fu(\bm{x},\bm{v},t)= S(\bm{x},\bm{v},t),
		\qquad
		& $(\bm{x},\bm{v}) \in X \times \bbR^d, ~t \in (t_n, t_{n+1})$, \label{eq:bgku_eqn}\\
		\Fu(\bm{x},\bm{v},t) = G^-(\bm{x},\bm{v},t),
		& $(\bm{x},\bm{v}) \in \Gamma^-, ~t \in (t_n, t_{n+1})$, \label{eq:bgku_bndy}\\
		\Fu(\bm{x},\bm{v},t_n^+) = G^n(\bm{x},\bm{v}),
		& $(\bm{x},\bm{v}) \in X \times \bbR^d$, \label{eq:bgku_init}
	\end{numcases} 
\end{subequations}
\begin{subequations}
	\label{eq:bgkc_system}
	\begin{numcases}{}
		\uppartial_t \Fc(\bm{x},\bm{v},t)+\bm{v}\cdot\nabla_{\bm{x}}{\Fc} + \frac{1}{\epsilon} \Fc(\bm{x},\bm{v},t)
		= \frac{1}{\epsilon} \mathcal{M}_F,
		& $(\bm{x},\bm{v}) \in X \times \bbR^d, ~t \in (t_n, t_{n+1})$, \label{eq:bgkc_eqn}\\
		\Fc(\bm{x},\bm{v},t) = 0,
		& $(\bm{x},\bm{v}) \in \Gamma^-, ~t \in (t_n, t_{n+1})$, \label{eq:bgkc_bndy}\\
		\Fc(\bm{x},\bm{v},t_n^+) = 0,
		& $(\bm{x},\bm{v}) \in X \times \bbR^d$, \label{eq:bgkc_init}
	\end{numcases} 
\end{subequations}
where $G^{0}(\bm{x},\bm{v})$ is given as in \eqref{eq:bgk_init} and for $n \geq 1$,
\begin{equation}
	\label{eq:G_n}
	G^n(\bm{x},\bm{v}) = \Fu(\bm{x},\bm{v},t_n^-) +  \Fc(\bm{x},\bm{v},t_n^-),
\end{equation}
and $S$ is a source term. We assume the source $S$ to be zero unless otherwise stated. The initial data and boundary conditions of the original system are inherited by the uncollided equation, while the collided equation is assigned zero inflow and initial condition.  For $s \geq 1$, the uncollided solution is reinitialized by adding together the current values of $\Fu$ and $\Fc$.  However, no splitting error is incurred by writing \eqref{eq:bgk_system} as the sum of \eqref{eq:bgku_system} and \eqref{eq:bgkc_system} until the equations are further discretized.   

The basic idea of the hybrid approach is to discretize \eqref{eq:bgku_system} with a high-resolution method in velocity with the expectation that collisions will allow \eqref{eq:bgkc_system} to be solved with a low-resolution method in velocity without significantly degrading the accuracy of the solution over a time step.  However, in the limit of zero collisions, the uncollided equation is exact and the collided equation plays no role.  At the end of the time step, the collided solution is reconstructed in the high-resolution space and then added to the uncollided solution in order to approximate $G^n$ in \eqref{eq:G_n}.  The combined solution is then used to reinitialize the uncollided equation at the next time step, while the initial condition for the collided equation is again reset to zero. To advance numerically within each time step, we employ a predictor-corrector method, the details of which are given in Section \ref{sec:time discretization}.

In the original hybrid formulation \cite{hauck}, the procedure summarized above resulted in a significant improvement in the time-to-solution for a fully implicit time-stepping scheme when applied to a linear, kinetic transport equation.  In the current context, we consider hyperbolic time scales of the bulk flow which typically do not require fully implicit approaches, unless the mesh is highly unstructured.  Instead, only \eqref{eq:bgku_system} is solved fully implicitly in order to step over hyperbolic time-scale restrictions imposed by high-velocity components of the distribution.  Meanwhile, \eqref{eq:bgkc_system} can be solved using an approximate model with a less restrictive time step.

\section{Discretization of the hybird method}\label{sec:discretization}
For the purposes of this paper, we consider the one-dimensional case, which implies that the gas constant $\gamma=3$.  We set $X = (\xl,\xr)$ and restrict the velocity domain to a bounded set $V = [-\vmax,\vmax]$. 
\subsection{Velocity discretization}
Many velocity discretization methods can easily fit into the hybrid strategy outlined above.  In the current paper, we choose to solve \eqref{eq:bgku_eqn} with a discrete velocity method \cite{mieussens,mieussens1,palczewski1997,bobylev1995} and \eqref{eq:bgkc_eqn} using the Euler equations as an moment-based velocity approximation.
\subsubsection{Discrete velocity model (DVM) for the uncollided equation}
Let $\{v_k\}_{k=1}^{N_v}$ and $\{\omega_k\}_{k=1}^{N_v}$ be the points and weights, respectively, of a one-dimensional,  $N_v$-point Gauss-Legendre quadrature set scaled to a truncated velocity domain $[-\vmax, \vmax]$, and for each $k \in {1,\dots,N_v}$, let $\gamma_k^{\pm}= \{ x \in X : \pm n(x) v_k  > 0\}$.  Then the discrete velocity model for \eqref{eq:bgku_eqn} which governs $f_{{\rm u},k}(x,t) \simeq F_{\rm{u}}(x,v_k,t)$ is given by
\begin{subequations}
\label{eq:dvm}
	\begin{numcases}{}
	\label{eq:dvm_eqn}
	\uppartial_t f_{{\rm u},k}(x,t)+v_k \uppartial_x{f_{{\rm u},k}(x,t)} + \frac{1}{\epsilon}f_{{\rm u},k}(x,t) = 0,
		& $x \in X, ~~t \in (t_n,t_{n+1})$, \\
		 \label{eq:dvm_bc}
		f_{{\rm u},k}(x,t) = G^-(x,{v_k},t),
		& $x \in \gamma_k^-,  ~~t \in (t_n,t_{n+1})$,\\
		 \label{eq:dvm_ic}
		f_{{\rm u},k}(x,t_n^+) = 	g^n_k(x)  ,
		& $x \in X$.
	\end{numcases} 
\end{subequations}
where $g^0_k(x)   = G^0(x,v_k)$ and, for $n \geq 1$,  
\begin{equation}
g^n_k(x) = f_{{\rm u},k}(x,t_n^-) + f_{{\rm c},k}(x, t_n^-), 
\end{equation}
with $\fc$ defined in \eqref{eq:fc_approx}  below.
The approximate uncollided moments ${\bm{q}}_{\rm{u}} \approx \bm{q}_{F_{\rm_u}}$ are obtained from $f_{{\rm u},k}$ via the quadrature formula
\begin{equation}
	\label{eq:qu_quad}
	\bm{q}_{\rm{u}}(x,t) 
	= \sum_{k=1}^{N_v}{\omega}_k \bm{e}_k f_{{\rm{u}},k}(x,t),
\end{equation}
where $\bm{e}_k=(1,v_k,\frac12|v_k|^2)^{\Tr}$.  In vectorized form, the discrete velocity equation \eqref{eq:dvm_eqn} and the initial condition \eqref{eq:dvm_ic}  can be written as
\footnote{Unfortunately it is difficult to specify \eqref{eq:dvm_bc} cleanly in vectorized form since $\gamma^-_k$ depends on $k$.}
\begin{subequations}
	\label{eq:dvm_vec}
	\begin{numcases}{}
		\p_t \bfu(x,t)+ \bm{V} \p_x \bm{f}_{\rm{u}} (x,t)+ \frac{1}{\epsilon} \bfu(x,t) = 0,
		& $x \in X, ~~t \in (t_n,t_{n+1})$, \label{eq:dvm_vec_eqn}\\
		\bfu(x,t_n^+) = 	\bm{g}^n,
		& $x \in X$, \label{eq:dvm_vec_ic}
	\end{numcases} 
\end{subequations}
where $\bm{f}_{\rm{u}} = (\bm{f}_{{\rm{u}},1} , \dots , \bm{f}_{{\rm{u}},{N_v}})^\Tr$, 
$\bm{g}^n= (g^n_1 , \dots , g_{N_v}^n)^\Tr$, 
and $\bm{V}= \operatorname{diag}(v_1 , \dots , v_{N_v})$. %

\subsubsection{Moment equations for the collided equation}
The exact equations for the mass, momentum, and energy moments of ${F_{\rm{c}}}$ are given by the quadrature formula
\begin{equation}\label{eq:collided_moments}
    \uppartial_t \bm{q}_{F_{\rm{c}}}+ \uppartial_x \int_{\bbR} {v \bm{e} F_c dv} =\frac{1}{\epsilon}\bm{q}_{F_{\rm{u}}},
\end{equation}
where $\bm{q}_{F_{\rm{c}}} = \int_{\bbR} \bm{e} {F_{\rm{c}}} dv$.
We approximate the evolution of  $\bm{q}_{F_{\rm{c}}}$ by a function $\bm{q}_c$ that satisfies the Euler approximation to  \eqref{eq:collided_moments};  that is, 
\begin{equation}\label{eq:collided_moments1}
    \uppartial_t \bm{q}_{\rm{c}}+ \uppartial_x \bm{\Phi}(\bm{q}_{\rm{c}}) =\frac{1}{\epsilon}\bm{q}_{\rm{u}},
    \quad \text{where} \quad    
    \bm{\Phi}(\bm{q}_{\rm{c}}) 
    = \int_{\bbR} {v \bm{e} \cE({\bm{q}_c}) dv} ,
\end{equation}
and $\bm{q}_{\rm{u}}$ is computed according to \eqref{eq:qu_quad}.  With initial and boundary conditions, the complete system is
\begin{subequations}
	\label{eq:Eulerc}
	\begin{numcases}{}
		 \uppartial_t \bm{q}_{\rm{c}}(x,t)+ \uppartial_x \bm{\Phi}(\bm{q}_{\rm{c}}(x,t)) =\frac{1}{\epsilon}\bm{q}_{\rm{u}}(x,t),
		& \quad $x \in X, ~~t \in (t_n,t_{n+1})$, \label{eq:Eulerc_eqn}\\
		\bm{\Phi}(\bm{q}_{\rm{c}}(x,t))=   \int_{\{v:n(x) v > 0 \}}  v \bm{e} \cE(\bm{q}_{\rm{c}}(x,t)) dv,
		& \quad $x \in \p X,  ~~t \in (t_n,t_{n+1})$,  \label{eq:Eulerc_bc}\\
		\bm{q}_{\rm{c}}(x,t_n^+) = 0,
		& \quad $x \in X$,\label{eq:Eulerc_ic}
	\end{numcases} 
\end{subequations}
and the approximate kinetic solution is 
\begin{equation}
\label{eq:fc_approx} 
	 f_{{\rm c}, k}(x, t) =  \cE(\bm{q}_{\rm{c}}(x, t))(v_k).
\end{equation}
\begin{remark}
	The boundary condition in \eqref{eq:Eulerc_bc} is derived by separating the flux integral in \eqref{eq:collided_moments} into incoming and outgoing data, using the kinetic data (which is zero) for the incoming data and the approximation $\Fc \approx \cE(\bm{q}_{\Fc})$ for the outgoing data; that is, for $x \in \uppartial X$, 
	\begin{equation}
		\begin{split}
	\bm{\Phi}(\bm{q}_{\rm{c}} (x,t))
	&\approx	\int_{\bbR} {v \bm{e}(v) \Fc (x,v,t)} dv \\
	&=
		\int_{\{v:n(x) v < 0\} }  v \bm{e}(v) \Fc(x,v,t) dv+ \int_{\{v:n(x) v > 0\} }  v \bm{e}(v) \Fc(x,v,t) dv \\
		 &\approx \int_{\{v:n(x) v > 0\} }  {v \bm{e} \cE(\bm{q}_{\Fc}(x,t))(v)}  dv
		 \approx \int_{\{v:n(x) v > 0\} }  {v \bm{e} \cE(\bm{q}_{\rm{c}}(x,t))(v) } dv.
		 \end{split}
	\end{equation}
There are other approaches to assigning boundary conditions, including asymptotic conditions derived via half-space problems; see for example \cite{bardos2006half}.
\end{remark}


\subsubsection{Euler Limit}
\label{subsec:continous-time-Euler-limit}

An important property of the hybrid is that it recovers (at least formally) the Euler equations in the limit $\epsilon \to 0$.   To investigate this limit, we apply the quadrature rule in \eqref{eq:qu_quad} to the discrete velocity equation in \eqref{eq:dvm} and add the result to the moment equation for $\bm{q}_{\rm{c}}$ in \eqref{eq:collided_moments1}.  This gives the following conservation law for $\bm{q} = \bm{q}_{\rm{u}} + \bm{q}_{\rm{c}}$:
\begin{subequations}
\label{eq:adding_moments}	
\begin{numcases}{}
\uppartial_t \bm{q} (x,t)  
 +    \uppartial_x \left( \sum_{k=1}^{N_v}{\omega}_k v_k \bm{e}_k  {f_{{\rm u},k}(x,t)} \right)
 +  \p_x \bm{\Phi}(\bm{q}_{\rm{c}} (x,t) ) = 0,  & \, $x \in X, t \in (t_n,t_{n+1})$, \label{eq:adding_moments_eqn}	 \\
 \bm{\Phi}(\bm{q}(x,t))=   \sum_{\{k:  n(x) v_k < 0\}}{\omega}_k v_k \bm{e}_k  G^-(x,v_k) 
 + \int \limits_{\{v:n(x) v > 0\} }  v \bm{e} \cE(\bm{q}_{\rm{c}}(x,t))(v) dv,
 & \, $x \in \p X,  t \in (t_n,t_{n+1})$,  \label{eq:eq:adding_moments_bc}
 \label{eq:adding_moments_bc}\\
 \bm{q}(x,t_n^+) = \sum_{k=1}^{N_v}{\omega}_k  \bm{e}_k \left( f_{{\rm u},k}(x,t_n^-) + \cE(\bm{q}_{\rm{c}}(x, t))(v_k) \right),
 & \, $x \in X$ , \label{eq:eq:adding_moments_ic}
 \label{eq:adding_moments_ic}
\end{numcases}
\end{subequations}
which depends on $\epsilon$ via  $\bm{f}_{{\rm u}}$ and $\bm{q}_{\rm{c}}$.  To assess the limiting behavior of \eqref{eq:adding_moments}, we use the exact solution for  \eqref{eq:dvm}.  For each  $x \in X$, $k \in \{1,\dots,N_v\}$, and $n \geq0$ fixed, let 
\begin{equation}
	\label{eq:t-star}
	t^*:= t^*(x,v_k,t) = \min_{\tau > 0} \{ \tau: x - v_k(t-\tau)  \in \gamma^-_k\}.
\end{equation}
Then for $t \in (t_n,t_{n+1})$ and a.e. $x \in X$,
\begin{equation}
f_{{\rm u},k}(x,t) 
= \begin{cases}
e^{-\frac{t-t_n}{\epsilon}} g^n (x-v_k t_n,v_k) ,&  t - t^* < t_n,\\
e^{-\frac{t - t^*}{\epsilon}}G^-(x-v_k t^*,v_k, t - t^*) ,&  t -  t^* > t_n.
\end{cases}
\end{equation}
In either case above, $f_{{\rm u},k} \to 0$ as $\epsilon \to 0$ in the interior of the domain.  Thus the contribution of $\fu$ to the flux in \eqref{eq:adding_moments_eqn}	and \eqref{eq:adding_moments_ic} vanishes.  Furthermore, $\bm{q}_c  \equiv \bm{q} - \bm{q}_u \to  \bm{q}$ as $\epsilon \to 0$, so that \eqref{eq:adding_moments} recovers the Euler system for $\bm{q}$; that is, in the limit as $\epsilon \to 0$,
\begin{subequations}
	\label{eq:Euler_limit}	
	\begin{numcases}{}
		\label{eq:Euler_limit_eqn}	
		\uppartial_t \bm{q} (x,t)  
		+  \p_x \bm{\Phi}(\bm{q} (x,t) ) = 0,  & \quad $x \in X, t \in (t_n,t_{n+1})$,  \\
		\bm{\Phi}(\bm{q}(x,t))=   \sum_{\{k:  n(x) v_k < 0\}}{\omega}_k v_k \bm{e}_k  G^-(x,v_k) 
		+ \int \limits_{\{v:n(x) v > 0\} }  v \bm{e} \cE(\bm{q}(x,t))(v) dv,
		& \quad $x \in \p X,  t \in (t_n,t_{n+1})$,  \label{eq:Euler_limit_bc}
	\\
		\bm{q}(x,t_n^+) = \begin{cases}
		    \sum_{k=1}^{N_v}{\omega}_k \bm{e}_k  G_0(x,v_k),
		& n = 0, \\
        \sum_{k=1}^{N_v}{\omega}_k \bm{e}_k   \cE(\bm{q}(x,t))(v) 
        &  n > 1,
		\end{cases}\label{eq:Euler_limit_ic} 
		& \quad $x \in X$. 
	\end{numcases}
\end{subequations}
The system \eqref{eq:Euler_limit}	is a consistent discretization of the Euler equations.   However, 
\begin{enumerate}
	\item There is an error introduced at the beginning of each time step due to the quadrature error in velocity.  At time $t^{n+1}$ the moment $\bm{q}(x,t_n^-)$, which  is extracted  at the end of the interval $(t_n, t_{n+1})$ from the solution of \eqref{eq:Euler_limit_eqn}, is replaced by the reinitialization value $\bm{q}(x,t_n^+) = \sum_{k=1}^{N_v}{\omega}_k \bm{e}_k   \cE(\bm{q}_{\rm{c}}(x,t))(v)$.
	\item In practice (i.e. for finite $\epsilon$), $\bfu$ will decay to zero over a boundary layer of width $\epsilon$.  If this boundary layer is not well-resolved by the spatial and velocity discretization, the boundary condition \eqref{eq:Euler_limit_bc} will introduce errors into the overall solution.  Such issues are well-known in the transport \cite{bensoussan1979boundary} and gas kinetic theory literature \cite{bardos2006half} and are not specifically induced by the hybrid method.  In the numerical results of the current paper, solutions are computed on domains with equilibrium boundary conditions that are imposed far away from the interesting dynamics of the solution.   In such situations, these boundary issues will not pollute the solution.
\end{enumerate}

\subsection{Time Discretization} \label{sec:time discretization}
In  this section, we introduce two different time discretizations methods used to evolve the hybrid algorithm, formed by the components in  \eqref{eq:dvm} and \eqref{eq:Eulerc}, over a fixed time interval $(t_n,t_{n+1})$.   
\subsubsection{{The hybrid BERK2 time discretization }}\label{sec:BERK2_time_relabeling}

The base method combines two backward Euler steps for \eqref{eq:dvm} with a second-order predictor-corrector method for  \eqref{eq:Eulerc}.  This  Backward Euler Runge-Kutta (BERK2) method for a single time step is given in  Algorithm \ref{algorithm1}. The fully implicit treatment of \eqref{eq:dvm} removes the time step restrictions otherwise induced by the advection operator.  When combined with upwind DG discretizations in space, the resulting triangular linear system can be solved using sweep techniques that are especially easy to implement on regular Cartesian grids.  Moreover, because the discrete velocity components of $\bfu$  are independent, the implementation is easy to parallelize. The explicit predictor-corrector method for \eqref{eq:Eulerc} then induces a time step that scales inversely with $\Lambda$, the magnitude of the largest  eigenvalue in the Euler system.  Meanwhile, standard IMEX methods such as \cite{xiong} require a time step that scales inversely with $\vmax$.  
 
\begin{algorithm}[ht!]
\caption{Hybrid BERK2 temporal update from $t_n$ to $t_{n+1}$}
\label{algorithm1}
\textbf{Require}:  $G^-$ \hfill \Comment{Boundary data (used in (3) and (6)} )\\
\textbf{Require}:  $\{v_k\}_{k=1}^{N_v}$\hfill \Comment{Discrete velocities}\\
\textbf{Require}:  $\{\bm{e}_k\}_{k=1}^{N_v}$\hfill \Comment{Discrete collision invariants} \\
\textbf{Require}:  $\{w_k\}_{k=1}^{N_v}$\hfill \Comment{Quadrature weights } \\
\textbf{Input}: $\bm{g}^n$,  $\bm{q}_c^n (= 0)$  \hfill \Comment{Initial data}
\begin{algorithmic}[1]
\State $\bm{E} \leftarrow [w_1 \bm{e}_1 | \cdots | w_{N_v} \bm{e}_{N_v}]$
\State $\bfu[n] \leftarrow \bm{g}^n$
\State Solve $\dfrac{\bfu[n+1/2]-\bfu[n]}{\Delta{t}/2}
	+\bm{V}\uppartial_x \bfu[n+1/2]
	+ \dfrac{1}{\epsilon}{\bfu[n+1/2]}=0$ for ${\bfu[n+1/2]}$
\hfill \Comment{Uncollided backward Euler}
\State $\bqu[n+1/2]\leftarrow \bm E \bfu[n+1/2]$ \hfill\Comment{Uncollided moments}
\State Solve $\dfrac{{\bqc[n+\frac12]}-{\bqc[n]}}{\Delta{t}/2}
+\uppartial_x \bm{\Phi}(\bm{q}_{\rm{c}}^n)
=\dfrac{1}{\epsilon}{\bqu[n+\frac12]}$ for $\bqc[n+\frac12]$
\hfill \Comment{Collided moment  predictor}
\State Solve $\dfrac{\bfu[n+1/2]-\bfu[n]}{\Delta{t}}
+\bm{V}\uppartial_x \bfu[n+1]
=-\dfrac{1}{\epsilon}{\bfu[n+1]}$ for ${\bfu[n+1]}$
\hfill \Comment{Uncollided backward Euler}
\State $\bqu[n+1]\leftarrow \bm E \bfu[n+1]$ \hfill \Comment{Uncollided moments}
\State Solve $\dfrac{{\bqc[n+1]}-{\bqc[n]}}{\Delta{t}}+\uppartial_x \bm{\Phi}(\bm{q}_{\rm{c}}^{n+\frac{1}{2}})
=\dfrac{1}{\epsilon}{\bqu[n+1]}$ for $\bqc[n+1]$
\hfill \Comment{Collided moment corrector}
\For{$k=1,\dots,N_v$}
\State $\fck^{n+1}\leftarrow{\cE}(\bqc[n+1])(v_k)$ \hfill  \Comment{Collided reconstruction}
\EndFor
\State $\bm{g}^{n+1} = {\bfu}^{n+1} + {\bfc}^{n+1}$ \hfill   \Comment{Uncollided relabeling}
\State $\bqc[n+1]\leftarrow  0$						 \hfill   \Comment{Collided relabeling}
\end{algorithmic} 
\textbf{Output}: $\bm{g}^{n+1}$,  $\bm{q}_c^{n+1}$ 
\end{algorithm}

\subsubsection{The Euler limit}
\label{subsec:discrete-time-Euler-limit}
We investigate the formal limit of the BERK2 scheme when $\epsilon \to 0$.  Similar to the time-continuous setting in Section \ref{subsec:continous-time-Euler-limit}, we take moments of $f_{{\rm u}}$ in Lines 3 and 6 of Algorithm \ref{algorithm1} and add to Lines 5 and 8, respectively.  
Let $\bm{q}^{\alpha} = \bm{q}_u^{\alpha} +\bm{q}_c^{\alpha}$, for ${\alpha} \in \{n,n+1/2,n+1\}$.  Then

\begin{subequations}
	\label{eq:Eulerc_predict}
	\begin{numcases}{}
		\frac{\bm{q}^{n+1/2}(x) - \bm{q}^n(x)}{\Delta t /2}+  \sum_{k=1}^{N_v}{\omega}_k \bm{e}_k v_k \uppartial_x{f_{{\rm u},k}^{n+1/2}(x)} 
		+ \p_x\bm{\Phi}(\bm{q}_c^n(x)) = 0, & $\quad x \in X$, \label{eq:Eulerc_predict_eqn}\\
		\bm{\Phi}(\bm{q}^n(x))=   \sum_{\{k:  n(x) v_k < 0\}}{\omega}_k v_k \bm{e}_k  G^-(x,v_k) 
		+ \int_{\{v:n(x) v > 0 \}}  v \bm{e} \cE(\bqc[n](x))(v)dv,
		& \quad $x \in \p X$,  \label{eq:Eulerc_predict_bc}
	\end{numcases} 
\end{subequations}

\begin{subequations}
	\label{eq:Eulerc_correct}
	\begin{numcases}{}
	\frac{\bm{q}^{n+1}(x) - \bm{q}^{n}(x)}{\Delta t }+  \sum_{k=1}^{N_v}{\omega}_k \bm{e}_k v_k \uppartial_x{f_{{\rm u},k}^{n+1}(x)} 
	+ \p_x\bm{\Phi}(\bm{q}_c^{n+1/2}(x)) = 0 , & $\quad x \in X$,  \label{eq:Eulerc_correct_eqn}\\
		\bm{\Phi}(\bm{q}^n(x))=   \sum_{\{k:  n(x) v_k < 0\}}{\omega}_k v_k \bm{e}_k  G^-(x,v_k)  
		+ \int_{\{v:n(x) v > 0 \}}   v \bm{e} \cE(\bqc[n+1/2](x)) (v)dv,
		& \quad $x \in \p X$.  \label{eq:Eulerc_correct_bc}
	\end{numcases} 
\end{subequations}
The update equations for $\bfu$ in Lines 3 and 6 of Algorithm \ref{algorithm1} take the form
\begin{equation}
	\label{eq:fu-ss}
	\frac{f_{{\rm u},k}^{{\alpha}}(x)- f_{{\rm u},k}^{n}(x)}{\tau} + v_k \uppartial_x f_{{\rm u},k}^*(x)  = -\frac{1}{\epsilon}f_{{\rm u},k}^*(x),
\end{equation}
for some $\tau>0$ and ${\alpha} \in \{ n+1/2,n +1\}$.  Let 
\begin{equation}
	\label{eq:s-star}
	s^*:= s^*(x,v_k) = \min_{s >0} \{ s: x - v_k s  \in \gamma^-_k\}
	\quad \text{and} \quad
	x_*:= x^*(x,v_k) = x - v_k s_*.
\end{equation}
Then the solution of \eqref{eq:fu-ss} is
\begin{equation}
\label{eq:fuk-ss}
	f_{{\rm u},k}^{\alpha}(x) = \begin{cases}
		e^{-\frac{\sigma(x - x_*)}{v_k}} G^-(x_*,v_k,t_\alpha) +  \frac{1}{\tau v_k}   \int_{x_*}^{x} 	e^{-\frac{\sigma(x - \xi)}{v_k}} f_{{\rm u},k}^{n}(\xi) d \xi, &v_k  \ne 0, \\
		\frac{\epsilon}{\tau + \epsilon} f_{{\rm u},k}^{n}(x), &v_k =0,  \\
	\end{cases}
\end{equation}
where $\sigma =\frac{1}{\tau} + \frac{1}{\epsilon}$.  Assume that  $\sup_\xi | f_{{\rm u},k}^{n}(\xi) | \leq M < \infty$.  Then
\begin{equation}
	\begin{split}
	\frac{1}{\tau v_k}   \int_{x_*}^{x} 	e^{-\frac{\sigma(x - \xi)}{v_k}} f_{{\rm u},k}^{n}(\xi) d \xi 
	& \leq 	\frac{M}{\tau v_k}  \int_{x_*}^{x} 	e^{-\frac{\sigma(x - \xi)}{v_k}} d \xi 
	= \frac{M}{\tau \sigma}	\left(1 - e^{-\frac{\sigma(x -x_*)}{v_k}} \right) \leq  \frac{M \epsilon}{1+{\epsilon}},
	\end{split}
\end{equation}
whereas the other terms in \eqref{eq:fuk-ss} clearly decay as $\epsilon \to 0$.  Thus if $f_{{\rm u},k}^{n}$ is bounded in $\xi$, it follows that $f_{{\rm u},k}^{\alpha}(x)  \to 0$ and $\ \bm{q}_c^{\alpha} \to \bm{q}^{\alpha} $ as $\epsilon \to 0$ for any $x \in X$.  As a result, in the limit $\epsilon \to 0$, \eqref{eq:Eulerc_predict} and \eqref{eq:Euler_limit_correct} transition to a second-order predictor-corrector method for the Euler equations:
\begin{subequations}
	\label{eq:Euler_limit_predict}
	\begin{numcases}{}
		\frac{\bm{q}^{n+1/2} - \bm{q}^n}{\Delta t /2}
		+ \p_x\bm{\Phi}(\bm{q}^n) = 0, & $\quad x \in X$, \label{eq:Euler_limit_predict_eqn}\\
		\bm{\Phi}(\bm{q}^n(x))=   \sum_{\{k:  n(x) v_k < 0\}}{\omega}_k v_k \bm{e}_k  G^-(x,v_k,t^{n+1/2}) 
		+ \int_{\{v:n(x) v > 0 \}}  v \bm{e} \cE(\bm{q}^n(x))(v)dv,
		& \quad $x \in \p X$,  \label{eq:Euler_limit_predict_bc}
	\end{numcases} 
\end{subequations}

\begin{subequations}
	\label{eq:Euler_limit_correct}
	\begin{numcases}{}
		\frac{\bm{q}^{n+1} - \bm{q}^{n}}{\Delta t }
		+ \p_x\bm{\Phi}(\bm{q}^{n+1/2}) = 0 , & $\quad x \in X$,  \label{eq:Euler_limit_correct_eqn}\\
		\bm{\Phi}(\bm{q}^n(x))=   \sum_{\{k:  n(x) v_k < 0\}}{\omega}_k v_k \bm{e}_k   G^-(x,v_k,t^{n+1}) 
		+ \int_{\{v:n(x) v > 0 \}}  v \bm{e} \cE(\bm{q}^{n+1/2}(x)) (v)dv,
		& \quad $x \in \p X$.  \label{eq:Euler_limit_correct_bc}
	\end{numcases} 
\end{subequations}

\subsubsection{{Correction and conservation fix}}\label{sec:BERK2_time_relabeling_correction_conservation_Fix}

A modified version of the BERK2 time discretization is given in Algorithm \ref{algorithm1}.  This modified version includes (i) a conservation fix to fix the deficiency caused by consistency errors in the discrete velocity quadrature and (ii) a relabeling step that solves the original BGK equation using a Maxwellian constructed from the hybrid strategy.  The details of these modifications are given in Algorithm \ref{algorithm2}.

Conservation of mass, momentum, and energy in the BGK equation is a consequence of \eqref{eq:collision_invariants}. 
To correct consistency errors that violate \eqref{eq:collision_invariants}, we adopt the technique first introduced in \cite{gamba2009spectral}.  
Given a discrete  velocity vector $\bm{f} = [f_1 ,\dots , f_{N_v}]^\Tr \in 
\bbR^{N_v},$ let  $ \bm{E}  = [w_1 \bm{e}_1 |  \cdots  | w_{N_v}\bm{e}_{N_v} ]$ and let
\begin{equation}
\bm{ q_{\bm{f}} } = \sum_{k=1}^{N_v} w_k \bm{e}_k f_k =  \bm{E}  \bm{f} \in \bbR^{3}
\end{equation}
be the moments associated to $\bm{f}$.  In order to match a moment $\widetilde{\bm{q}} \in \bbR^{3}$ where $\widetilde{\bm{q}} \ne \bm{q}_f$, the corrected discrete velocity vector $\widetilde{\bm{f}}$ is given by
\begin{equation}
	\widetilde{\bm{f}} = \operatorname*{argmin}_{\bm{g} \in 
		\bbR^{N_v} }
	\{ \| \bm{f} - \bm{g}\|^2_2 : \bm{E} \bm{g} =\widetilde{\bm{q}} \}.
\end{equation} 
The solution to this optimization problem satisfies $\bm{q}_{\widetilde{\bm{f}}} \equiv \bm{E} \bm{f} = \widetilde{\bm{q}}$ and is given by
\begin{align}
	\widetilde{\bm{f}}&=\bm{f}+\bm{E}^{\Tr}(\bm{E}\bm{E}^{\Tr})^{-1}(\widetilde{\bm{q}}-\bm{E} \bm{f}).
\end{align}

The relabeling step constructs a discrete Mawellian from the moments $\bm{q}^{n+1} = \bqu[n+1] + \bqc[n+1]$.  A final implicit update of the BGK equation uses this discrete Maxwellian as a fixed source.  Because there is a fixed source, this update can be efficiently completed without iteration. The implicit update has been done through the backward differentiation formula (BDF2)  discretization.  This approach to relabeling was introduced in \cite{crockatt2020improvements} and tends to yield more accurate answers because errors in the moments that are induced by the hybrid approximation are smaller than the errors in the kinetic distribution.

\begin{algorithm}[ht!]
\caption{Hybrid BERK2 temporal update from $t_n$ to $t_{n+1}$  with \rfbox{BDF2 correction} and \sfbox {conservation fix}}
	\label{algorithm2}
\textbf{Require}:  $G^-$ \hfill \Comment{Boundary data (used in (3) and (6)} )\\
\textbf{Require}:  $\{v_k\}_{k=1}^{N_v}$\hfill \Comment{Discrete velocities}\\
\textbf{Require}:  $\{\bm{e}_k\}_{k=1}^{N_v}$\hfill \Comment{Discrete collision invariants} \\
\textbf{Require}:  $\{w_k\}_{k=1}^{N_v}$\hfill \Comment{Quadrature weights } \\
\textbf{Input}: $\bm{g}^n$,  $\bm{q}_c^n (= 0)$  \hfill \Comment{Initial data}
	\begin{algorithmic}[1]
	\State $\bm{E} \leftarrow [w_1 \bm{e}_1 | \cdots | w_{N_v} \bm{e}_{N_v}]$
	\State $\bfu[n] \leftarrow \bm{g}^n$
	\State Solve $\displaystyle\frac{\bfu[n+1/2]-\bfu[n]}{\Delta{t}/2}
	+\bm{V}\uppartial_x \bfu[n+1/2]
	=-\frac{1}{\epsilon}{\bfu[n+1/2]}$ for ${\bfu[n+1/2]}$ \hfill \Comment{Uncollided backward Euler}
		\State $\bqu[n+\frac12]\leftarrow \bm{E} \bfu$ \hfill\Comment{Uncollided moments}
		\State Solve $\dfrac{{\bqc[n+\frac12]}-{\bqc[n]}}{\Delta{t}/2}
		+\uppartial_x \bm{\Phi}(\bm{q}_{\rm{c}}^n)
		=\frac{1}{\epsilon}{\bqu[n+\frac12]}$ for $\bqc[n+\frac12]$
		\State $\bfc[n+\frac12] \leftarrow {\cE}(\bm{q}_{\rm{c}}^{n+\frac12})$	\hfill \Comment{Discrete collided distribution}
  \State $\bm{g}^{n+\frac12} \leftarrow \bfu[n+\frac12]+\bfc[n+\frac12]$\hfill \Comment{midstep for BDF2 correction}
		\State Solve $\dfrac{\bfu[n+1]-\bfu[n]}{\Delta{t}}
	+\bm{V}\uppartial_x \bfu[n+1]
	=-\dfrac{1}{\epsilon}{\bfu[n+1]}$ for ${\bfu[n+1]}$
	\hfill \Comment{Uncollided backward Euler}
				\State $\bqu[n+1]\leftarrow \bm{E} \bfu$ \hfill\Comment{Uncollided moments}
		\State Solve $\dfrac{{\bqc[n+1]}-{\bqc[n]}}{\Delta{t}}+\uppartial_x \bm{\Phi}(\bm{q}_{\rm{c}}^{n+\frac{1}{2}})
		=\dfrac{1}{\epsilon}{\bqu[n+1]}$ for $\bqc[n+1]$
		\hfill \Comment{Collided moment corrector}
		{\State {$\bm{q}^{n+1}=\bm{q}_{\rm{u}}^{n+1}+\bm{q}_{\rm{c}}^{n+1}$}} \hfill \Comment{Full moments}
		\For{$k=1,\dots  ,N_v$}
			\State $M^{n+1}_k \leftarrow {\cE}(\bm{q}^{n+1})(v_k)$	\hfill \Comment{Discrete Maxwellian}
		\EndFor 
		\State $\bm{M}^{n+1} \leftarrow [M^{n+1}_1, \dots, M^{n+1}_{N_v}]^{\Tr} $
		{\State \sfbox{${\bm{M}^{n+1}}\leftarrow \bm{M}^{n+1}+\bm{E}^{\Tr}(\bm{E}\bm{E}^{\Tr})^{-1}(\bm{q}^{n+1}-\bm{E}\bm{M}^{n+1})$}\hfill \Comment{Conservation fix}}
\State \rfbox{{Solve
				$\dfrac{{\bm{f}^{n+1}}-\frac{4}{3}{\bm{g}^{n+\frac{1}{2}}}+\frac{1}{3}{\bm{g}^n}}{\frac{\Delta{t}}{3}}+\bm{V}\uppartial_x \bm{f}^{n+1} + \dfrac{1}{\epsilon} \bm{f}^{n+1} =  \dfrac{1}{\epsilon} \bm{M}^{n+1} $ for $\bm{f}^{n+1}$}}
		\hfill \Comment{BDF2 correction}  
		{\State \sfbox{$\bm{f}^{n+1} \leftarrow \bm{f}^{n+1}+\bm{E}^{\Tr}(\bm{E}\bm{E}^{\Tr})^{-1}(\bm{q}^{n+1}-\bm{E}\bm{f}^{n+1})$}
			\hfill \Comment{Conservation fix}}
		\State $\bm{g}^{n+1} \leftarrow   \bm{f}^{n+1}$	
		\State $\bm{q}_c^{n+1} \leftarrow 0$ 
	\end{algorithmic} 
	\textbf{Output}: $\bm{g}^{n+1}$,  $\bm{q}_c^{n+1} $
\end{algorithm}

%% file: sec3_implementation_details.tex
\subsection{Spatial discretization}\label{sec:spatial_discretization}

Both \eqref{eq:dvm} and \eqref{eq:collided_moments1} are discretized with a discontinuous Galerkin method \cite{reed1973triangular, shu2009discontinuous}.  Since this method is by now well-known, the presentation here will be naturally brief. For simplicity, we discretize the time continuous equations.

Given $N_x + 1$ grid points $\xl=x_{1/2} < x_{3/2} < \dots < x_{N_x +1/2} = \xr$, let $S_i = (x_{i-1/2},x_{i+1/2})$ be a cell of width $h_i = x_{i+1/2} - x_{i-1/2}$ and center $\frac12(x_{i+1/2} + x_{i-1/2})$.  Let ${\bbP}_{N}$ be the set
of all polynomials with maximum polynomial degree $N$, and let $h = \max_i h_i$. Define the broken finite element space
\begin{equation}
\label{eqn:broken_space1}
    \WS^{h}_N := \left\{ w^{h} \in L^{2}\left[ \xl, \xr \right] : \,
    w^{h} \bigl|_{S_i} \in  {\bbP}^{N} \quad \forall i \in \{1, \ldots,N_x \} \right\} .
\end{equation}
  A basis for $\WS^{h}_N $ is formed by the functions
\begin{equation}
    w_{i,\ell}(x) = p_{\ell}\left( \frac{x - x_i}{h_i / 2}\right), 
    \qquad i \in \{1 , \dots , N_x\},
    \quad \ell \in \{1 , \dots , N+1\},
\end{equation}
where $\{p_\ell \}_{\ell=1}^{N+1}$ is a set of degree $N$ Lagrange polynomials on $[-1,1]$.  These polynomials are constructed from the Gauss-Lobatto points $\{\xi_{\ell}\}_{\ell=1}^{N+1}$ on $[-1,1]$:
\begin{equation}
 p_\ell\left(\xi\right) = \prod_{\underset{\ell' \ne \ell}{\ell'=1}}^{N+1} \frac{\left( \xi - \xi_{\ell'} \right)}{\left(\xi_\ell- \xi_{\ell'} \right)} ,\qquad \ell \in \{1 ,\dots, N+1\}.
\end{equation}

\subsubsection{Uncollided equation}\label{sec:2.4}
To discretize \eqref{eq:dvm}, we seek for each $k \in \{1, \dots, N_v\}$, a function 
\begin{equation}
f^h_{{\rm{u}},k}(x,t) = \sum_{i=1}^{N_x} \sum_{\ell=1}^{N+1} f_{k,i,\ell}(t) w_{i,\ell}(x),
\end{equation}
such that
for all $i \in \{1 , \dots , N_x\}$ and $\ell \in \{1 , \dots , N+1\}$,
\begin{equation}\label{eq:semi-DG}
\begin{aligned}
   \uppartial_t\int_{S_i}{{f^h_{{\rm{u}},k}(x,t)}w_{i,\ell}}\,dx
    &+v_k\left({f}^{h,*}_{{\rm{u}},k}(x_{i+1/2},t)w_{i,\ell}(x^{-}_{i+1/2})
            -{f}^{h,*}_{{\rm{u}},k}(x_{i-1/2},t)w_{i,\ell}(x^{+}_{i-1/2})\right)\\
    &-v_k\int_{S_i}{{f^h_{{\rm{u}},k}(x,t)}w_{i,\ell}'(x)}\,dx
    +\frac{1}{\varepsilon}\int_{S_i}{f^h_{{\rm{u}},k}(x,t)w_{i,\ell}(x)}\,dx =0,
\end{aligned}
\end{equation}
where ${f}^{h,*}_{{\rm{u}},k}$ is the upwind trace:
\begin{subequations}
\begin{numcases}{{f}^{h,*}_{{\rm{u}},k}(x_{i+1/2},t)=}
{f}^{h}_{{\rm{u}},k}(x^-_{i+1/2},t), & $v_k > 0$, \\
{f}^{h}_{{\rm{u}},k}(x^+_{i+1/2},t), & $v_k < 0$.
\end{numcases}
\end{subequations}
In terms of $\bm{f}_{k,i} = (f_{k,i,1}, \dots, f_{k,i,N+1})^\Tr$, \eqref{eq:semi-DG} takes the form, for $i \in \{1,\dots,N_x$\}, 
\begin{subequations}
\label{eq:uncollided_dg_simple}
\begin{numcases}{}
\label{eq:uncollided_dg1_simple}
\frac{h_i }{2} \bm{J} \, \frac{d \bm{f}_{k,i}}{dt}(t)
    +v_k\left[\bm{L}_1\,{\bm{f}_{k,i}(t)}
    -\bm{L}_2\,{\bm{f}_{k,i-1}(t)}
    -\bm{K}\,{\bm{f}_{k,i}(t)}\right]
    +\frac{h_i}{2\varepsilon}\bm{J}_i\,{\bm{f}_{k,i}(t)} =0, & $\nu_k > 0,$ \\
    \label{eq:uncollided_dg2_simple}
\frac{h_i}{2}      \bm{J}\,\frac{d \bm{f}_{k,i}}{dt}(t)
    +v_k\left[\bm{L}_3\,{\bm{f}_{k,i+1}(t)}
    -\bm{L}_4\,{\bm{f}_{k,i}(t)}
    -\bm{K}\,{\bm{f}_{k,i}(t)}\right]
    +\frac{h_i}{2\varepsilon}\bm{J}\,{\bm{f}_{k,i}(t)} =0, & $\nu_k < 0$,
\end{numcases}
\end{subequations}
where $\bm{J}\in R^{(N+1)\times (N+1)}$ and $\bm{K}\in R^{(N+1)\times (N+1)}$ have components
\begin{align}
\bm{J}_{\ell,\ell'} = \int_{-1}^{1}\Phi_{\ell}(\xi) \, \Phi_{\ell'}(\xi) \, \diff \xi \qquad \text{and} \qquad
\bm{K}_{\ell,\ell'} = \int_{-1}^{1}\Phi_{\ell}(\xi) \, \Phi'_{\ell'}(\xi) \, \diff \xi,
\end{align}
and the matrices $\bm{L}_r$, $r \in \{1,2,3,4\}$, have components
\begin{subequations}
\begin{gather}
    (\bm{L}_1)_{\ell,\ell'}
        =\Phi_{\ell}(1)\Phi_{\ell'}(1) 
        = \delta_{\ell,N} \delta_{\ell',N},
    \qquad
    (\bm{L}_2)_{\ell,\ell'}
        =\Phi_{\ell}(-1)\Phi_{\ell'}(1) 
        = \delta_{\ell,1} \delta_{\ell',N} ,
    \\
    (\bm{L}_3)_{\ell,\ell'}
        =\Phi_{\ell}(1)\Phi_{\ell'}(-1)
        = \delta_{\ell,N} \delta_{\ell',1} ,
    \qquad
    (\bm{L}_4)_{\ell,\ell'}
        =\Phi_{\ell}(-1)\Phi_{\ell'}(-1)
        = \delta_{\ell,1} \delta_{\ell',1} ,
\end{gather}
\end{subequations}
with $\delta$ the usual Kronecker delta.  

The implementation of \eqref{eq:uncollided_dg_simple} requires the specification of ghost values $\bm{f}_{k,0}$ when $\nu_k>0$ and $\bm{f}_{k,N_x+1}$ when $\nu_k<0$.  In the numerical examples of Section \ref{sec:results}, either periodic or Dirichlet boundary conditions are used.  For the periodic case, we set
\begin{numcases}{}
    \bm{f}_{k,0} (t)   = \bm{f}_{k,0}, & $\nu_k > 0$,\\
    \bm{f}_{k,N_x+1}(t)  = \bm{f}_{k,0}(t), & $\nu_k < 0$.
\end{numcases}
For the Dirichlet case, we set the nodal values equal to the boundary condition on each end:
\begin{numcases}{}
    f_{k,0,\ell}(t)  = G^-(\xl,\nu_k,t), & $\nu_k > 0, \quad \ell=1,\dots, N+1$,\\
    f_{k,N_x+1,\ell}(t)  = G^-(\xr,\nu_k,t), & $\nu_k < 0, \quad \ell=1,\dots, N+1$.
\end{numcases}
In practice, the equations in \eqref{eq:uncollided_dg_simple} are solved using backward Euler or BDF2.  In either case, an effective steady-state problem is created which can be solved efficiently using sweeps that update $\bm{f}_{k,i}$ using the information from $\bm{f}_{k,i-1}$ when $\nu_k >0$ and information from 
$\bm{f}_{k,i+1}$ when $\nu_k<0$.  The sweeps are  independent across angles and require only the inversion of an $(N+1) \times (N+1)$ per angle to update the degrees of freedom in each spatial cell. For problems with incoming boundary data given, only one sweep across the spatial mesh is needed for each angle.  For problems with periodic or reflecting boundaries, several sweep iterations are needed to reach convergence. For the numerical examples in Section \ref{sec:results} with periodic boundary conditions, we update the ghost cells after each sweep and observe that only one or two iterations are required for convergence.

\subsubsection{Collided equation}\label{sec:DG_collided}
To discretize \eqref{eq:collided_moments1}, we seek for each a vector-valued function 
\begin{equation}
	\bm{q}^h_{\rm{c}} = \sum_{i=1}^{N_x} \sum_{\ell=1}^{N+1}  \bm{q}_{i,\ell}(t) w_{i,\ell}(x),
\end{equation}
 such that
for all $i \in \{1 , \dots , N_x\}$ and $\ell \in \{1 , \dots , N+1\}$,
\begin{equation}
    \begin{split}
      \uppartial_t\int_{S_i} \bm{q}^h_{\rm{c}}(x,t) w_{i,\ell}(x) dx 
      &+ \bm{\Phi}^*(\bm{q}^h_{\rm{c}}(x^-_{i+1/2},t),(\bm{q}^h_{\rm{c}}(x^+_{i+1/2},t) )w_{i,\ell}(x^-_{i+1/2}) \\
     & - \bm{\Phi}^*(\bm{q}^h_{\rm{c}}(x^-_{i-1/2},t) ,\bm{q}^h_{\rm{c}}(x^+_{i-1/2},t)  )w_{i,\ell}(x^+_{i-1/2}) \\
      &+\int_{S_i} \bm{\Phi}(\bm{q}^h_{\rm{c}}(x,t))w'_{i,\ell}(x)\,dx 
      = \frac{1}{\varepsilon}\int_{S_i}{\bm{q}^h_{\rm{u}}(x,t) w_{i,\ell}}(x)\,dx ,
    \end{split}
\end{equation}
where
\begin{equation}
    \bm{q}^h_{\rm{u}}(x,t)
    = \sum_{k=1}^{N_v}{\omega_k} \bm{e}_k f^h_{\rm{u},k}(x,t),
\end{equation}
and the numerical flux $\bm{\Phi}^*$ is the local Lax-Friedrichs flux:
\begin{equation}
    \bm{\Phi}^*(\bm{q})=\bm{\Phi}^*(\bm{q}^{+},\bm{q}^{-})=\frac{1}{2}\left({\bm{\Phi}(\bm{q}^{+})+\bm{\Phi}(\bm{q}^{-})-\max_{\bm{q}\in\{\bm{q}^{-},\bm{q}^{+}\}}{|\lambda(\bm{\Phi}'(\bm{q}))|(\bm{q}^{+}-\bm{q}^{-})}}\right).
\end{equation}
The flux Jacobian matrix
\begin{equation}
\bm{\Phi}'(\bm{q})
=\begin{pmatrix}
0&1&0\\
\frac{(\gamma-3)}{2}\frac{q_2^2}{q_1^2}&(3-\gamma)\frac{q_2}{q_1}&(\gamma-1)\\
-\gamma\frac{q_2 q_3}{{q_1}^2}+(\gamma-1)\frac{{q_2}^3}{{q_1}^3}&\gamma\frac{q_3}{q_1}-\frac{3}{2}(\gamma-1)\frac{{q_2}^2}{{q_1}^2}&\gamma\frac{q_2}{q_1}
\end{pmatrix}
=\begin{pmatrix}
0&1&0\\
0&0&2\\
-3\frac{q_2 q_3}{{q_1}^2}+2\frac{{q_2}^3}{{q_1}^3}&3\frac{q_3}{q_1}-3\frac{{q_2}^2}{{q_1}^2}&3\frac{q_2}{q_1}
\end{pmatrix}
\end{equation}
has eigenvalues $
    \lambda
    =\{u,u-c,u+c\}
$,
with
\begin{equation}
    u = \frac{q_2}{q_1} 
    \quad \text{and} \quad
    c = \sqrt{\frac{\gamma(2q_1 q_3 -{q_2}^2)}{{q_1}^2}} = \sqrt{\gamma\theta}=\sqrt{3\theta}.
\end{equation}
Thus
\begin{equation}\label{eq:lambda}
	 \Lambda=\max_{\bm{q}\in\{\bm{q}^{-},\bm{q}^{+}\}} | \lambda (\bm{\Phi}'(\bm{q}))|
	 =\max_{\bm{q}\in\{\bm{q}^{-},\bm{q}^{+}\}}\big(\max{\left(|u|,|u-c|,|u+c|\right)}\big). 
\end{equation}
Once $\bm{q}^h_{\rm{c}}$ is computed, the associated Maxwellian  $\cE(\bm{q}^h_{\rm{c}})$ can be computed at each DG node.

In order to prevent spurious oscillations, we post-process the moments after each update $\bm{q}^h_{\rm{c}}$ using the total variation bounded (TVB) slope limiter introduced in \cite{cockburn2}.  This limiter requires the specification of a TVB parameter $M$ which is chosen to be $M=20$.  However, it is only applied to certain test problems; see details in the next section.

%% file: sec4_results.tex
\section{Numerical Results}\label{sec:results}
In this section, we present numerical results.   For all tests, we compare the BERK2 hybrid methods to the solutions obtained using the IMEX-DG scheme.  Both implementations utilize the same discrete velocity quadratures for velocity discretization.  Additionally, the IMEX-DG scheme relies on the stiffly accurate IMEX-SSP2(3,2,2) scheme \cite{pareschi2} for 2nd order and the IMEX-ARS(4,4,3) scheme \cite{ascher1997} for the 3rd-order in the temporal domain.  See \cite{pareschi2,ascher1997} for more details about these IMEX schemes. Note that DG$k$ refers to the $k$-th order DG scheme, contrasting with the conventional $(k+1)$-th order scheme.   

The time step for BERK2-DG is $\Delta{t}=C \Lambda^{-1} \Delta{x}$, where 
\begin{equation}
	\Lambda \equiv \max_{0 \leq i \leq N_x} \max_{ \bm{q}\in\{\bm{q}(x^-_{i+1/2}),\bm{q}(x^+_{i+1/2}) \} } 
	| \lambda (\bm{\Phi}' (\bm{q}) ) |
\end{equation}
approximates the maximum wave speed associated to the flux matrix ${\Phi'(\bm{q})}$ for \eqref{eq:Eulerc_eqn} over the spatial domain $X$. Meanwhile, the time step for the IMEX-DG is  $\Delta{t}=C \vmax^{-1} \Delta{x}$, where $\vmax$ is chosen to ensure that the amount of mass lost in the  distribution function is within acceptable bounds. For the tests in Sections \ref{subsec:asymptotics} and \ref{subsec:accuracy}, the CFL constant $C$ for the BERK2-DG4 hybrid is given in \tableref{tab:cfl_4.1,4.2}. In the accuracy test case (Section \ref{subsec:accuracy}), the constants are chosen to be smaller to postpone saturation in convergence due to temporal error. In Sections \ref{subsec:sod}-\ref{subsec:gas_injection} a DG3 scheme is used for all spatial discretizations and the CFL constant $C$ is given in \tableref{tab:comparisons_dt}.

DG slope limiters are not used for the tests in Sections \ref{subsec:asymptotics} and \ref{subsec:accuracy}, which focus on smooth solutions.  In the remaining tests, the TVB slope limiter  \cite{cockburn2} is used in the simulation of the collided moments $\bqc$ and for simulations of the Euler equations. 

In all the numerical tests, the initial distribution function $f$ is given by a Maxwellian, as defined in \eqref{eq:Maxwellian}.  That is, $f(x,v,t=0)=\cE(\bm{q}(x,t=0))(v)$, where $\bm{q}(x,t=0)$ is determined from the variables $\rho(x, t=0)$, $u(x, t=0)$, and $p(x, t=0)$.

Reference solutions for the Euler equations are computed with the DoGPack software package \cite{dogpack} using a 4th-order DG scheme with $N_x=3000$ spatial points and a 4th-order low-storage strong stability Runge-Kutta method \cite{ketcheson} for time integration.   Reference solutions for the BGK equations are computed using the 3rd-order IMEX3 scheme with refined grid points. To be more specific, the reference solution named IMEX3H is used in Sections \ref{subsec:sod}-\ref{subsec:shu-osher} and in Section \ref{subsec:gas_injection}. This is a highly resolved IMEX3-DG3 solution with a grid size of $N_x=N_v=1000$ in Sections \ref{subsec:sod}-\ref{subsec:shu-osher} and $N_x=800$ and $N_v=1100$ in Section \ref{subsec:gas_injection}. 

The purpose of each numerical test is as follows:  The asymptotic test in  Section \ref{subsec:asymptotics} is to assess convergence to the Euler limit as $\epsilon \to 0$.  The accuracy test in Section \ref{subsec:accuracy} is to investigate the accuracy of the BERK2-DG method for different values of $\epsilon$.  In particular, it is expected that high-resolution, high-order schemes will experience temporal accuracy saturation; understanding where this saturation occurs is important for practical applications.  The tests in Sections \ref{subsec:sod}-\ref{subsec:shu-osher} are standard benchmarks from gas dynamics, meant to assess how the BERK2-DG scheme captures shocks, contacts, and rarefactions.   Finally, the test in Section  \ref{subsec:gas_injection} is designed to emphasize the main advantage of the hybrid BERK2 scheme which is the efficiency afforded by the relaxed time step restriction when $\Lambda < \vmax$.

\begin{table}[hb!]
\caption{CFL constant $C$, used in sections \ref{subsec:asymptotics} and \ref{subsec:accuracy}.}
\centering
\begin{tabular}{|l||l|l|l|l|l|}
\hline
section& CFL constant & BERK2-DG2 & BERK2-DG3 & BERK2-DG4 \\ 
\hline\hline
\ref{subsec:asymptotics}& $C$  & - & - & 0.1\\ 
\hline
\ref{subsec:accuracy}& $C$  & 0.2 & 0.1 & 0.05\\ 
\hline
\end{tabular}
\label{tab:cfl_4.1,4.2}
\end{table}

\begin{table}[ht!]
\caption{Comparisons of the number of time steps, CFL constant $C$, maximum velocity $\vmax$, and the wave speed $\Lambda$ in each of the benchmark problems from Sections \ref{subsec:sod}-\ref{subsec:gas_injection}. The IMEX3L is an IMEX-ARS(4,4,3), IMEX2L is an IMEX-SSP2(3,2,2), BERK2L is a BERK2 solution, and BERK2LS is BERK2 solutions with reduced CFL constants. Since the wave-speed $\Lambda$ changes dynamically in each time-step, the time-step size for the BERK2 methods also changes, especially for the gas-injection problem. Note that the ratio $\vmax/\Lambda$ indicates the efficiency of BERK2 scheme in general.}
\centering
\begin{tabular}{|l||l|rrr|l|c|c|c|}
\hline
section              & time-stepping & \multicolumn{3}{c|}{$\#$ of time steps}               & $C$& $\vmax$  & $\Lambda$  & $\vmax/\Lambda$ \\ \cline{3-5}
                     &               & \multicolumn{1}{r}{$\epsilon=1$}&\multicolumn{1}{r}{$10^{-2}$}&\multicolumn{1}{r|}{$10^{-6}$} &  &          &            &                  \\ \hline\hline 
\multirow{3}{*}{\begin{tabular}[c]{@{}l@{}}\ref{subsec:sod}\\ (Sod)\end{tabular}} & IMEX3L & 429& 429& 429 & 0.14 & \multirow{3}{*}{6.0} & \multirow{3}{*}{1.732-2.934} & \multirow{3}{*}{2.045-3.464}\\ 
                     & IMEX2L & 300& 300& 300   & 0.2 &                 &          &         \\ 
                     & BERK2L & 134& 137& 138   & 0.2 &                 &          &         \\ 
                     & BERK2LS & 274& 274& 274   & 0.1 &                 &          &         \\ \hline\hline
\multirow{3}{*}{\begin{tabular}[c]{@{}l@{}}\ref{subsec:lax}\\ (Lax)\end{tabular}} & IMEX3L & 536& 536& 536 & 0.14 & \multirow{3}{*}{15.0} & \multirow{3}{*}{5.575-8.110} & \multirow{3}{*}{1.850-2.691}\\ 
                     & IMEX2L & 375& 375& 375   & 0.2 &                 &          &         \\ 
                     & BERK2L & 140& 140& 189   & 0.2 &                 &          &         \\ 
                     & BERK2LS & 280& 280& 280   & 0.1 &                 &          &         \\ \hline\hline
\multirow{3}{*}{\begin{tabular}[c]{@{}l@{}}\ref{subsec:shu-osher}\\ (Shu-Osher)\end{tabular}} & IMEX3L & 1800& 1800& 1800 & 0.14 & \multirow{3}{*}{14.0} & \multirow{3}{*}{6.206-8.085} & \multirow{3}{*}{1.732-2.256}\\ 
                     & IMEX2L & 1260& 1260& 1260   & 0.2 &                 &          &         \\ 
                     & BERK2L & 559& 583& 686   & 0.2 &                 &          &         \\ 
                     & BERK2LS & 1148& 1148& 1148   & 0.1 &                 &          &         \\ \hline\hline
\multirow{3}{*}{\begin{tabular}[c]{@{}l@{}}\ref{subsec:gas_injection}\\ (gas-injection)\end{tabular}} & IMEX3L & 1000& 1000& 1000 & 0.1 & \multirow{3}{*}{110.0} & \multirow{3}{*}{0.548-8.549} & \multirow{3}{*}{12.867-200.730}\\ 
                     & IMEX2L & 1000& 1000& 1000   & 0.1 &                 &          &         \\ 
                     & BERK2L & 14& 27& 49   & 0.1 &                 &          &       \\ 
                     & BERK2LS & 61& 117 & 210   & 0.025 &                 &          &       \\ \hline
\end{tabular}
\label{tab:comparisons_dt}
\end{table}
\subsection{Asymptotic test}\label{subsec:asymptotics}
In this test, we confirm numerically the Euler limit established in Sections \ref{subsec:continous-time-Euler-limit} and \ref{subsec:discrete-time-Euler-limit}.  The domain is $X \times V=[-\pi,\pi] \times[-7,7]$ and the boundary conditions are periodic.  The initial data is $G^0(x,v) = \cE(\bm{q}(x,0))(v)$, where $\bm{q}$ is determined from the variables
\begin{equation}
\label{eq:asymp_test_ic}	\rho(x,t=0)=1+0.2\sin10x, \qquad u(x,t=0) = 1,  \qquad p(x,t=0) = 1,
\end{equation} 
and $\epsilon=10^{-12}$. 
Here, the Euler solution is used as the reference solution with the initial condition $\rho_0(x)=1+0.2\sin10x$.
The results for the mass density at time  $t=0.1$ are given in  \tableref{tab:asymptotic_results_BERK2_all}. All three versions---standard BERK2, BERK2 with BDF2 correction step, and BERK2 with BDF2 correction and the conservation fix---exhibit the 2nd-order temporal convergence rate that is expected based on the limiting equations in \eqref{eq:Euler_limit_predict} and \eqref{eq:Euler_limit_correct}. 
\begin{table}[ht!]
\caption{Asymptotic test from Section \ref{subsec:asymptotics}.  The $L_2$ errors and order of accuracy for the density at $t=0.1$ with $\epsilon=10^{-12}$.  The spatial discretization is DG4.  Left: BERK2 without the correction step.  Middle:  BERK2 with BDF2 correction step.  Right BERK2 with both BDF2 correction step and conservation fix. The bulk velocity $u$ is excluded since $u=1$. The influences of the BDF2 correction and the conservation fix on accuracy are very small in the asymptotic test.}
\centering
 \begin{tabular}{|c|c| c|| c| c|| c|c||c|c|}
  \hline
 \multicolumn{1}{|c}{} &\multicolumn{2}{|c||}{}&\multicolumn{2}{c||}{No correction}&\multicolumn{2}{c||}{BDF2 correction}&\multicolumn{2}{c|}{ BDF2 Correction + conservation fix}\\
 \hline
  &$\epsilon$ & $N_x$ & $L^2(\rho-\rho_{{\rm{Euler}}})$ & order& $L^2(\rho-\rho_{{\rm{Euler}}})$ & order& $L^2(\rho-\rho_{{\rm{Euler}}})$ & order\\
\hline\hline  
&&    8 & 2.376e-01 &    -    & 2.376e-01 &    -    & 2.376e-01 &    -    \\ 
&&   16 & 3.455e-02 &   2.781 & 3.455e-02 &   2.781 & 3.455e-02 &   2.781 \\ 
&&   32 & 2.612e-03 &   3.726 & 2.612e-03 &   3.726 & 2.612e-03 &   3.726 \\ 
$\rho$&$10^{-12}$&   64 & 1.078e-04 &   4.599 & 1.078e-04 &   4.599 & 1.078e-04 & 4.599 \\ 
&&  128 & 1.427e-05 &   2.917 & 1.425e-05 &   2.919 & 1.425e-05 &   2.919 \\ 
&&  256 & 3.261e-06 &   2.129 & 3.242e-06 &   2.136 & 3.242e-06 &   2.136 \\ 
&&  512 & 8.136e-07 &   2.003 & 7.960e-07 &   2.026 & 7.960e-07 &   2.026 \\ 
\hline\hline 
& &    8 & 2.542e-01&    -     & 2.542e-01 &    -     & 2.542e-01 &    -    \\ 
& &    16 &4.384e-02 & 2.536 & 4.384e-02 &   2.536    & 4.384e-02 &   2.536 \\ 
& &    32 &4.088e-03 &  3.423   & 4.088e-03 &   3.423 & 4.088e-03 &   3.423 \\ 
$\theta$&$10^{-12}$&    64 &3.189e-04 &    3.680 & 3.189e-04 &   3.680 & 3.189e-04 &   3.680 \\ 
& &    128 &6.083e-05 & 2.390  & 6.083e-05 &   2.390   & 6.083e-05 &   2.390 \\ 
& &    256 &1.467e-05 & 2.052  & 1.467e-05 &   2.052   & 1.467e-05 &   2.052 \\ 
& &    512 &3.487e-06 & 2.073 & 3.487e-06 &   2.073    & 3.487e-06 &   2.073 \\ 
\hline
 \end{tabular}
\label{tab:asymptotic_results_BERK2_all}
\end{table}
\subsection{Accuracy test}\label{subsec:accuracy}
In this test, we consider the smooth initial condition from \cite{xiong} to test the order of accuracy of the hybrid for arbitrary various $\epsilon$.  The domain is $X \times V=[-\pi,\pi] \times[-7,7]$, and the boundary conditions are periodic.   The initial data is $G^0(x,v) = \cE(\bm{q}(x,0))(v)$, where $\bm{q}$ is determined from the variables
\begin{equation}
	\label{eq:accuracy_test_ic}
	\rho(x,t=0)=1+0.2\sin10x, \qquad u(x,t=0) = 1,  \qquad p(x,t=0) =1.
\end{equation} 
The final time is $t=0.1$. A velocity grid with $N_v=100$ is used throughout.  Because the analytic solution is not available, the solution with the half mesh size is used as the reference solution.  We consider the $L^2$-errors for $\rho$, given by $E^{h}_2:=\|{\rho}^{h}-{\rho}^{h/2}\|_{L^{2}}$.  The convergence order $\nu$ is defined by $\nu={\log{\left({E^{h}_2}/{E^{h/2}_2}\right)}}/{\log{2}}$. 

The results are given in \tableref{tab:BERK2_all_accuracy_test_t01}.  While the standard BERK2 method exhibits degradation in convergence order, especially for $\epsilon \in \{1,0.1,0.01 \}$, the order improves significantly when the BDF2 correction step and the conservation fix are employed.
\begin{table}[ht!]\footnotesize
\caption{Accuracy test from Section \ref{subsec:accuracy}.  The  $L_2$-errors $E_h^2$ and convergence order $\nu$ of  the density for various $\epsilon$.  The time step size is $\Delta{t}=C \Delta{x} \Lambda^{-1}$, where $C = 0.2$ for DG2,  $C = 0.1$ for DG3, and $C = 0.05$ for DG4.  The final time is $t=0.1$.  BERK2 without the correction step is on the left and BERK2 with both the BDF2 correction step and the conservation fix is on the right.}
\centering
 \begin{tabular}{|c| c|| c| c| c| c| c| c||c| c| c| c| c| c|}
  \hline
  \multicolumn{2}{|c||}{}&\multicolumn{6}{c||}{No correction}&\multicolumn{6}{c|}{BDF2 correction + conservation fix}\\
  \hline
  \multicolumn{2}{|c||}{}&\multicolumn{2}{c|}{DG2}&\multicolumn{2}{c|}{DG3}&\multicolumn{2}{c||}{DG4}&\multicolumn{2}{c|}{DG2}&\multicolumn{2}{c|}{DG3}&\multicolumn{2}{c|}{DG4}\\
 \hline
  $\epsilon$ & $N_x$ & $E_h^2$ & $\nu$ & $E_h^2$ & $\nu$ & $E_h^2$ & $\nu$ & $E_h^2$ & $\nu$ & $E_h^2$ & $\nu$ & $E_h^2$ & $\nu$\\
\hline\hline
& 16 & 2.56e-02 & -  & 1.13e-03 & -  & 6.34e-05 & - & 2.61e-02 & -  & 1.16e-03 & -  & 6.40e-05 & - \\ 
& 32 & 1.04e-02 & 1.3  & 2.03e-04 & 2.5  & 4.84e-06 & 3.7 & 1.02e-02 & 1.3  & 1.93e-04 & 2.6  & 4.64e-06 & 3.8 \\ 
1& 64 & 3.76e-03 & 1.5  & 2.91e-05 & 2.8  & 3.39e-07 & 3.8 & 3.29e-03 & 1.6  & 2.46e-05 & 3.0  & 3.00e-07 & 3.9 \\ 
& 128 & 1.16e-03 & 1.7  & 4.38e-06 & 2.7  & 2.42e-08 & 3.8 & 8.34e-04 & 2.0  & 3.05e-06 & 3.0  & 1.89e-08 & 4.0 \\ 
& 256 & 3.82e-04 & 1.6  & 7.34e-07 & 2.6  & 1.88e-09 & 3.7 & 2.03e-04 & 2.0  & 3.80e-07 & 3.0  & 1.18e-09 & 4.0 \\ 
 \hline\hline
& 16 & 2.56e-02 & -  & 1.13e-03 & -  & 6.34e-05 & - & 2.61e-02 & -  & 1.16e-03 & -  & 6.40e-05 & - \\ 
& 32 & 1.03e-02 & 1.3  & 1.99e-04 & 2.5  & 4.85e-06 & 3.7 & 1.02e-02 & 1.3  & 1.95e-04 & 2.6  & 4.70e-06 & 3.8 \\ 
$10^{-1}$& 64 & 3.64e-03 & 1.5  & 2.81e-05 & 2.8  & 3.35e-07 & 3.9 & 3.31e-03 & 1.6  & 2.48e-05 & 3.0  & 3.04e-07 & 4.0 \\ 
& 128 & 1.11e-03 & 1.7  & 4.23e-06 & 2.7  & 2.37e-08 & 3.8 & 8.41e-04 & 2.0  & 3.07e-06 & 3.0  & 1.91e-08 & 4.0 \\ 
& 256 & 3.72e-04 & 1.6  & 7.18e-07 & 2.6  & 1.83e-09 & 3.7 & 2.03e-04 & 2.0  & 3.78e-07 & 3.0  & 1.18e-09 & 4.0 \\ 
 \hline\hline 
& 16 & 2.56e-02 & -  & 1.13e-03 & -  & 6.34e-05 & - & 2.61e-02 & -  & 1.16e-03 & -  & 6.40e-05 & - \\ 
& 32 & 1.08e-02 & 1.2  & 1.89e-04 & 2.6  & 5.07e-06 & 3.6 & 1.03e-02 & 1.3  & 1.94e-04 & 2.6  & 4.94e-06 & 3.7 \\ 
$10^{-2}$& 64 & 3.77e-03 & 1.5  & 3.23e-05 & 2.5  & 3.89e-07 & 3.7 & 3.37e-03 & 1.6  & 2.54e-05 & 2.9  & 3.20e-07 & 4.0 \\ 
& 128 & 1.75e-03 & 1.1  & 7.91e-06 & 2.0  & 3.77e-08 & 3.4 & 9.17e-04 & 1.9  & 3.27e-06 & 3.0  & 2.01e-08 & 4.0 \\ 
& 256 & 9.78e-04 & 0.8  & 2.09e-06 & 1.9  & 4.42e-09 & 3.1 & 2.37e-04 & 2.0  & 4.10e-07 & 3.0  & 1.24e-09 & 4.0 \\ 
 \hline\hline 
& 16 & 2.56e-02 & -  & 1.13e-03 & -  & 6.34e-05 & - & 2.61e-02 & -  & 1.16e-03 & -  & 6.40e-05 & - \\ 
& 32 & 1.28e-02 & 1.0  & 1.79e-04 & 2.7  & 5.17e-06 & 3.6 & 1.28e-02 & 1.0  & 1.79e-04 & 2.7  & 5.17e-06 & 3.6 \\ 
$10^{-6}$& 64 & 3.49e-03 & 1.9  & 2.34e-05 & 2.9  & 3.26e-07 & 4.0 & 3.49e-03 & 1.9  & 2.34e-05 & 2.9  & 3.26e-07 & 4.0 \\ 
& 128 & 8.67e-04 & 2.0  & 3.19e-06 & 2.9  & 2.21e-08 & 3.9 & 8.67e-04 & 2.0  & 3.19e-06 & 2.9  & 2.21e-08 & 3.9 \\ 
& 256 & 2.25e-04 & 1.9  & 6.75e-07 & 2.2  & 4.14e-09 & 2.4 & 2.25e-04 & 1.9  & 6.73e-07 & 2.2  & 4.09e-09 & 2.4 \\ 
 \hline 
 \end{tabular}
\label{tab:BERK2_all_accuracy_test_t01}
\end{table}
\begin{table}[ht!]
\caption{Comparisons of CPU time and speed-up for each problem. I3L, I2L, B2L, and B2LS denote IMEX3L, IMEX2L, BERK2L, and BERK2LS, respectively.}
\centering
\begin{tabular}{|l|l||rrrr||rr||rr|}
\hline
                               &            & \multicolumn{4}{c||}{CPU time (s)}                                                                       & \multicolumn{2}{c||}{Speed-up (B2L)}     & \multicolumn{2}{c|}{Speed-up (B2LS)}    \\ \hline
                               & $\epsilon$ & \multicolumn{1}{r|}{I3L}    & \multicolumn{1}{r|}{I2L}     & \multicolumn{1}{r|}{B2L}    & B2LS     & \multicolumn{1}{r|}{vs I3L}    & vs I2L    & \multicolumn{1}{r|}{vs I3L}    & vs I2L    \\ \hline\hline
\multirow{3}{*}{Sod}           & 1          & \multicolumn{1}{r|}{19.87}  & \multicolumn{1}{r|}{8.64}    & \multicolumn{1}{r|}{1.06}   & 2.29     & \multicolumn{1}{r|}{18.75}     & 8.15      & \multicolumn{1}{r|}{8.68}      & 3.77      \\ \cline{2-10} 
                               & $10^{-2}$  & \multicolumn{1}{r|}{20.94}  & \multicolumn{1}{r|}{8.71}    & \multicolumn{1}{r|}{1.11}   & 2.35     & \multicolumn{1}{r|}{18.86}     & 7.85      & \multicolumn{1}{r|}{8.91}      & 3.71      \\ \cline{2-10} 
                               & $10^{-6}$  & \multicolumn{1}{r|}{18.71}  & \multicolumn{1}{r|}{8.57}    & \multicolumn{1}{r|}{1.05}   & 2.31     & \multicolumn{1}{r|}{17.82}     & 8.16      & \multicolumn{1}{r|}{8.10}      & 3.71      \\ \hline\hline
\multirow{3}{*}{Lax}           & 1          & \multicolumn{1}{r|}{25.56}  & \multicolumn{1}{r|}{11.91}   & \multicolumn{1}{r|}{1.74}   & 3.71     & \multicolumn{1}{r|}{14.69}     & 6.84      & \multicolumn{1}{r|}{6.89}      & 3.21      \\ \cline{2-10} 
                               & $10^{-2}$  & \multicolumn{1}{r|}{24.99}  & \multicolumn{1}{r|}{12.16}   & \multicolumn{1}{r|}{1.74}   & 3.77     & \multicolumn{1}{r|}{14.36}     & 6.99      & \multicolumn{1}{r|}{6.63}      & 3.23      \\ \cline{2-10} 
                               & $10^{-6}$  & \multicolumn{1}{r|}{25.21}  & \multicolumn{1}{r|}{11.37}   & \multicolumn{1}{r|}{1.71}   & 3.82     & \multicolumn{1}{r|}{14.74}     & 6.65      & \multicolumn{1}{r|}{6.60}      & 2.98      \\ \hline\hline
\multirow{3}{*}{Shu-Osher}     & 1          & \multicolumn{1}{r|}{292.33} & \multicolumn{1}{r|}{120.93}  & \multicolumn{1}{r|}{29.77}  & 62.1     & \multicolumn{1}{r|}{9.82}      & 4.06      & \multicolumn{1}{r|}{4.71}      & 1.95      \\ \cline{2-10} 
                               & $10^{-2}$  & \multicolumn{1}{r|}{287.11} & \multicolumn{1}{r|}{121.02}  & \multicolumn{1}{r|}{29.50}  & 61.2     & \multicolumn{1}{r|}{9.73}      & 4.10      & \multicolumn{1}{r|}{4.69}      & 1.98      \\ \cline{2-10} 
                               & $10^{-6}$  & \multicolumn{1}{r|}{291.56} & \multicolumn{1}{r|}{120.36}  & \multicolumn{1}{r|}{29.92}  & 61.7     & \multicolumn{1}{r|}{9.74}      & 4.02      & \multicolumn{1}{r|}{4.73}      & 1.95      \\ \hline\hline
\multirow{3}{*}{gas-injection} & 1          & \multicolumn{1}{r|}{697.96} & \multicolumn{1}{r|}{399.83}  & \multicolumn{1}{r|}{3.09}   & 10.33    & \multicolumn{1}{r|}{225.88}    & 129.39    & \multicolumn{1}{r|}{67.57}     & 38.71     \\ \cline{2-10} 
                               & $10^{-2}$  & \multicolumn{1}{r|}{699.51} & \multicolumn{1}{r|}{404.01}  & \multicolumn{1}{r|}{6.09}   & 18.84    & \multicolumn{1}{r|}{114.86}    & 66.34     & \multicolumn{1}{r|}{37.13}     & 21.44     \\ \cline{2-10} 
                               & $10^{-6}$  & \multicolumn{1}{r|}{700.4}  & \multicolumn{1}{r|}{405.88}  & \multicolumn{1}{r|}{10.21}  & 33.85    & \multicolumn{1}{r|}{68.60}     & 39.75     & \multicolumn{1}{r|}{20.69}     & 11.99     \\ \hline
\end{tabular}
\label{tab:speedup}
\end{table}
\subsection{Sod shock tube problem}\label{subsec:sod}
The Sod shock tube problem \cite{sod} is used to determine how well a numerical solver captures fundamental flow features.  The initial conditions are $G^0(x,v) = \cE(\bm{q}(x,0))(v)$, where $\bm{q}$ is determined from the variables
\begin{equation}
(\rho(x,0),u(x,0),p(x,0))=\begin{cases}
(1,0,1), \hspace{18mm}0\leq x\leq 0.5,\\
(0.125,0,0.1), \hspace{8mm}0.5< x\leq 1.
\end{cases}
\end{equation}
The spatial
domain $[0, 1]$ is discretized with $N_x = 100$ spatial cells and the velocity domain $[-6, 6]$ is
discretized with a $N_v = 100$ point Gauss-Legendre quadrature set, scaled to fit the velocity domain. The boundary conditions are chosen to coincide with the initial constant moments at both ends of the spatial domain. The final time is $t = 0.1$.  Numerical solutions are computed with the hybrid BERK2 time discretization and a 3rd-order DG discretization in space.

\figref{fig:comparison_Sod_diff_plot} includes five solutions. EulerH is the Euler solution that is obtained by the Euler solver from DoGPack \cite{dogpack} with $N_x=3000$.  IMEX3H is a highly resolved IMEX3-DG solution with $N_x=N_v=1000$; it is used as a reference.  IMEX3L is an IMEX3-DG solution with $N_x=N_v=100$. IMEX2L is an IMEX3-DG solution with $N_x=N_v=100$. BERK2L is the hybrid BERK2 solution with $N_x=N_v=100$.  All schemes are plotted with the CFL constant $C$ given in \tableref{tab:comparisons_dt}. The differences in wave speeds, time-step size, and the size of the velocity domain are presented in \tableref{tab:comparisons_dt} to facilitate an effective comparison of efficiency for all of the test problems in Sections \ref{subsec:sod}-\ref{subsec:gas_injection}.

We plot three differences to compare the accuracy and CPU time for each method of the IMEX-DG and BERK2 methods:  (i) IMEX3H-IMEX3L, (ii) IMEX3H-IMEX2L, and (iii) IMEX3H-BERK2L.  From the results in \figref{fig:comparison_Sod_diff_plot}, we conclude that the BERK2L scheme provides a good approximation within a smaller CPU time than IMEX2L and IMEX3L.  In this particular velocity domain setting, BERK2L shows the CPU time corresponding to approximately $1/8$ of the CPU time for IMEX2L and $1/18$ of the CPU time for IMEX3L as shown in 
 \tableref{tab:speedup}.  The only significant differences in accuracy occur when $\epsilon=0.01$.  In that instance, the IMEX-DG schemes provide better accuracy.  For this reason, we have included a plot of the three solutions together to highlight that these differences occur primarily near the sharp transitions in the solution profile. The fourth row (BERK2LS) of \figref{fig:comparison_Sod_diff_plot} emphasizes the BERK2L solutions with a smaller CFL constant of $C=0.1$, which yields comparable accuracy with less computation time than the IMEX solutions. Within the BERK2LS configuration, the CPU time is roughly 1/3 that of IMEX2L and approximately 1/8 the time of IMEX3L.
\subsection{Lax shock tube problem}\label{subsec:lax}
The Lax shock tube problem \cite{lax} is another Riemann problem that is similar to the previous Sod problem, but uses a non-zero initial bulk velocity. The initial conditions are $G^0(x,v) = \cE(\bm{q}(x,0))(v)$, where $\bm{q}$ is determined from the variables
\begin{equation}
(\rho(x,0),u(x,0),p(x,0))=\begin{cases}
(0.445,0.698,3.528), \hspace{8mm}-0.5\leq x\leq 0.5,\\
(0.5,0,0.571), \hspace{21.1mm}0.5< x\leq 1.5.
\end{cases}
\end{equation}
The spatial domain $[-0.5,1.5]$ is discretized with $N_x = 100$ spatial cells and the velocity domain $[-15,15]$ is discretized with $N_v = 100$ points. Numerical results are obtained through the 3rd order nodal DG scheme up to time $t = 0.1$. The boundary conditions are set to be the same way as the previous Sod shock tube problem.

Similar to the Sod problem, \figref{fig:comparison_Lax_diff_plot} includes five solutions. EulerH is the Euler solution that is obtained by the Euler solver from DoGPack \cite{dogpack} with $N_x=3000$.  IMEX3H is a highly resolved IMEX3-DG solution with $N_x=N_v=1000$ that will be used as a reference solution.  IMEX2L is an IMEX3-DG solution with $N_x=N_v=100$. IMEX2L is an  IMEX3-DG solution with $N_x=N_v=100$. BERK2L is the hybrid BERK2 solution with $N_x=N_v=100$.  All schemes are plotted with the CFL constant C given in \tableref{tab:comparisons_dt}.

We plot three differences to compare the accuracy and CPU time for each method of the IMEX-DG and BERK2 methods:  (i) IMEX3H-IMEX3L, (ii) IMEX3H-IMEX2L, and (iii) IMEX3H-BERK2L.  The results in  \figref{fig:comparison_Lax_diff_plot} are similar to those from the Sod problem.  Again, we conclude that the BERK2L scheme provides a good approximation within a smaller CPU time than IMEX2L and IMEX3L.  In this particular velocity domain setting, BERK2L shows the CPU time corresponding to approximately $1/6$ of the CPU time for IMEX2L and $1/14$ of the CPU time for IMEX3L  as shown in \tableref{tab:speedup}.  Again, the only significant differences in accuracy occur when $\epsilon=0.01$.  In that instance, the IMEX-DG schemes provide better accuracy.  To fix this, we have included the fourth row in which the BERK2LS solutions are presented with a reduced CFL constant of $C=0.1$, resulting in comparable accuracy with less computation time compared to the IMEX solutions. In the context of the BERK2LS framework, the CPU time is estimated to be about 1/3 of the time required for IMEX2L and roughly 1/6 of the time required for IMEX3L.
\subsection{Shu-Osher problem}\label{subsec:shu-osher}
The Shu-Osher problem \cite{shu} is a simulation of a shock-turbulence interaction where a shock propagates into a perturbed density field, thus it can be used to determine the ability of a numerical solver to capture a shock wave, its interaction with an unsteady density field, and the waves propagating downstream of the shock \cite{burger}. The initial conditions are $G^0(x,v) = \cE(\bm{q}(x,0))(v)$, where $\bm{q}$ is determined from the variables
\begin{equation}
(\rho(x,0),u(x,0),p(x,0))=\begin{cases}
(1.756757, 2.005122, 10.333333), & \quad  x\leq -4,\\
(1+0.2\sin{5x},0,1),  \quad  x>-4.
\end{cases}
\end{equation}
The initial shock is located at $x=-4$. In this problem, more spatial and angular grids are used due to the increased size of the spatial and velocity domain. The spatial domain $[-10,10]$ and the velocity domain $[-14,14]$ are discretized with $N_x = 200$ uniform grid points and $N_v = 200$ Gauss-Legendre quadrature points, respectively. We compute the solutions with the 3rd order DG scheme up to time $t = 1.8$. The same boundary conditions are used as before. 

Similar to the Sod and Lax problems, \figref{fig:comparison_ShuOsher_diff_plot} includes five solutions. EulerH is the Euler solution that is obtained by the Euler solver from DoGPack \cite{dogpack} with $N_x=3000$.  IMEX3H is a highly resolved IMEX3-DG solution with $N_x=N_v=1000$ that will be used as a reference solution. All IMEX2L, IMEX3L, and BERK2L solutions are computed with $N_x=N_v=200$. All schemes are plotted with the CFL constant C given in \tableref{tab:comparisons_dt}.

We plot three differences to compare the accuracy and CPU time for each method of the IMEX-DG and BERK2 methods in the same order as before.  The results in  \figref{fig:comparison_ShuOsher_diff_plot} are similar to those from the Sod problem.  Again, we conclude that the BERK2L scheme provides a good approximation within a smaller CPU time than IMEX2L and IMEX3L.  In this particular velocity domain setting, BERK2L shows the CPU time corresponding to $1/4$ of the CPU time for IMEX2L and $1/9$ of the CPU time for IMEX3L as shown in \tableref{tab:speedup}.  In this problem, we do not observe significant differences in accuracy found in the previous problems when $\epsilon=0.01$. However, we still include BERK2LS in the 4th row to demonstrate the efficiency of the BERK2 scheme. Operating within the BERK2LS environment, the computational time used is nearly half of what is used for IMEX2L and close to 1/4 of the time consumed by IMEX3L.
\subsection{Gas injection problems}\label{subsec:gas_injection}
The gas injection problem was designed to emphasize the benefit of the BERK2 method when $\vmax$ is very large relative to $\Lambda$, i.e., the ratio $\vmax/\Lambda$ is large. 
The initial condition is 
\begin{equation}
 G^0(x,v) = \cE(\bm{q}_1(x,0))(v),
\end{equation}
where $\bm{q}_1$ is determined from the variables 
\begin{equation}
	\rho_1(x,0) = 1, \qquad u_1(x,0) = 0 ,\qquad \theta_1(x,0) = 0.1.
\end{equation}
The source is 
\begin{equation}
 S(x,v,t) = \eta(x)\cE(\bm{q_2}(x,0))(v),
\end{equation}
where $\bm{q}_2$ is determined from the variables 
\begin{equation}
	\rho_2(x,0) = 0.01, \qquad u_2(x,0) = 100 ,\qquad \theta_2(x,0) = 100,
\end{equation}
and $\eta$ is
\begin{align}
    &\eta(x)=c_0\exp{\left(\frac{-(x-x_0)^2}{2\sigma^2}\right)},\quad x_0=0.5,\hspace{2mm}\sigma=0.1.
\end{align}
Here, $c_0$ is a normalization constant, defined such that $\int_0^1\eta(x)\,dx=1$. 
This problem requires a large velocity domain by construction. The spatial domain and the velocity domain are chosen to be $[-3,19]$ and $[-\vmax,\vmax]$, respectively, where $\vmax=110$. We present our results for BERK2L in Figures \ref{fig:comparison_gas_injection_diff_wide_rho_plot}-\ref{fig:comparison_gas_injection_diff_wide_T_plot} and BERK2LS in Figures \ref{fig:comparison_gas_injection_diff_S_wide_rho_plot}-\ref{fig:comparison_gas_injection_diff_S_wide_T_plot}.

The first rows of  Figures \ref{fig:comparison_gas_injection_diff_wide_rho_plot}-\ref{fig:comparison_gas_injection_diff_S_wide_T_plot} use the IMEX3H solutions as a reference. The remaining rows of \figref{fig:comparison_gas_injection_diff_wide_rho_plot}-\ref{fig:comparison_gas_injection_diff_S_wide_T_plot} show the $\rho$, $u$, and $\theta$ of the gas injection problem, and the leftmost column tells the computational time for the IMEX3L, IMEX2L, and hybrid BERK2L schemes. 
Setting $u_2=100$ greatly magnifies the computation time differences between BERK2L and IMEX2,3L schemes due to the increased size of the velocity domain with $v_{\text{max}}=110$. We choose more angular grids with $N_v=1000$ to correctly capture the velocity integration with minimal quadrature error. The spatial grid is chosen to be uniform with $N_x=200$. All schemes are plotted with the CFL constant C given in \tableref{tab:comparisons_dt}.

\tableref{tab:speedup} shows CPU time for the gas injection problem, and it is clear to see the time difference between the hybrid BERK2L and the rest. On average, for the $u_2 = 100$ case, the speed-up is about 225.88X, 114.86X, and 68.60X versus IMEX3L and 129.39X, 66.34X, and 39.75X versus IMEX2L when $\epsilon=1.0$, $\epsilon=10^{-2}$, and $\epsilon=10^{-6}$, respectively.  The corresponding solution profiles are shown in \figref{fig:comparison_gas_injection_diff_wide_rho_plot}-\ref{fig:comparison_gas_injection_diff_wide_T_plot}. The figure format in this section has been altered from earlier sections to enhance the visualization of the shock wave front. It is clear from these plots that the BERK2L solution has smeared out the high-speed structure due to stepping over the fastest time scales in the problem when $\epsilon=1$ and $\epsilon=10^{-2}$. This is the price to be paid for the efficiency and stability of the implicit calculation of the uncollided solution. The required number of time-step in BERK2 is significantly smaller than IMEX schemes as it is shown in \tableref{tab:speedup}, i.e., BERK2 can handle a much larger time-step size than IMEX schemes due to the relaxed CFL condition in this problem. 

For BERK2LS, which uses a small CFL constant $C=0.025$, the speed-up is 67.57X, 37.13X, and 20.69X versus IMEX3L and 38.71X, 21.44X, and 11.99X versus IMEX2L when $\epsilon=1.0$, $\epsilon=10^{-2}$, and $\epsilon=10^{-6}$, respectively.  The corresponding solution profiles are shown in \figref{fig:comparison_gas_injection_diff_S_wide_rho_plot}-\ref{fig:comparison_gas_injection_diff_S_wide_T_plot}.  In this case, the profile accuracy (BERK2LS) is more comparable to IMEX3H. Note that, for this particular problem, the explicit advection of IMEX schemes handle the bump-on-tail more accurately, though the errors are small in BERK2LS.

%% file: sec5_conclusion.tex
\section{Conclusion}\label{sec:conclusion}
 Our hybrid BERK2 scheme for the uncollided-collided decomposition works uniformly for all Knudsen numbers. The comparisons of the scheme with non-splitting IMEX schemes indicate that the hybrid BERK2 scheme is more effective than the IMEX scheme with a coarse spatial discretization especially when the ratio of $v_\text{max}$ to $\Lambda$ is large.
 
 Fully implicit treatment for the uncollided equation can be easily achieved by the linearity of the uncollided equation. This approach can remove the velocity CFL restriction existing in IMEX time discretizations, which is especially beneficial for problems that have a large ratio of $v_\text{max}$ to $\Lambda$. The obtained uncollided solution is updated as a source term in the collided equation. This approach makes the collided equation possible to be solved semi-implicitly with explicit calculations, which allows much less restrictive CFL conditions. The biggest potentiality of this scheme is that the kinetic equations can be coupled to existing high-order hydrodynamic codes without suffering from very restrictive CFL conditions. The splitting error in our hybrid scheme can be further reduced by the method, e.g., the integral deferred correction scheme \cite{crockatt}. Reducing the splitting error of the hybrid method is another interesting topic that can be studied in the future.

\begin{figure}[ht!]
  \centering
  \raisebox{35pt}{\parbox[b]{.15\textwidth}{IMEX3H}}%
  \subfloat[][$\rho$]{\includegraphics[width=.25\textwidth]{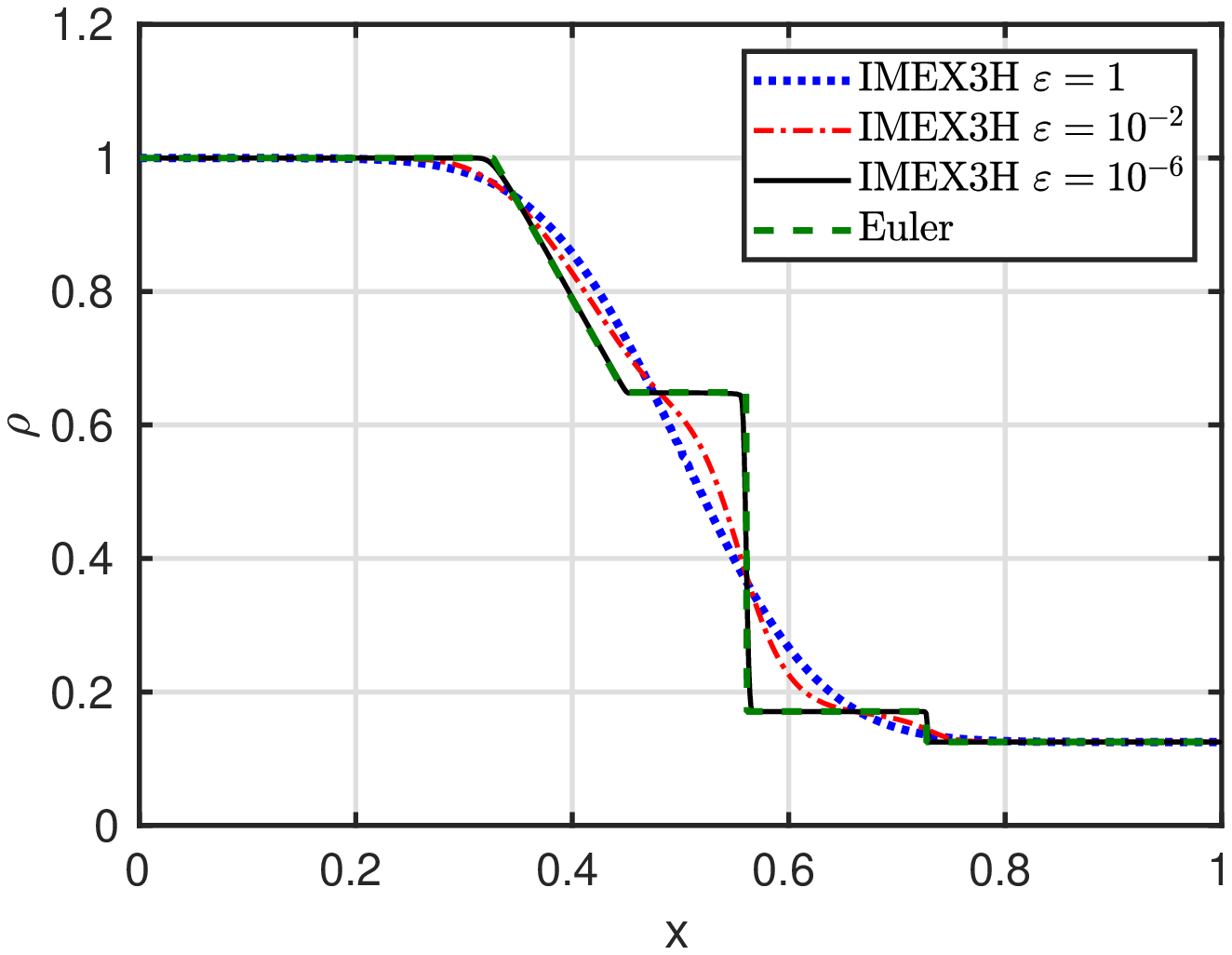}}
  \subfloat[][$u$]{\includegraphics[width=.25\textwidth]{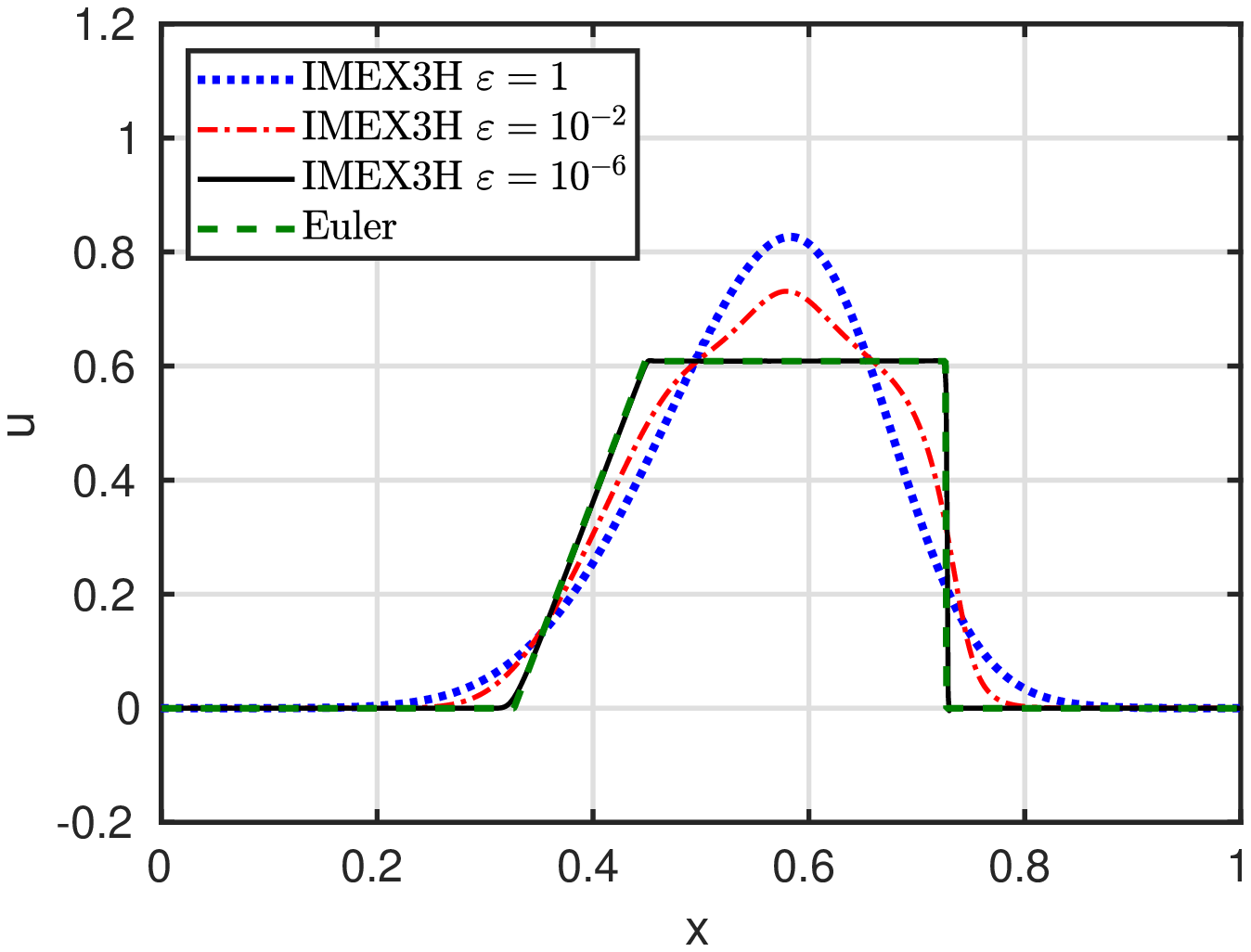}}
  \subfloat[][$\theta$]{\includegraphics[width=.25\textwidth]{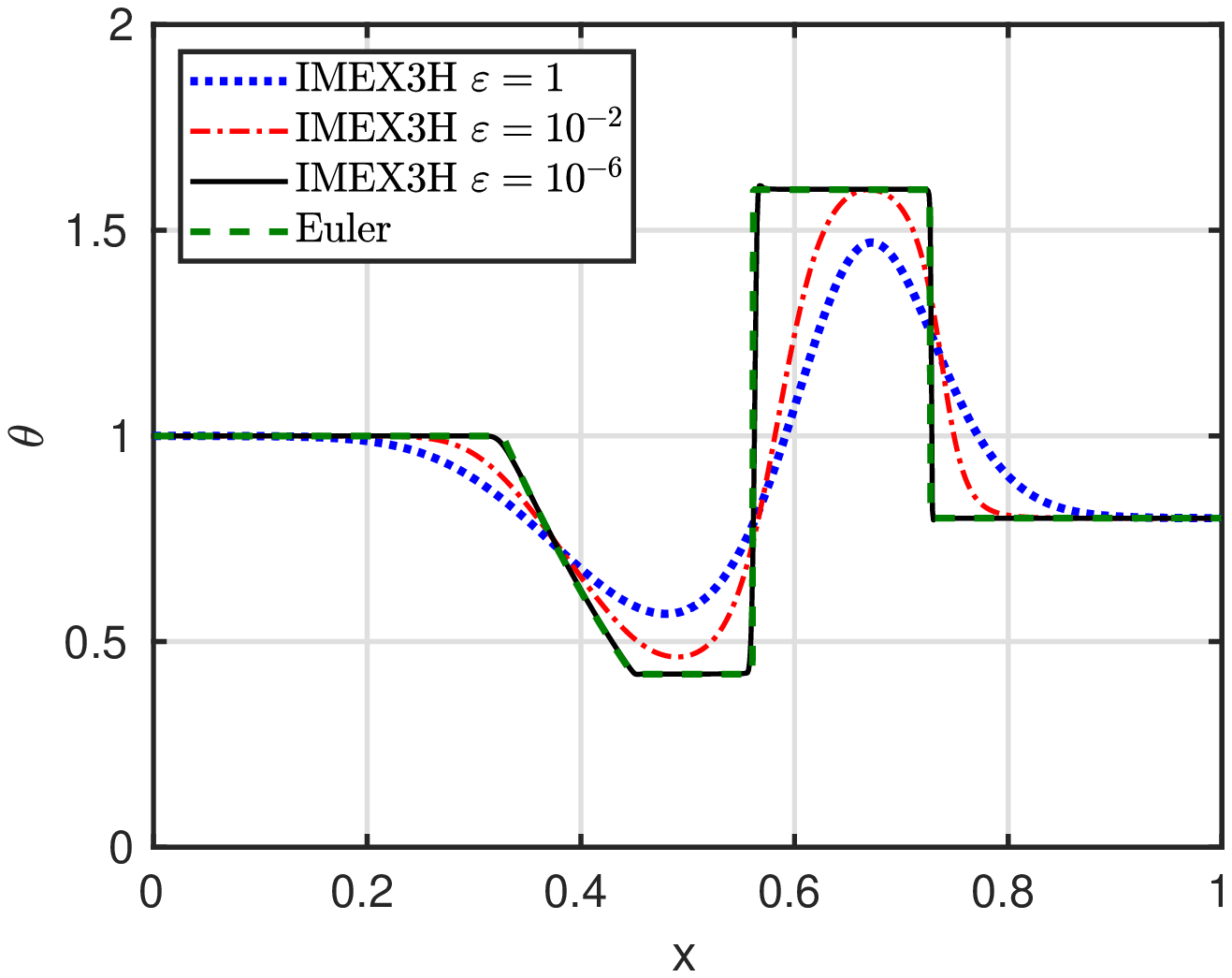}}\par
  \raisebox{35pt}{\parbox[b]{.15\textwidth}{\tiny$\epsilon=1$\\IMEX3L (C=0.14): 19.87s\\IMEX2L (C=0.2): 8.64s\\BERK2L (C=0.2): 1.06s}}%
  \subfloat[][$\rho$]{\includegraphics[width=.25\textwidth]{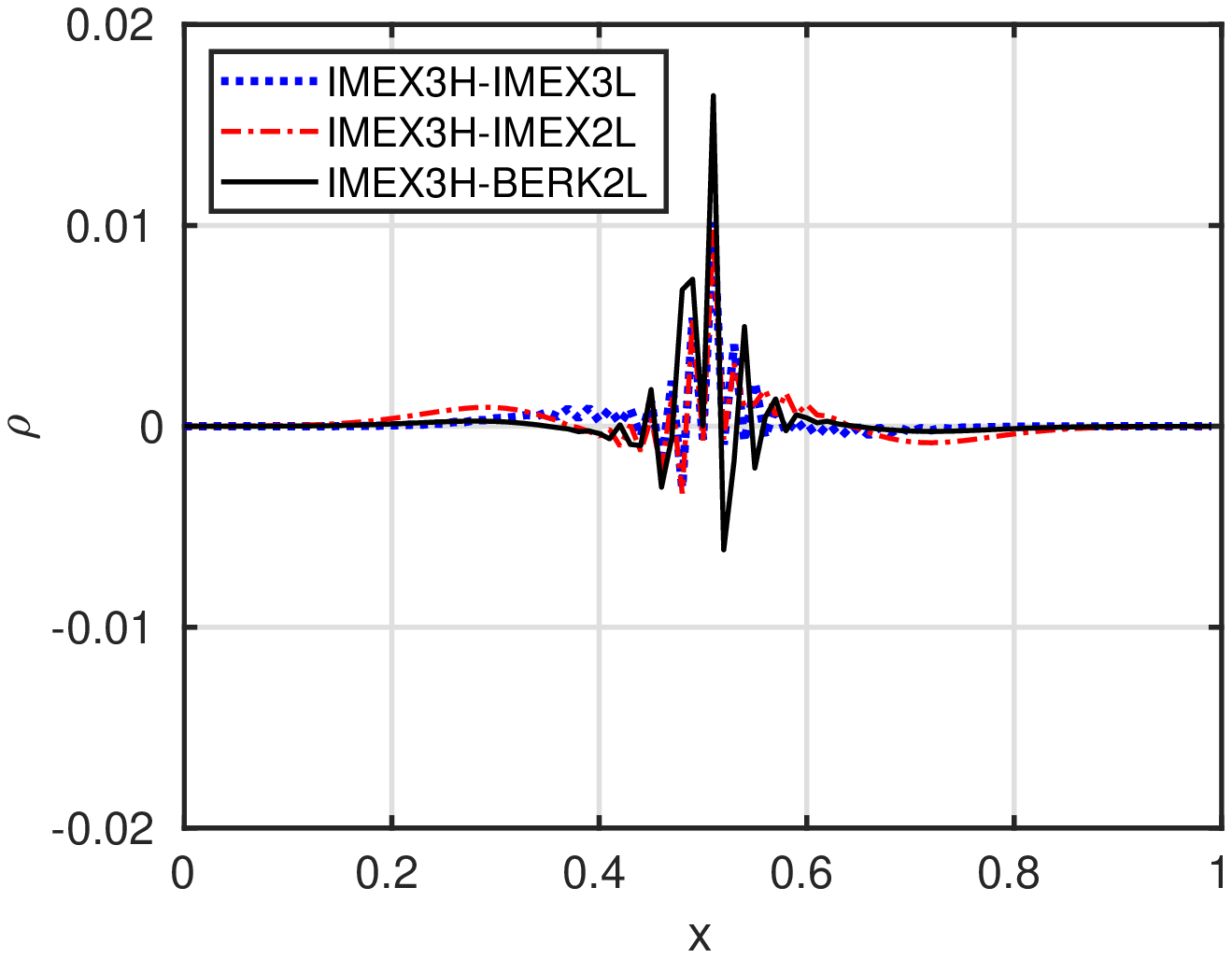}}
  \subfloat[][$u$]{\includegraphics[width=.25\textwidth]{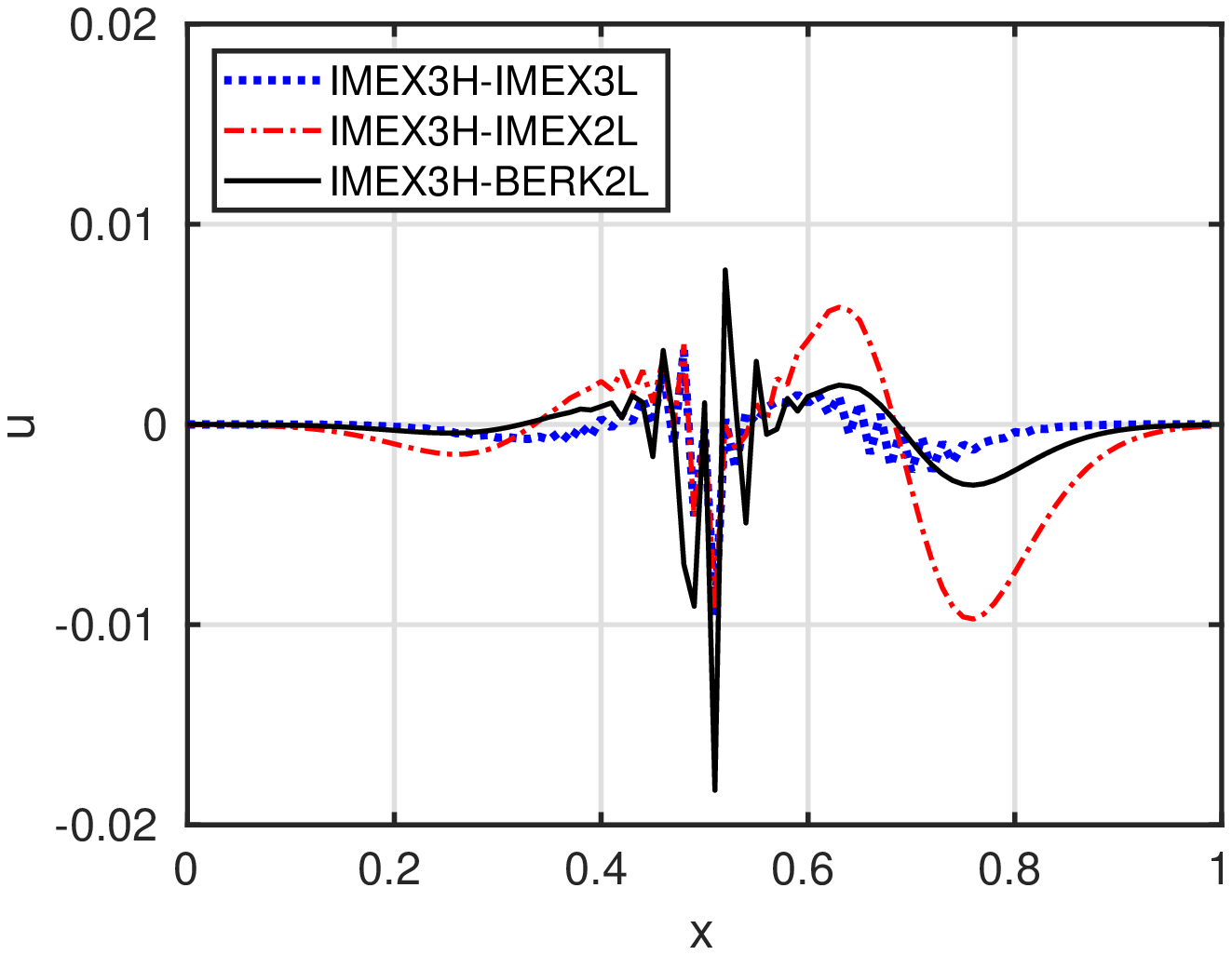}}
  \subfloat[][$\theta$]{\includegraphics[width=.25\textwidth]{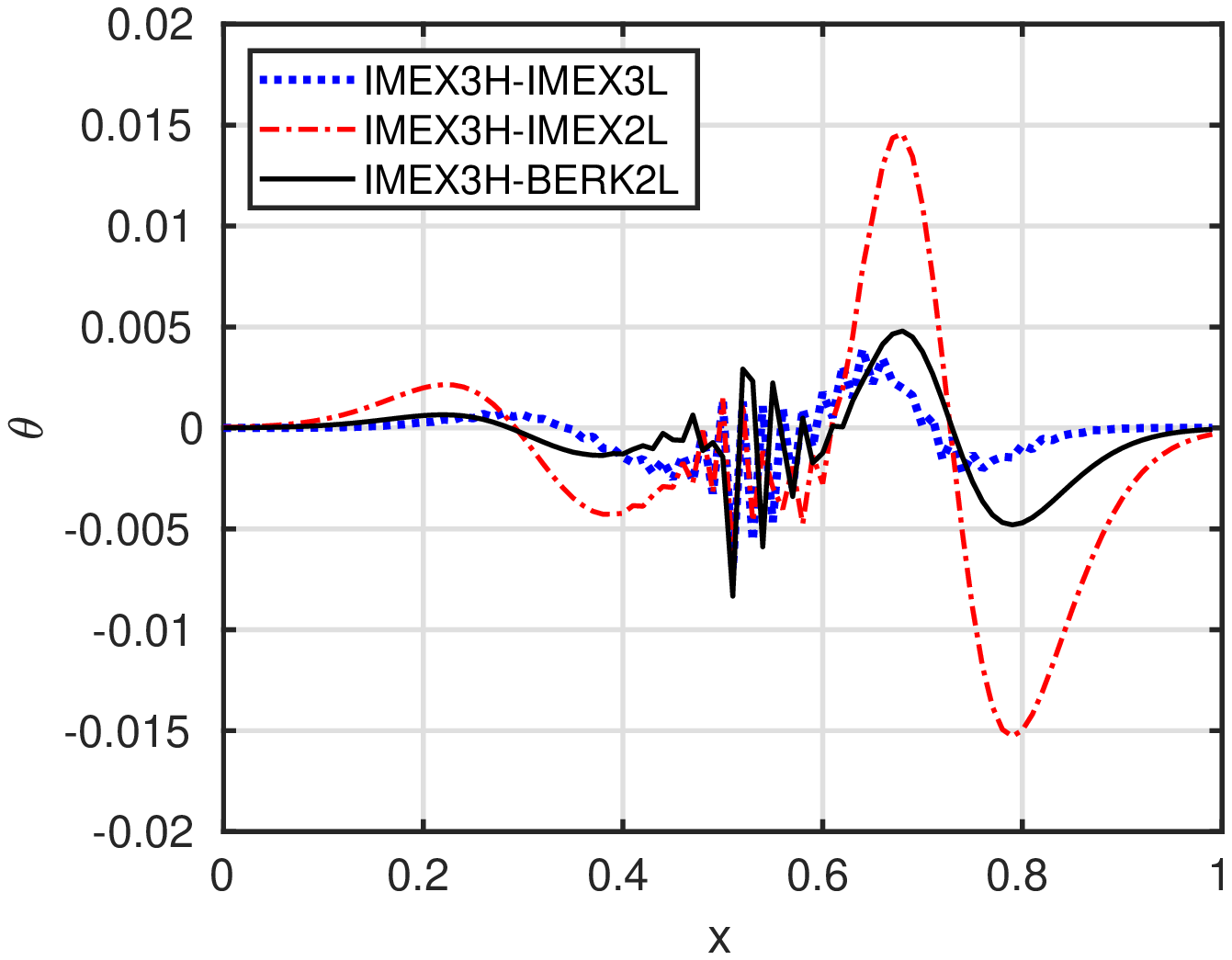}}\par
  \raisebox{35pt}{\parbox[b]{.15\textwidth}{\tiny$\epsilon=10^{-2}$\\IMEX3L (C=0.14): 20.94s\\IMEX2L (C=0.2): 8.71s\\BERK2L (C=0.2): 1.11s}}%
  \subfloat[][$\rho$]{\includegraphics[width=.25\textwidth]{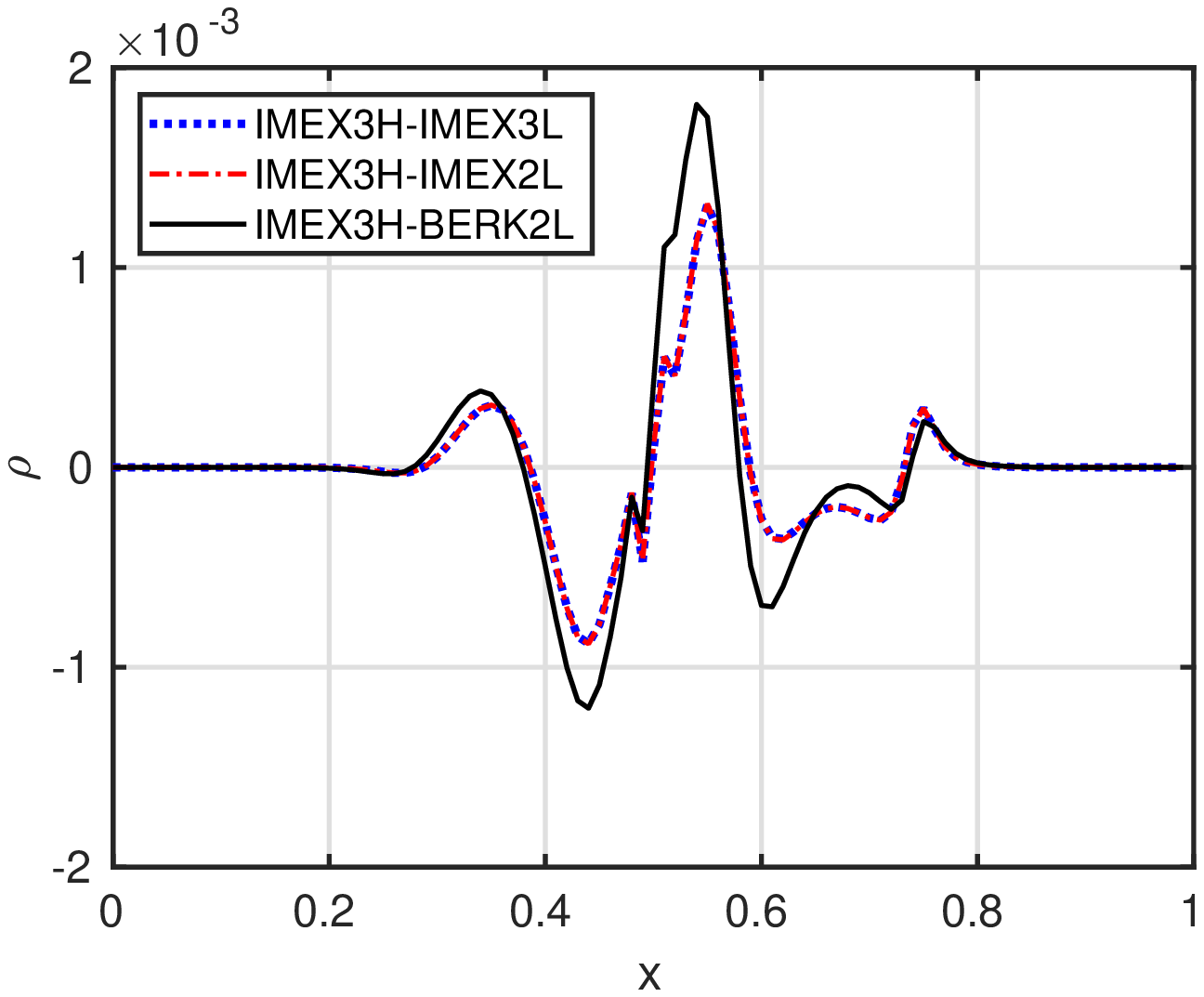}}
  \subfloat[][$u$]{\includegraphics[width=.25\textwidth]{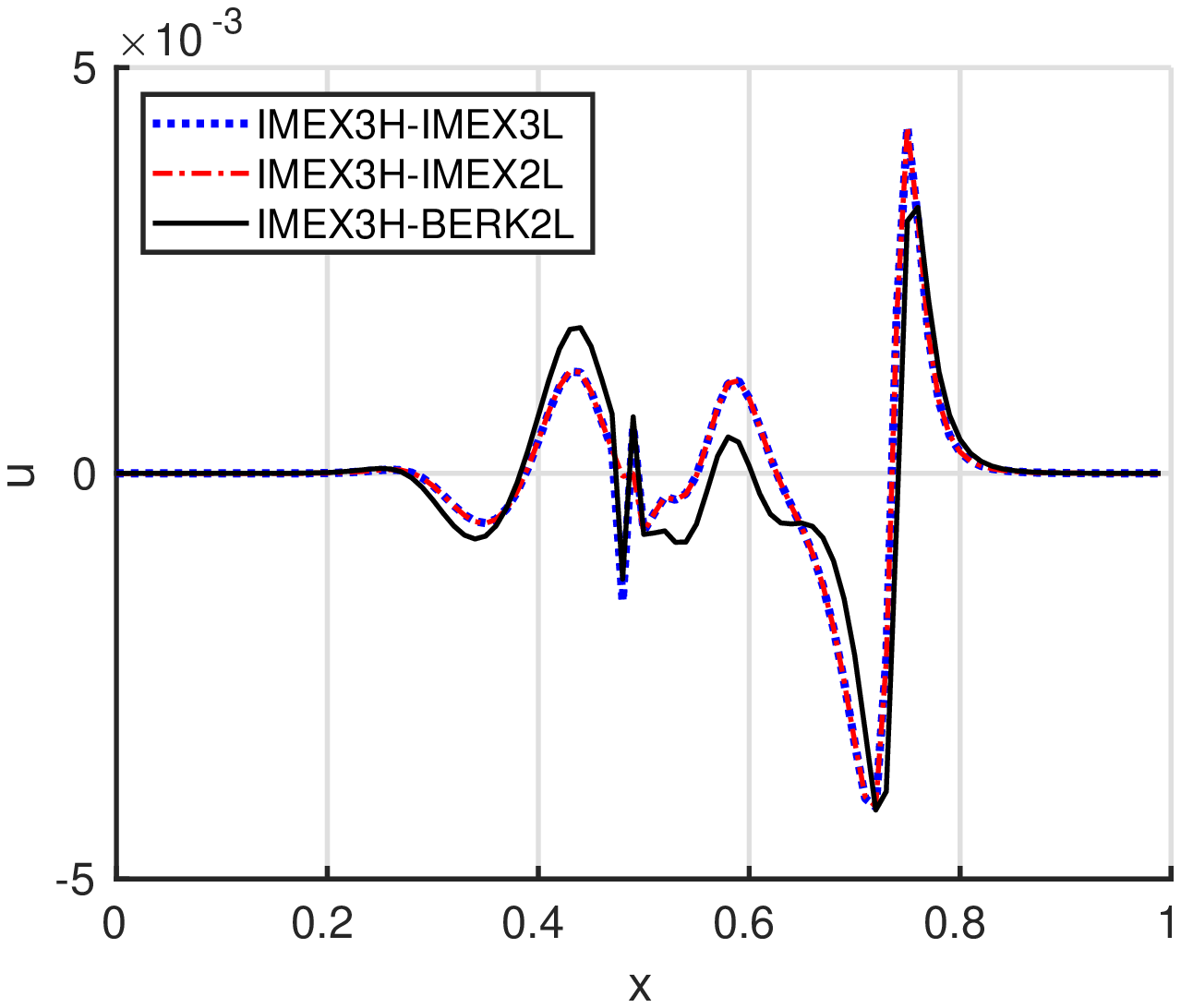}}
  \subfloat[][$\theta$]{\includegraphics[width=.25\textwidth]{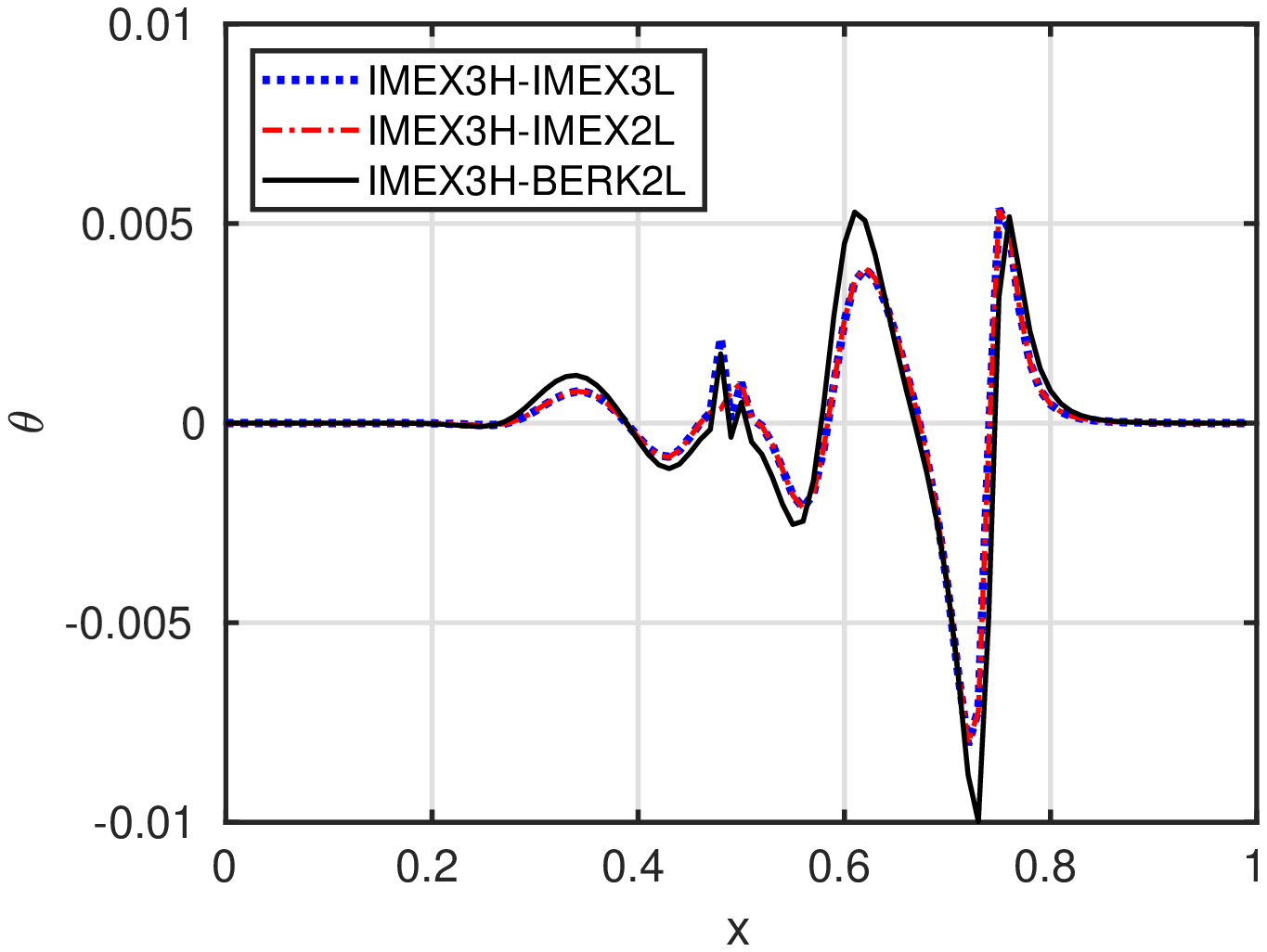}}\par
    \raisebox{35pt}{\parbox[b]{.15\textwidth}{\tiny$\epsilon=10^{-2}$\\IMEX3L (C=0.14): 20.94s\\IMEX2L (C=0.2): 8.71s\\BERK2LS (C=0.1): 2.35s}}%
  \subfloat[][$\rho$]{\includegraphics[width=.25\textwidth]{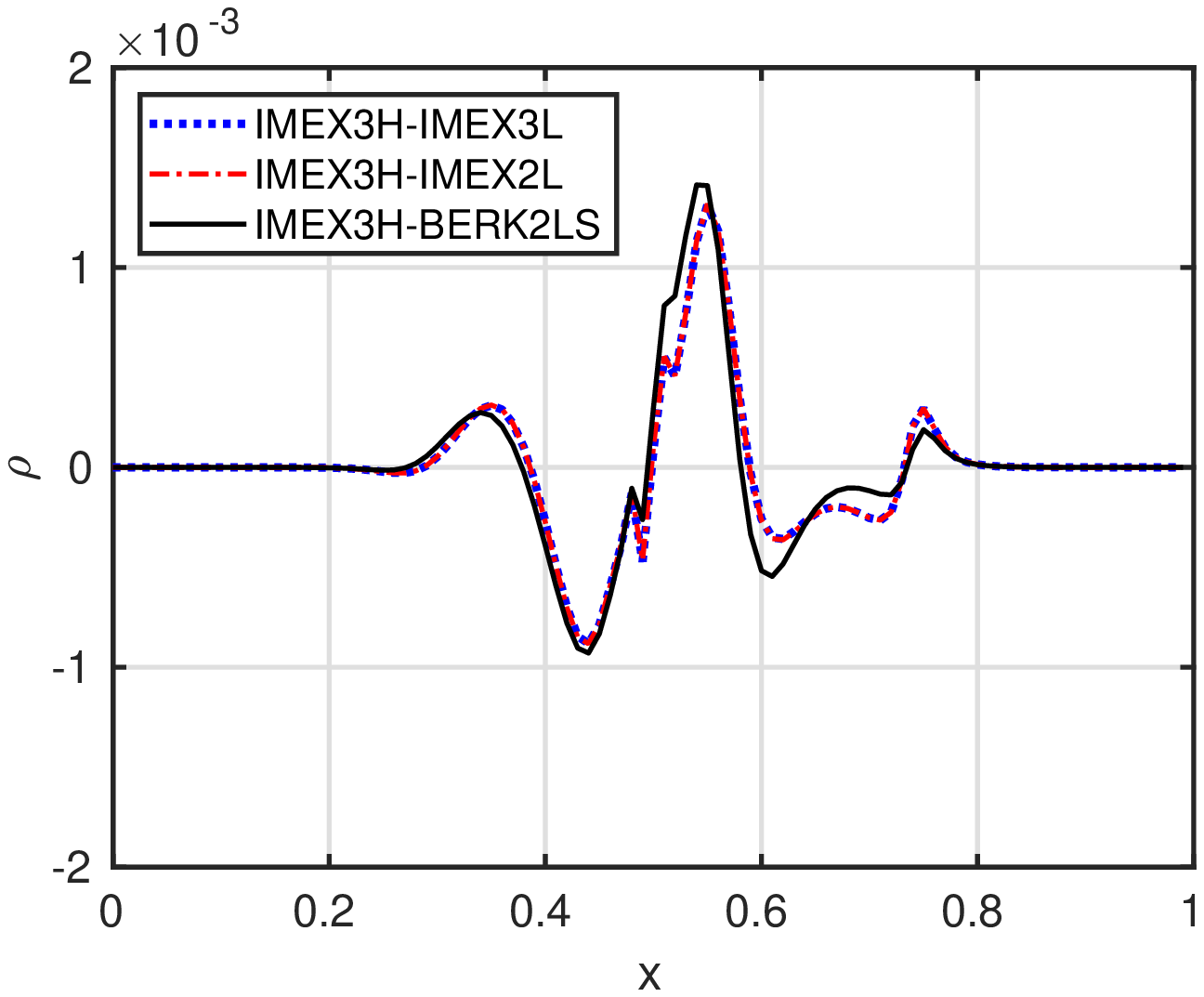}}
  \subfloat[][$u$]{\includegraphics[width=.25\textwidth]{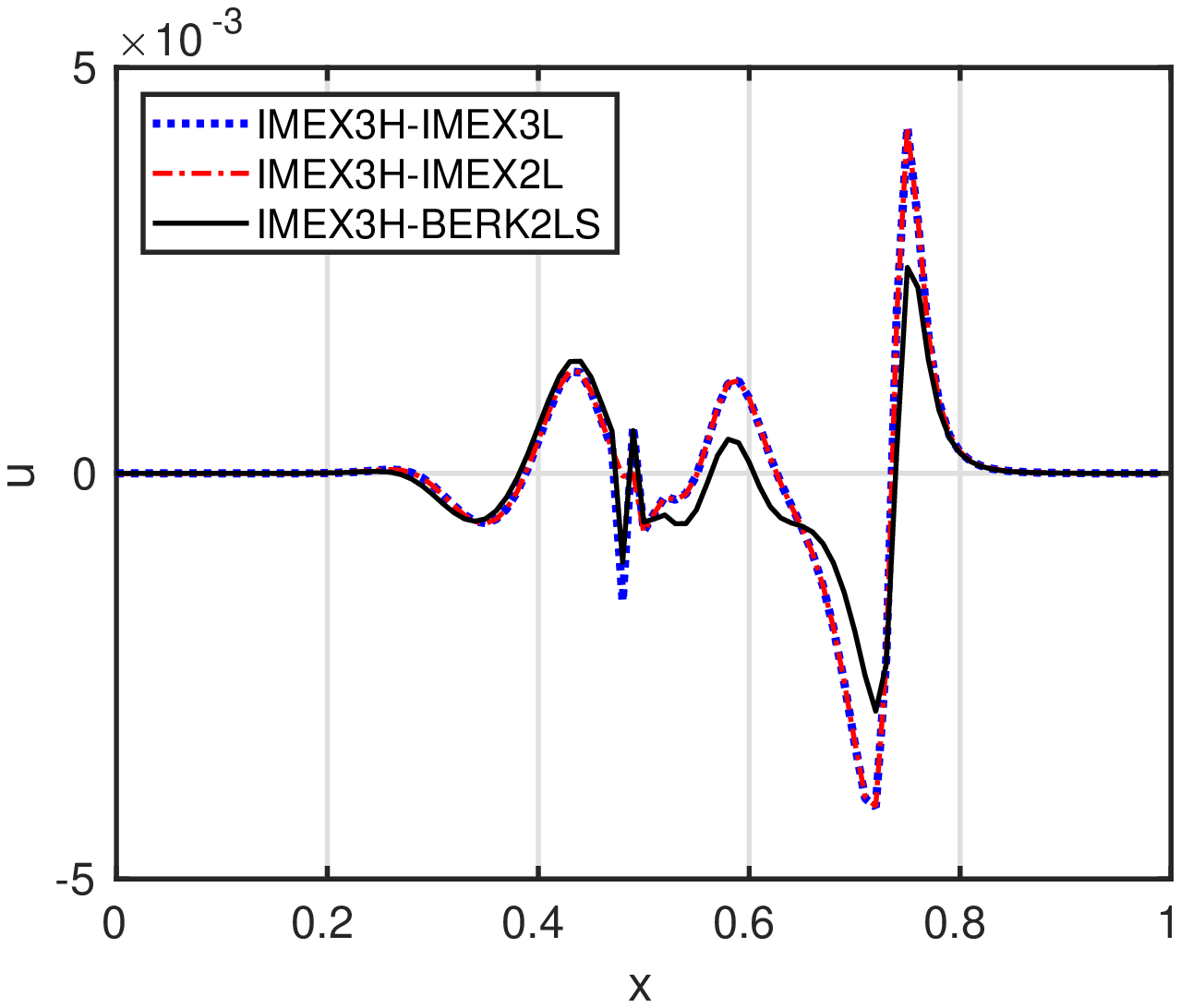}}
  \subfloat[][$\theta$]{\includegraphics[width=.25\textwidth]{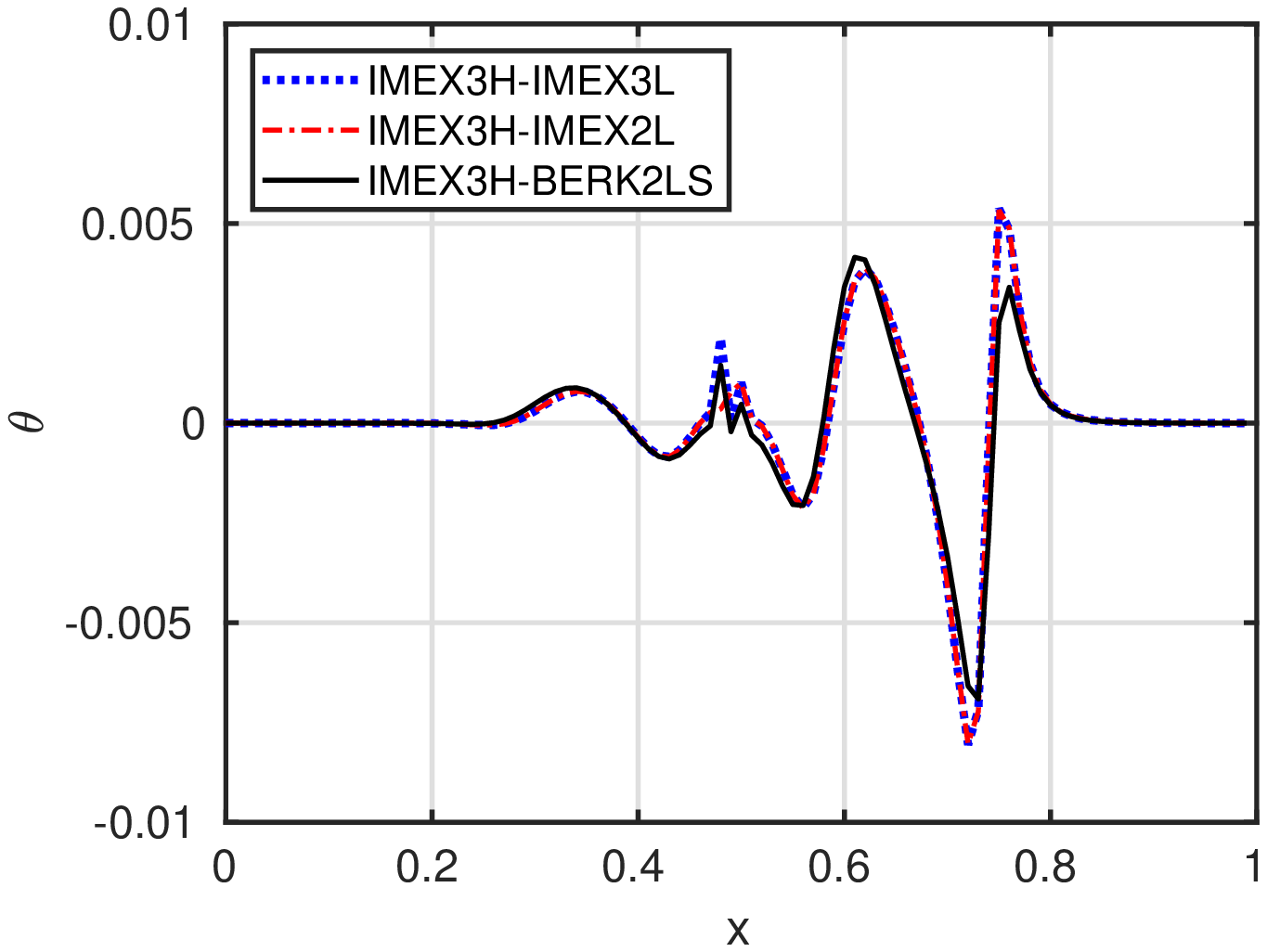}}\par
  \raisebox{35pt}{\parbox[b]{.15\textwidth}{\tiny$\epsilon=10^{-6}$\\IMEX3L (C=0.14): 18.71s\\IMEX2L (C=0.2): 8.57s\\BERK2L (C=0.2): 1.05s}}%
  \subfloat[][$\rho$]{\includegraphics[width=.25\textwidth]{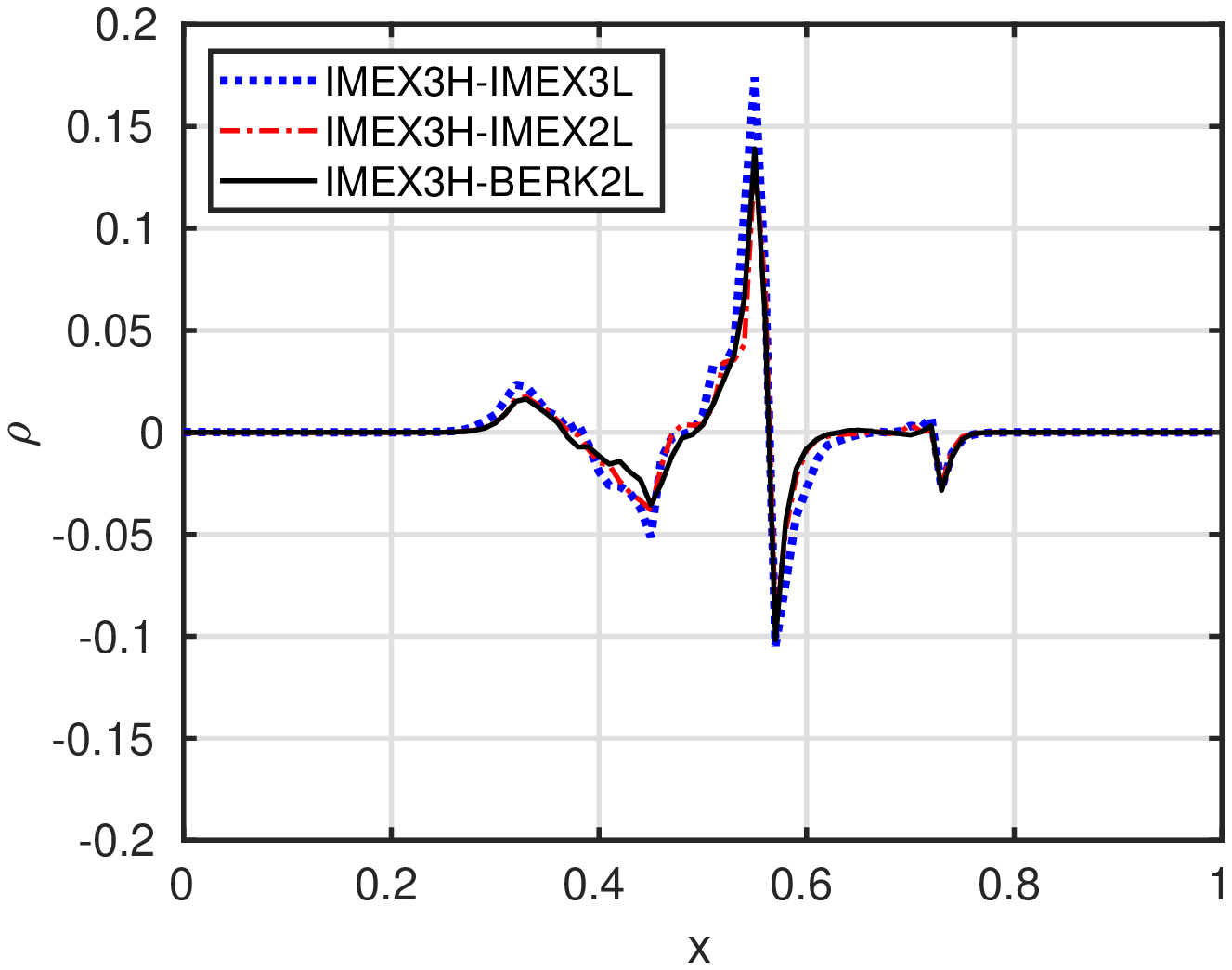}}
  \subfloat[][$u$]{\includegraphics[width=.25\textwidth]{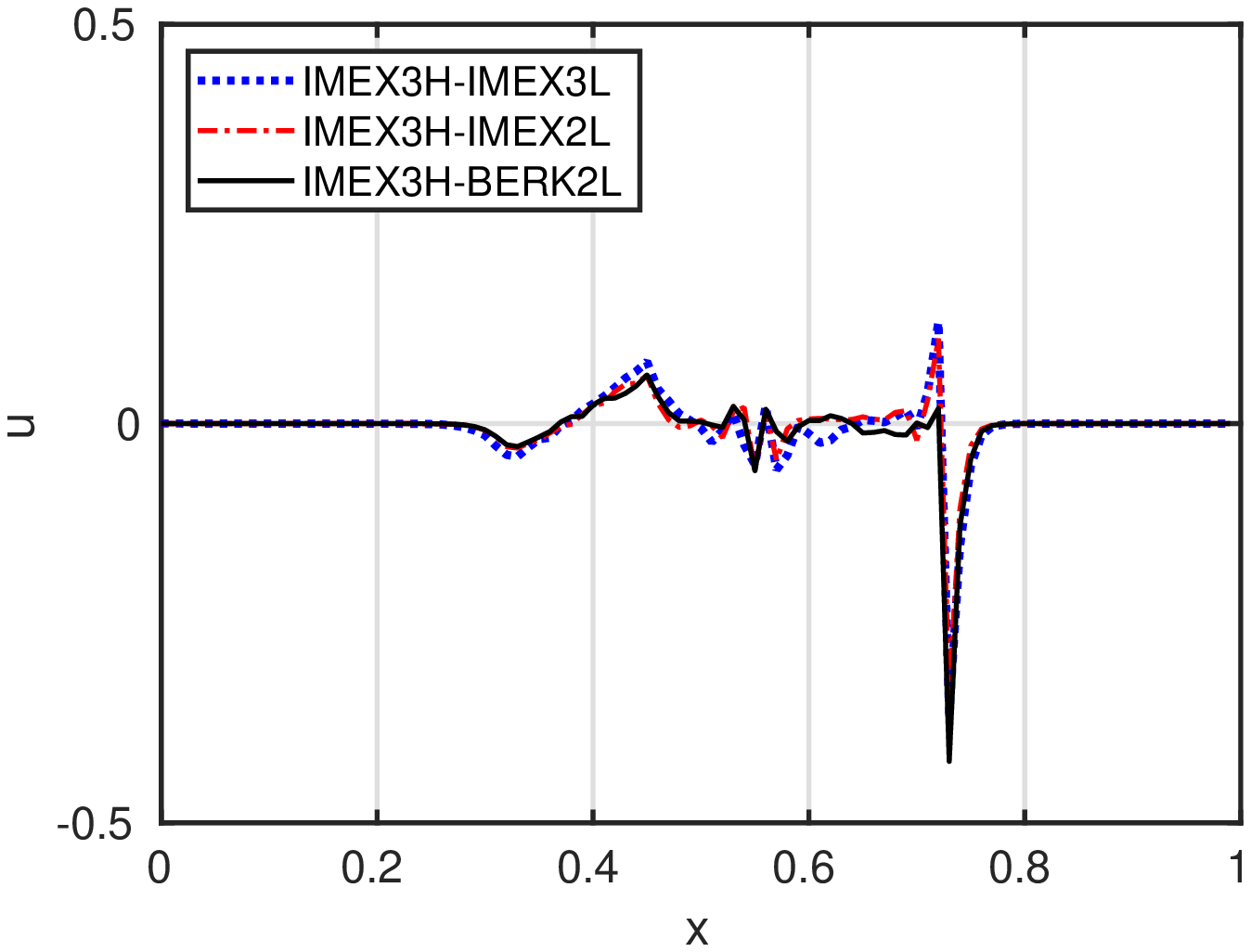}}
  \subfloat[][$\theta$]{\includegraphics[width=.25\textwidth]{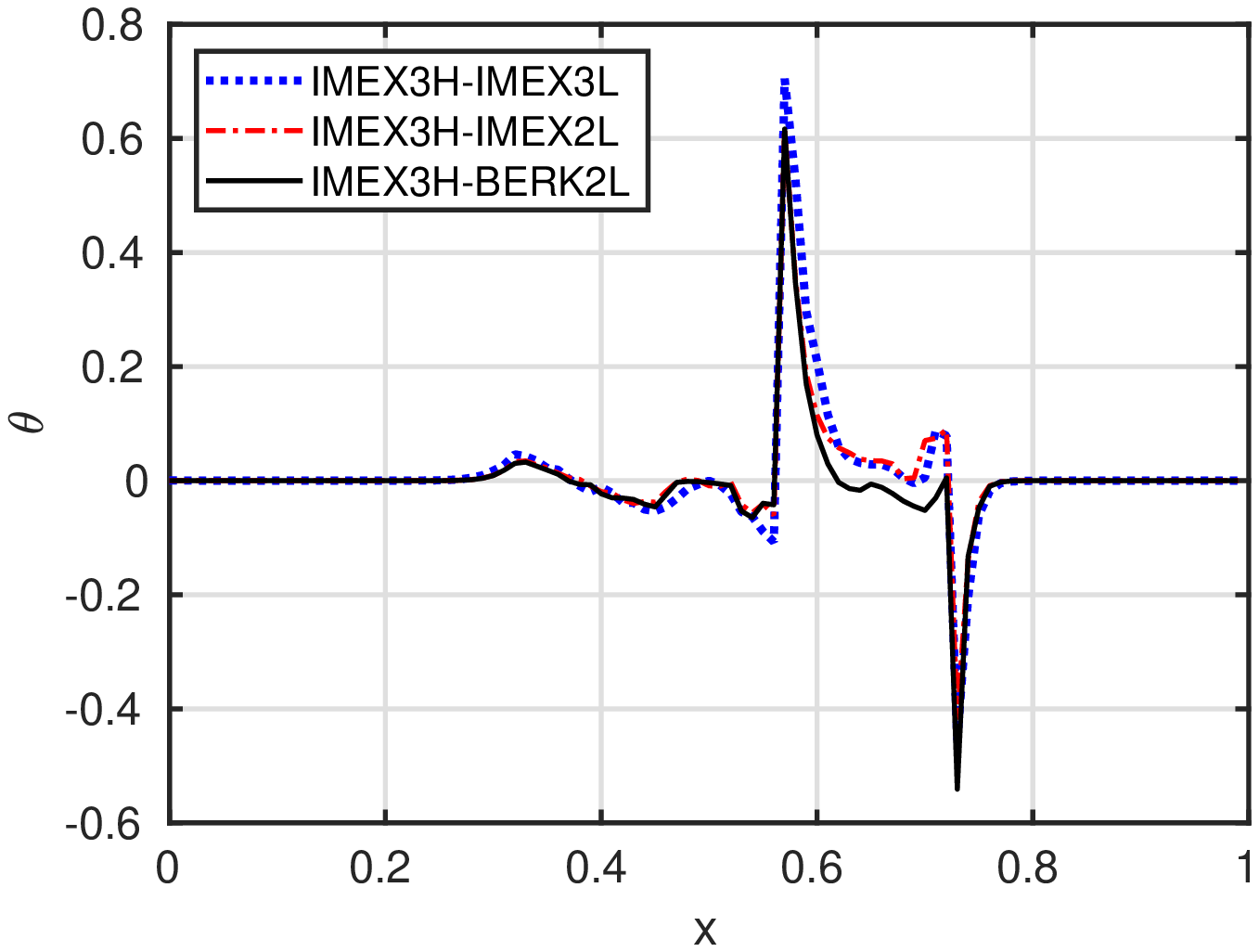}}\par
  \caption{Comparison of numerical solutions for the Sod problem at t=0.1 with various $\epsilon$. Row 1: Euler and IMEX3H at different values of $\epsilon$. Row 2: For $\epsilon =1$, and for each unknown ($\rho$, $u$, $\theta$), we plot the differences; IMEX3H-IMEX3L, IMEX3H–IMEX2L and IMEX3H–BERK2L. Row 3: Same as Row 1, but for $\epsilon=10^{-2}$, Row 4: Same as Row 3, but with a smaller CFL constant $C=0.1$. Row 5: Same as Row 1, but for $\epsilon=10^{-6}$. The computation time for each method is shown in the leftmost column. 
  (Parameters for each method - Euler: $N_x=3000$, $C=0.1$;\quad IMEX3H: $N_x=N_v=1000$, $C=0.14$;\quad IMEX3L: $N_x=N_v=100$, $C=0.14$;\quad IMEX2L: $N_x=N_v=100$, $C=0.2$;\quad BERK2L: $N_x=N_v=100$, $C=0.2$)}
  \label{fig:comparison_Sod_diff_plot}
\end{figure}
\begin{figure}[ht!]
  \centering
  \raisebox{35pt}{\parbox[b]{.15\textwidth}{IMEX3H}}%
  \subfloat[][$\rho$]{\includegraphics[width=.26\textwidth]{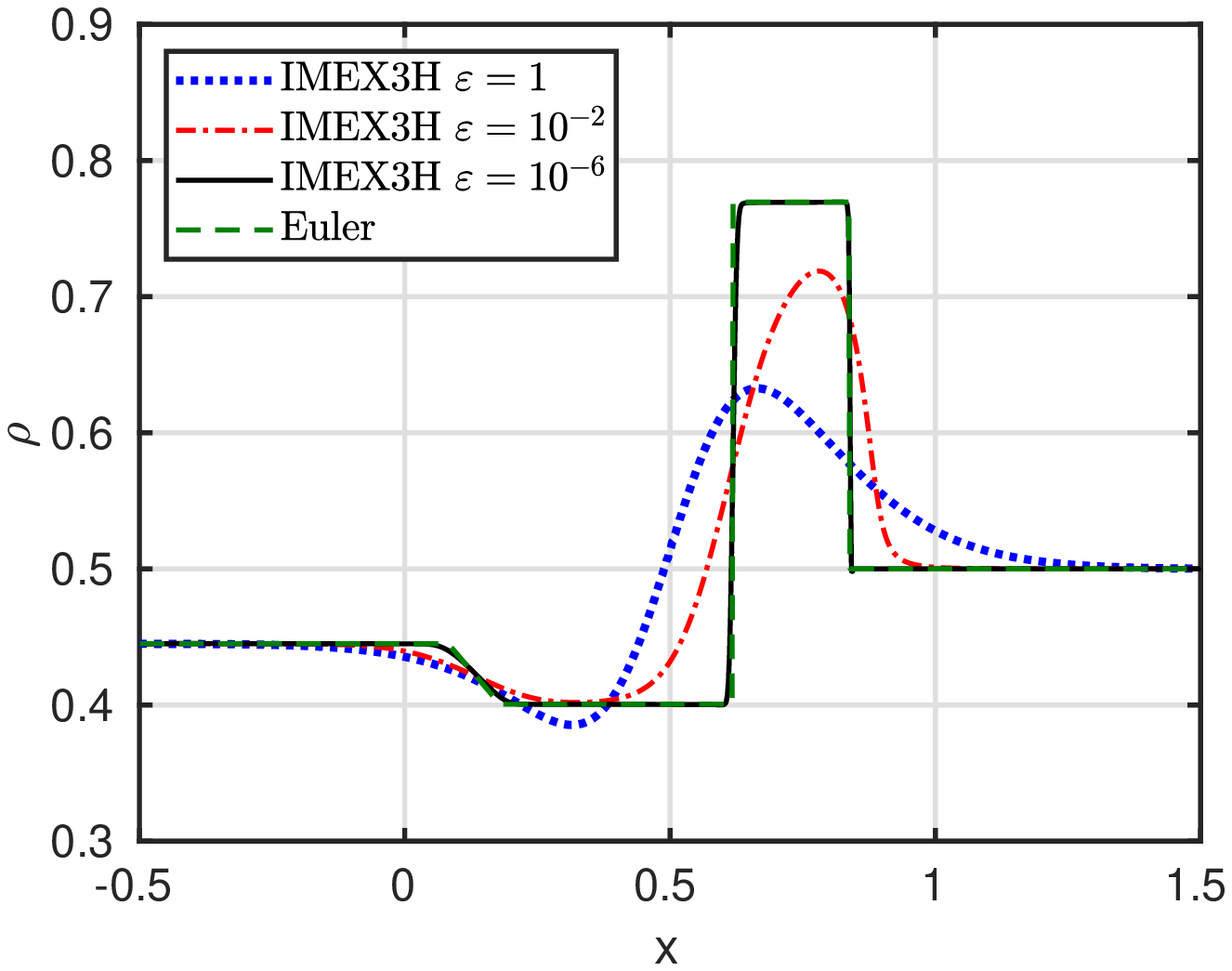}}
  \subfloat[][$u$]{\includegraphics[width=.26\textwidth]{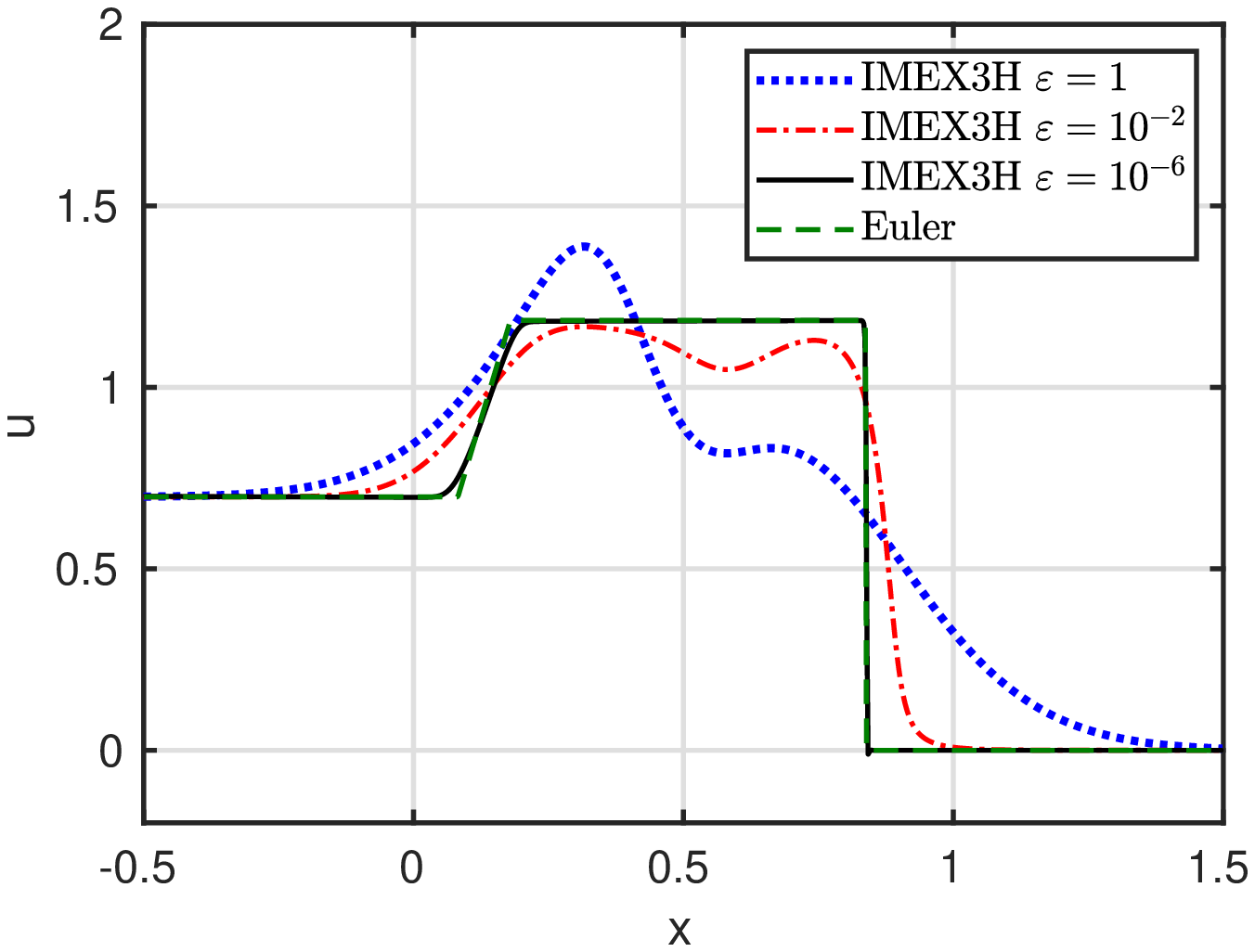}}
  \subfloat[][$\theta$]{\includegraphics[width=.26\textwidth]{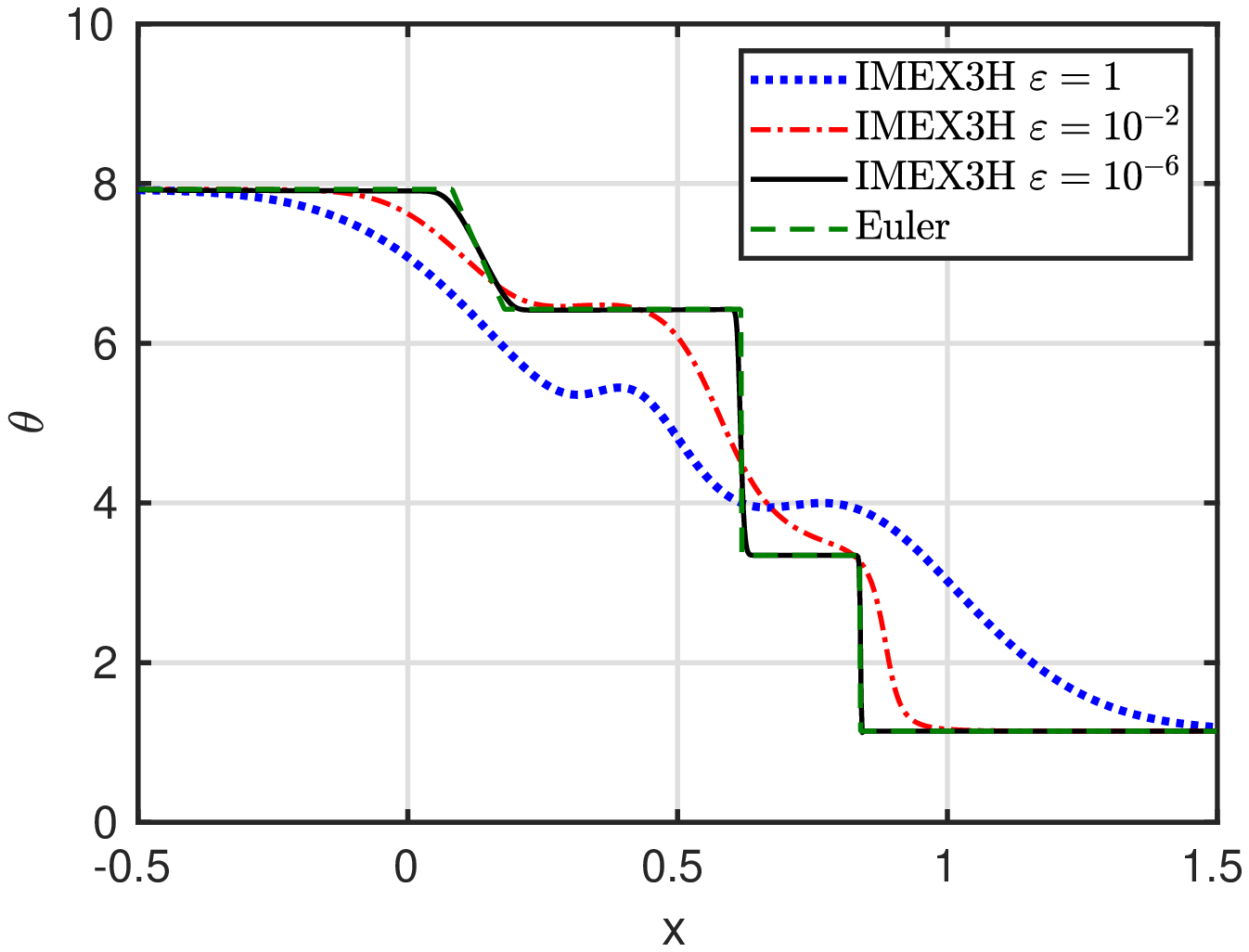}}\par  
  \raisebox{35pt}{\parbox[b]{.15\textwidth}{\tiny$\epsilon=1$\\IMEX3L (C=0.14): 25.56s\\IMEX2L (C=0.2): 11.91s\\BERK2L (C=0.2): 1.74s}}%
  \subfloat[][$\rho$]{\includegraphics[width=.26\textwidth]{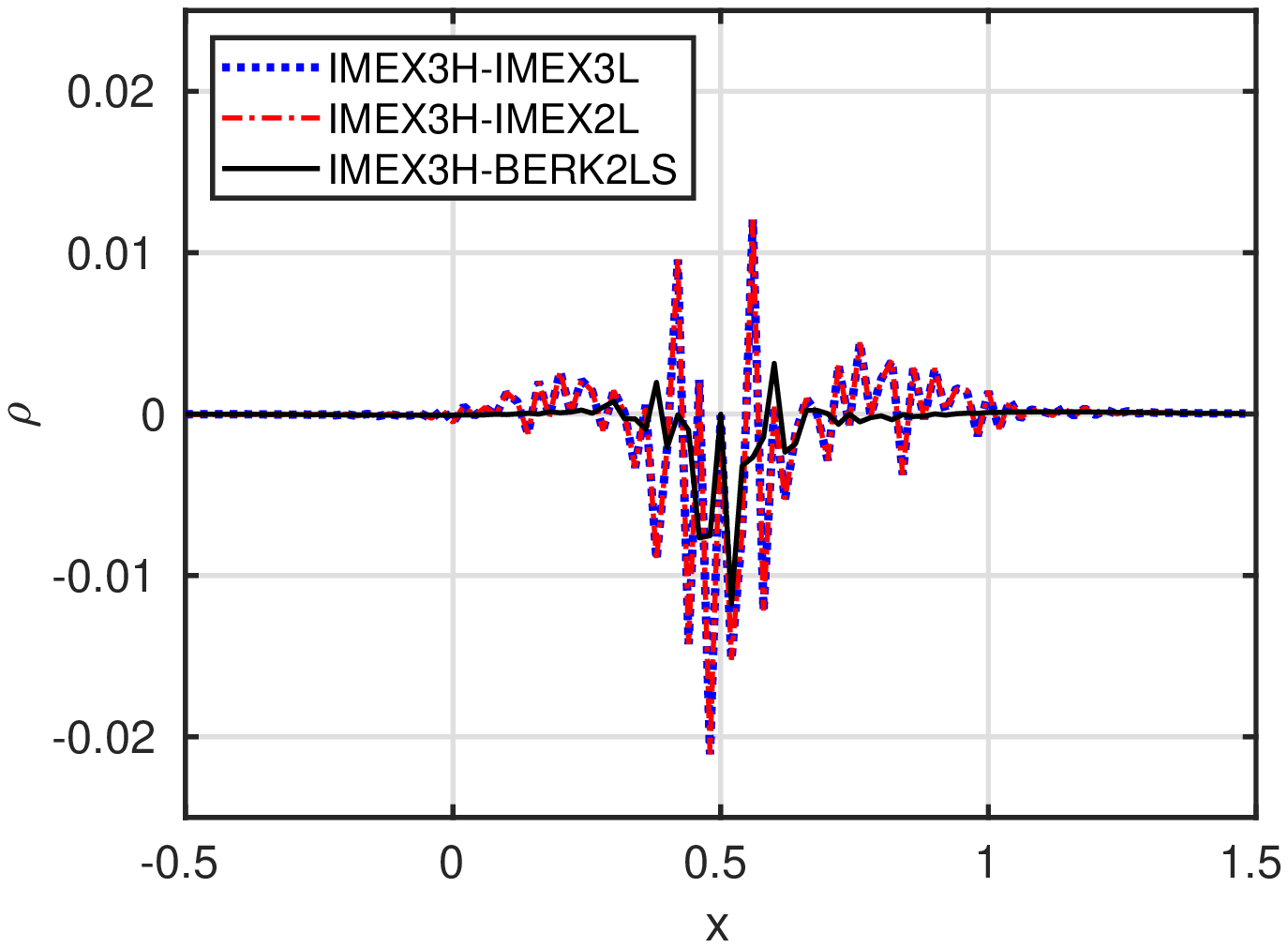}}
  \subfloat[][$u$]{\includegraphics[width=.26\textwidth]{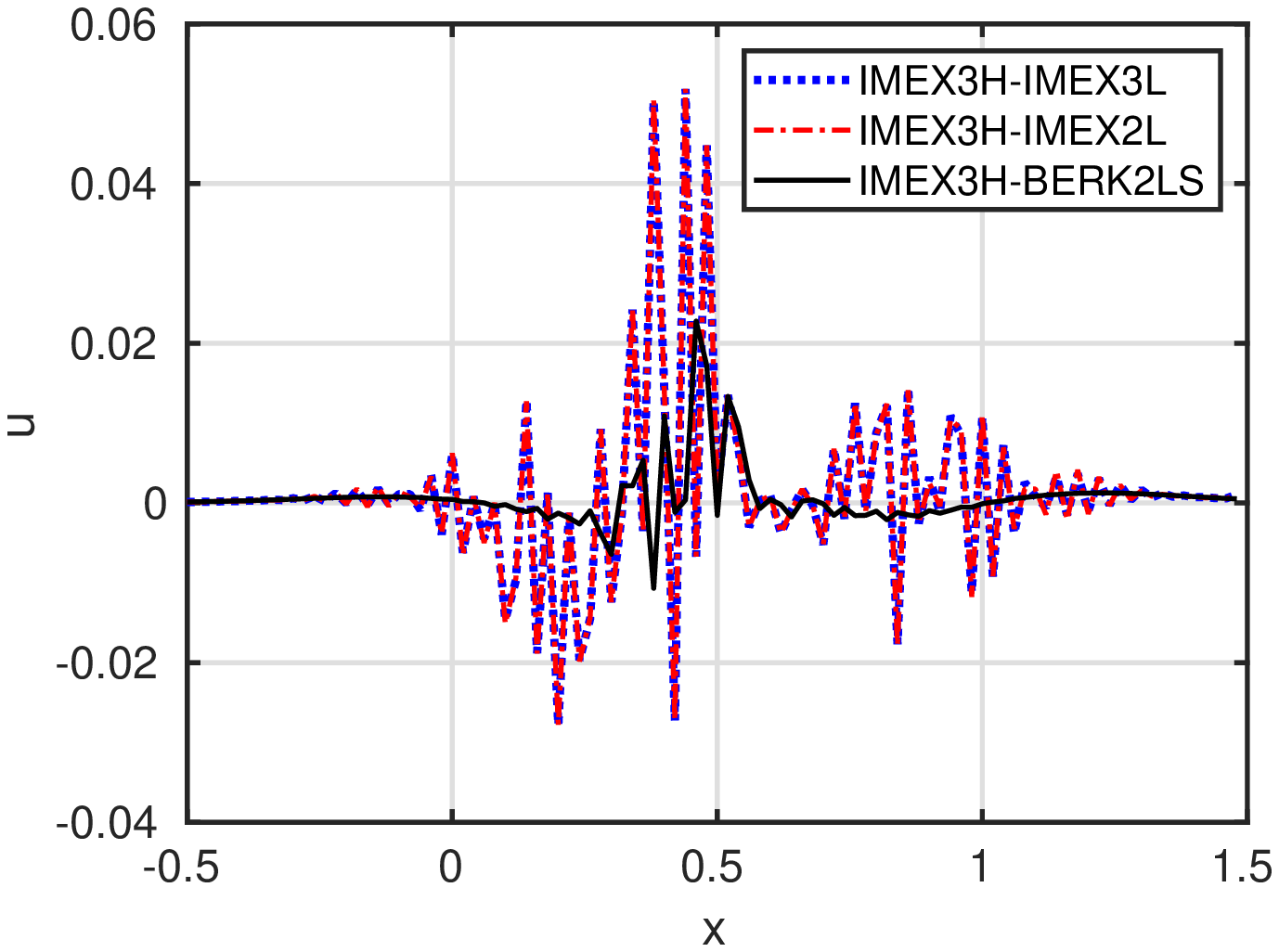}}
  \subfloat[][$\theta$]{\includegraphics[width=.26\textwidth]{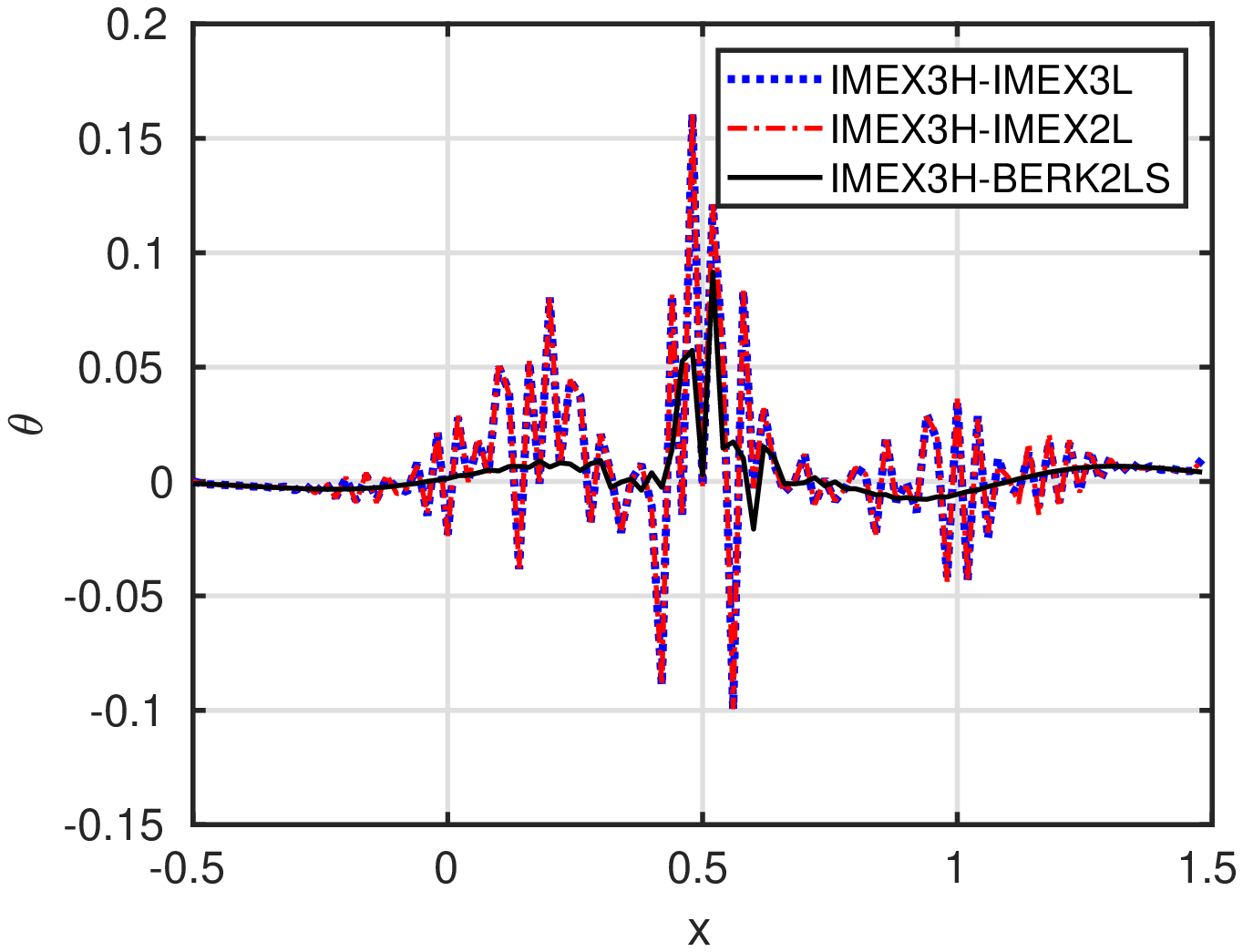}}\par
  \raisebox{35pt}{\parbox[b]{.15\textwidth}{\tiny$\epsilon=10^{-2}$\\IMEX3L (C=0.14): 24.99s\\IMEX2L (C=0.2): 12.16s\\BERK2L (C=0.2): 1.74s}}%
  \subfloat[][$\rho$]{\includegraphics[width=.26\textwidth]{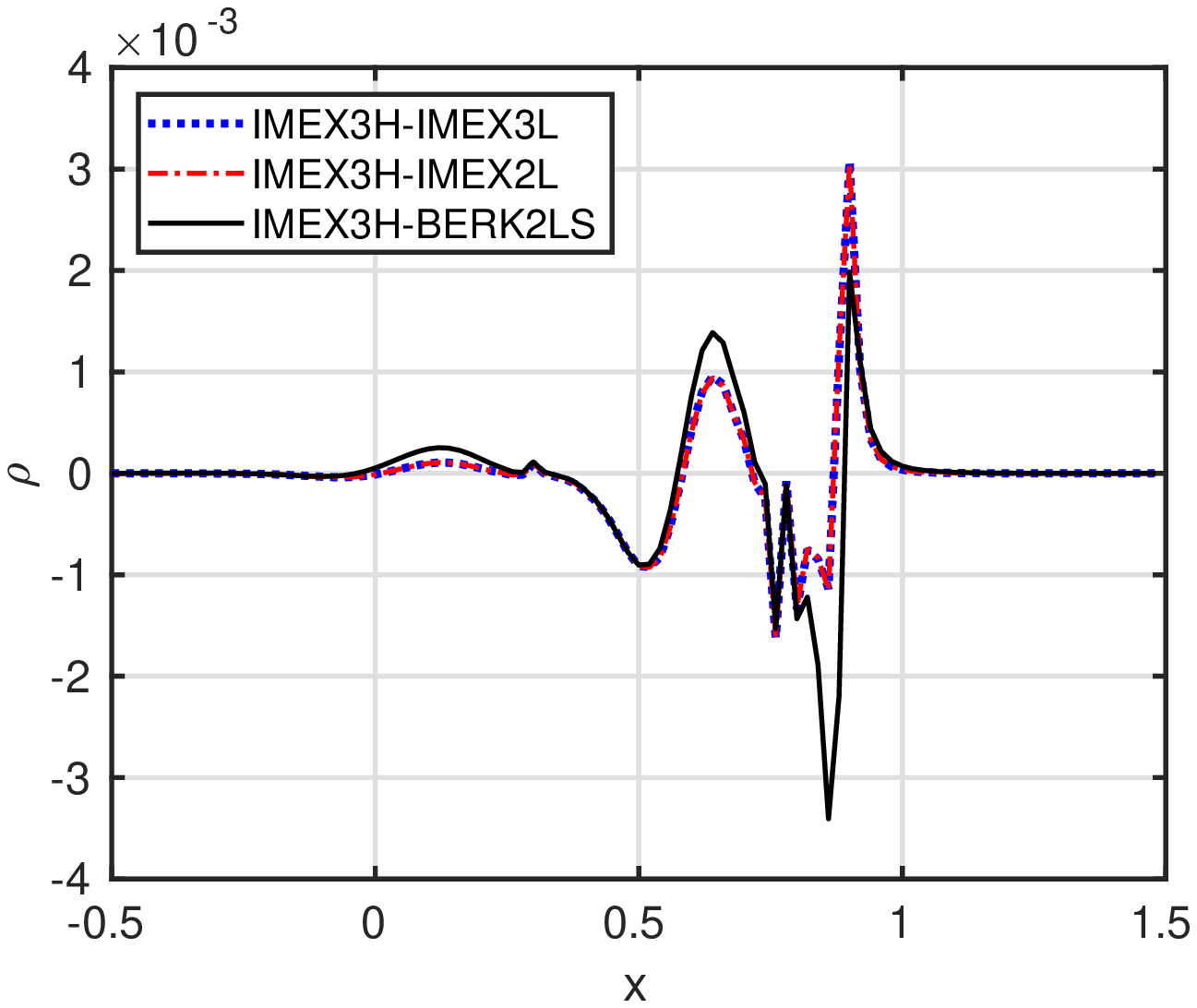}}
  \subfloat[][$u$]{\includegraphics[width=.26\textwidth]{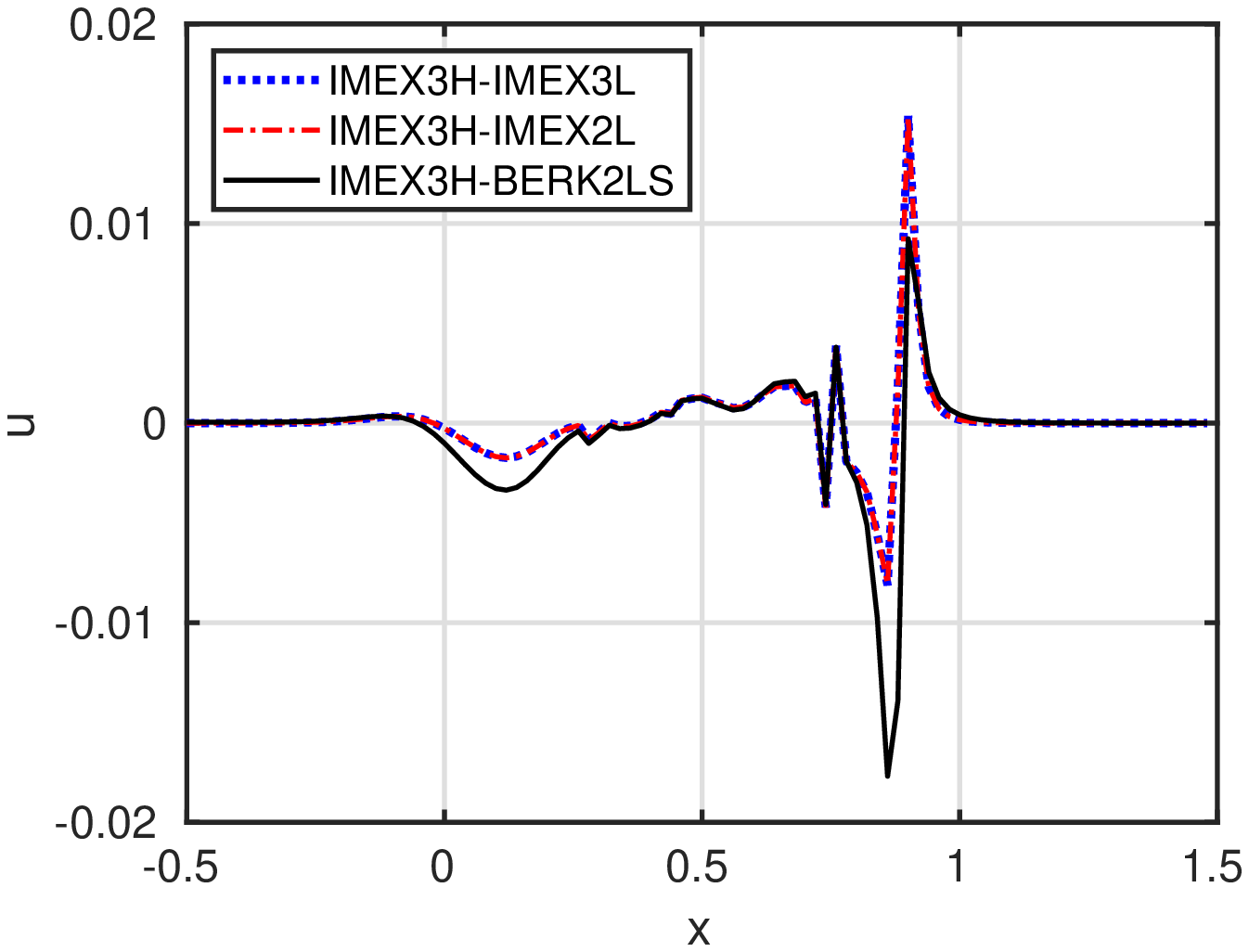}}
  \subfloat[][$\theta$]{\includegraphics[width=.26\textwidth]{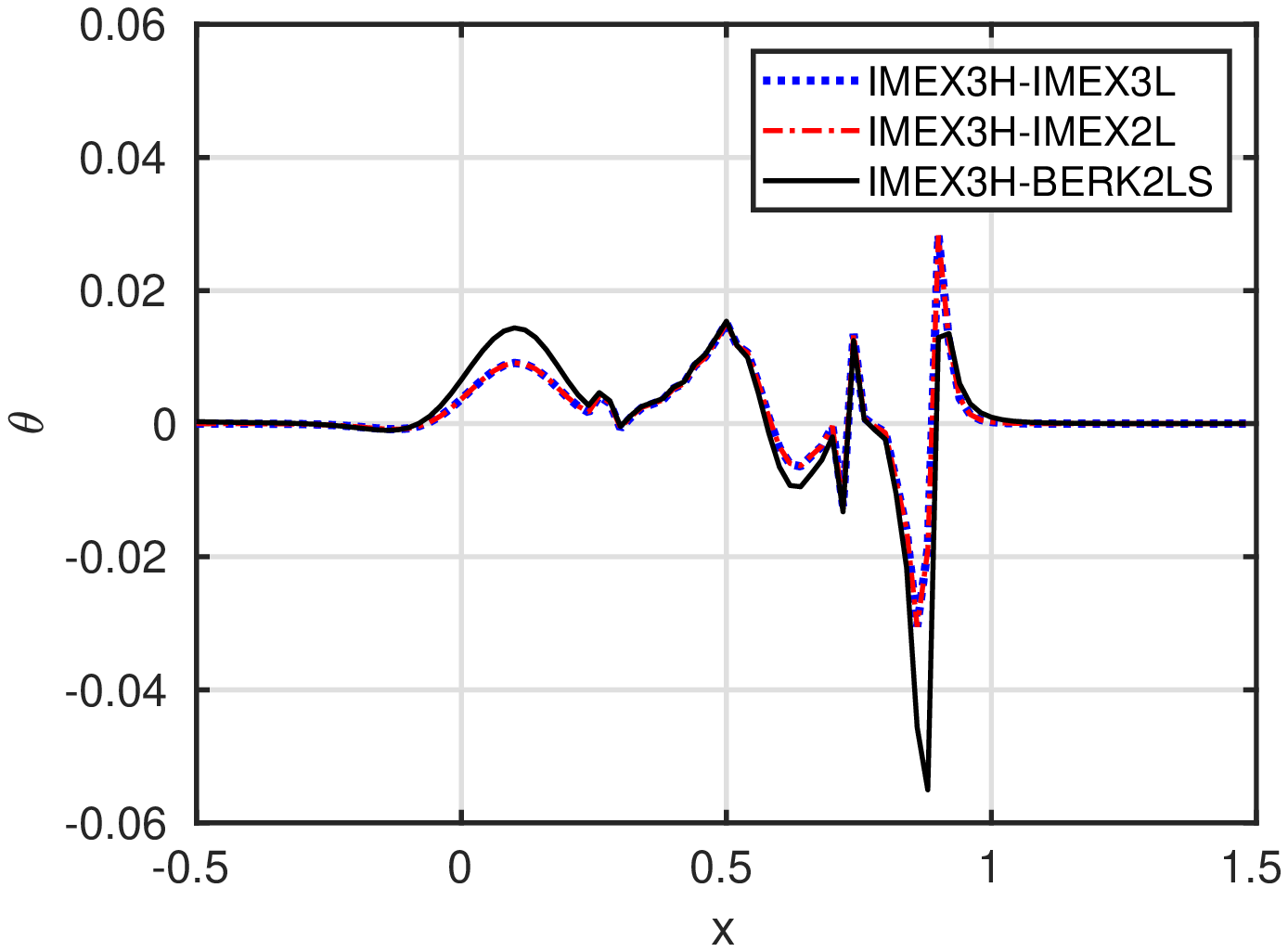}}\par
  \raisebox{35pt}{\parbox[b]{.15\textwidth}{\tiny$\epsilon=10^{-2}$\\IMEX3L (C=0.14): 24.99s\\IMEX2L (C=0.2): 12.16s\\BERK2LS (C=0.1): 3.77s}}%
  \subfloat[][$\rho$]{\includegraphics[width=.26\textwidth]{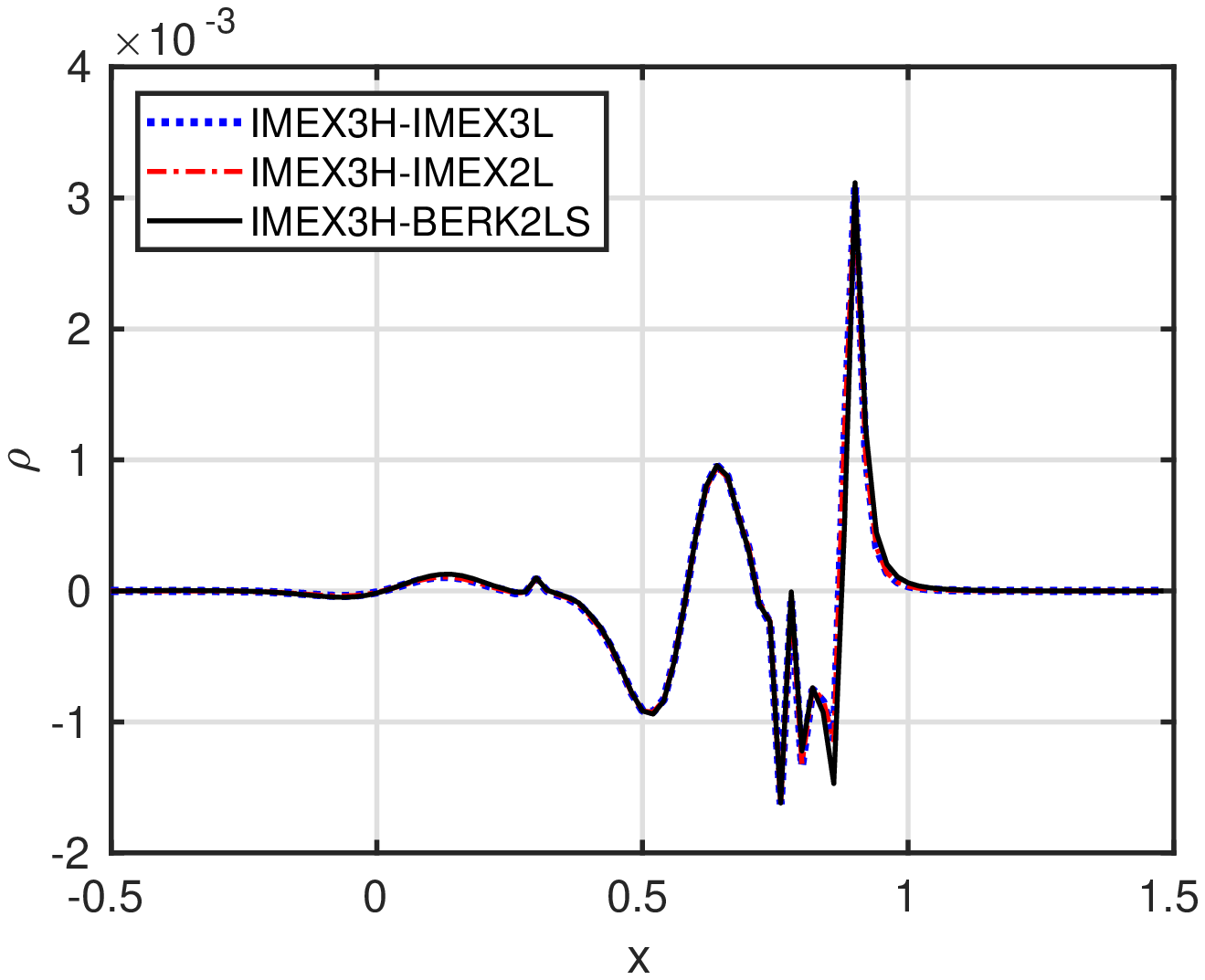}}
  \subfloat[][$u$]{\includegraphics[width=.26\textwidth]{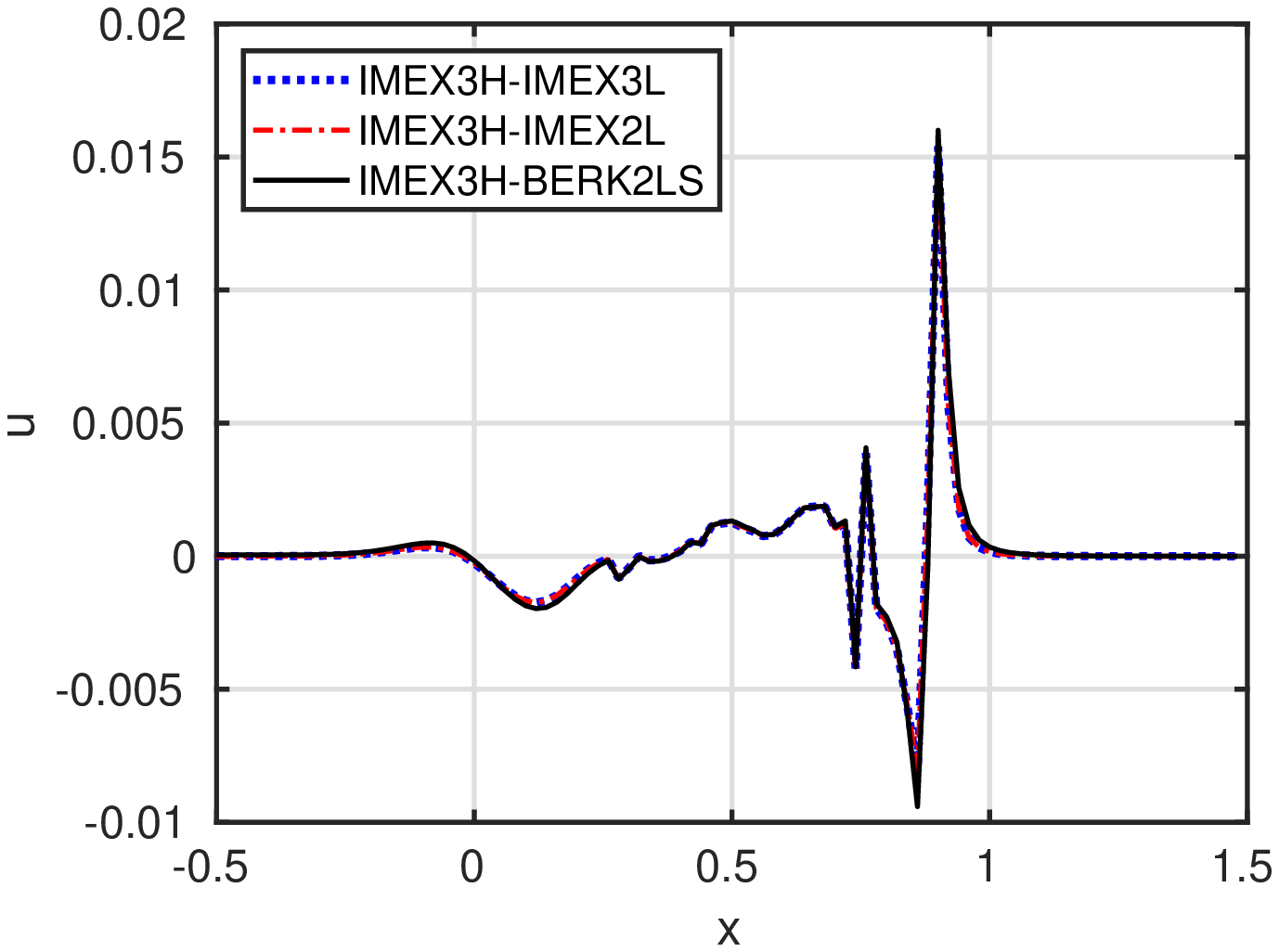}}
  \subfloat[][$\theta$]{\includegraphics[width=.26\textwidth]{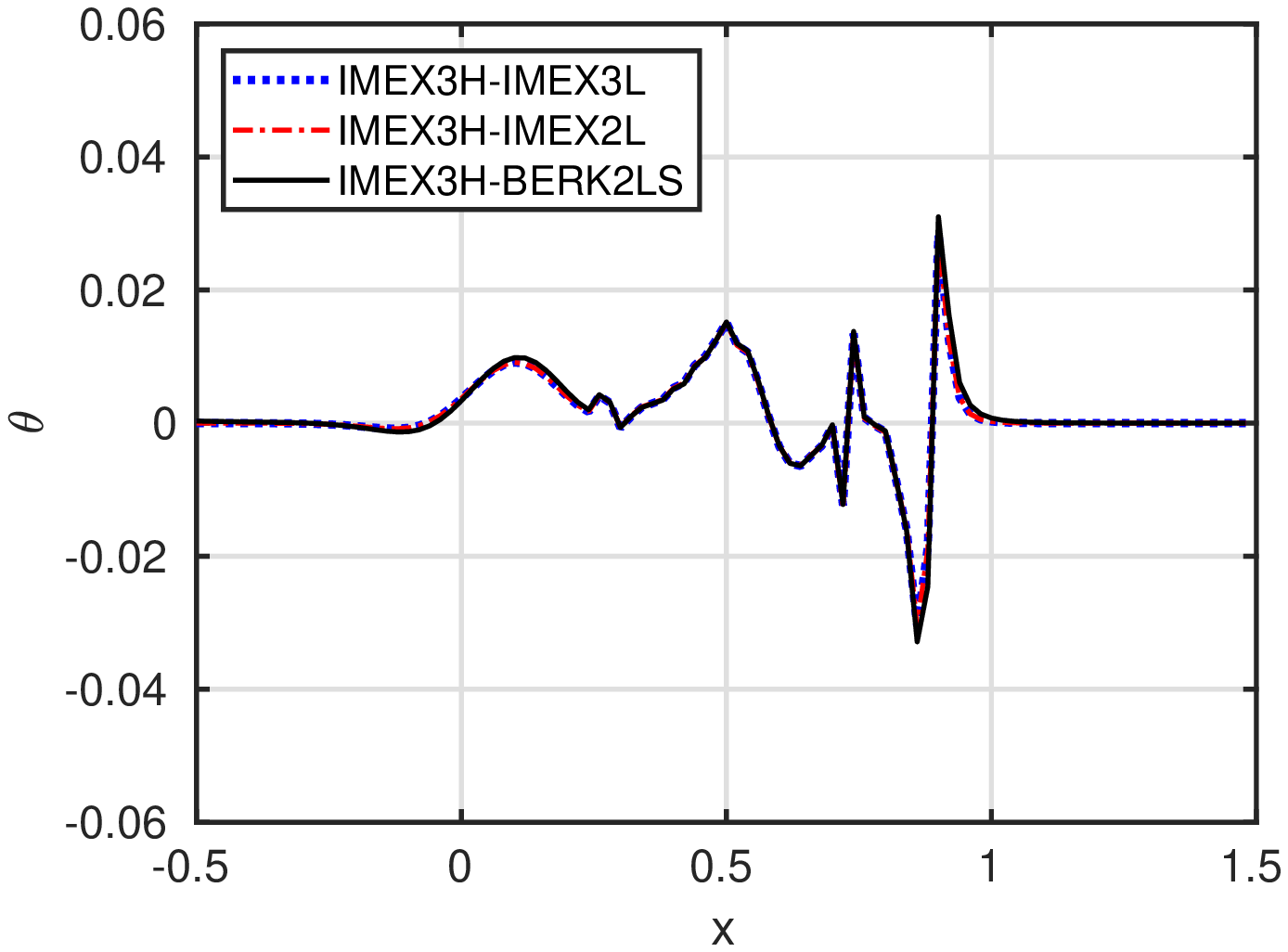}}\par
  \raisebox{35pt}{\parbox[b]{.15\textwidth}{\tiny$\epsilon=10^{-6}$\\IMEX3L (C=0.14): 25.21s\\IMEX2L (C=0.2): 11.37s\\BERK2L (C=0.2): 1.71s}}%
  \subfloat[][$\rho$]{\includegraphics[width=.26\textwidth]{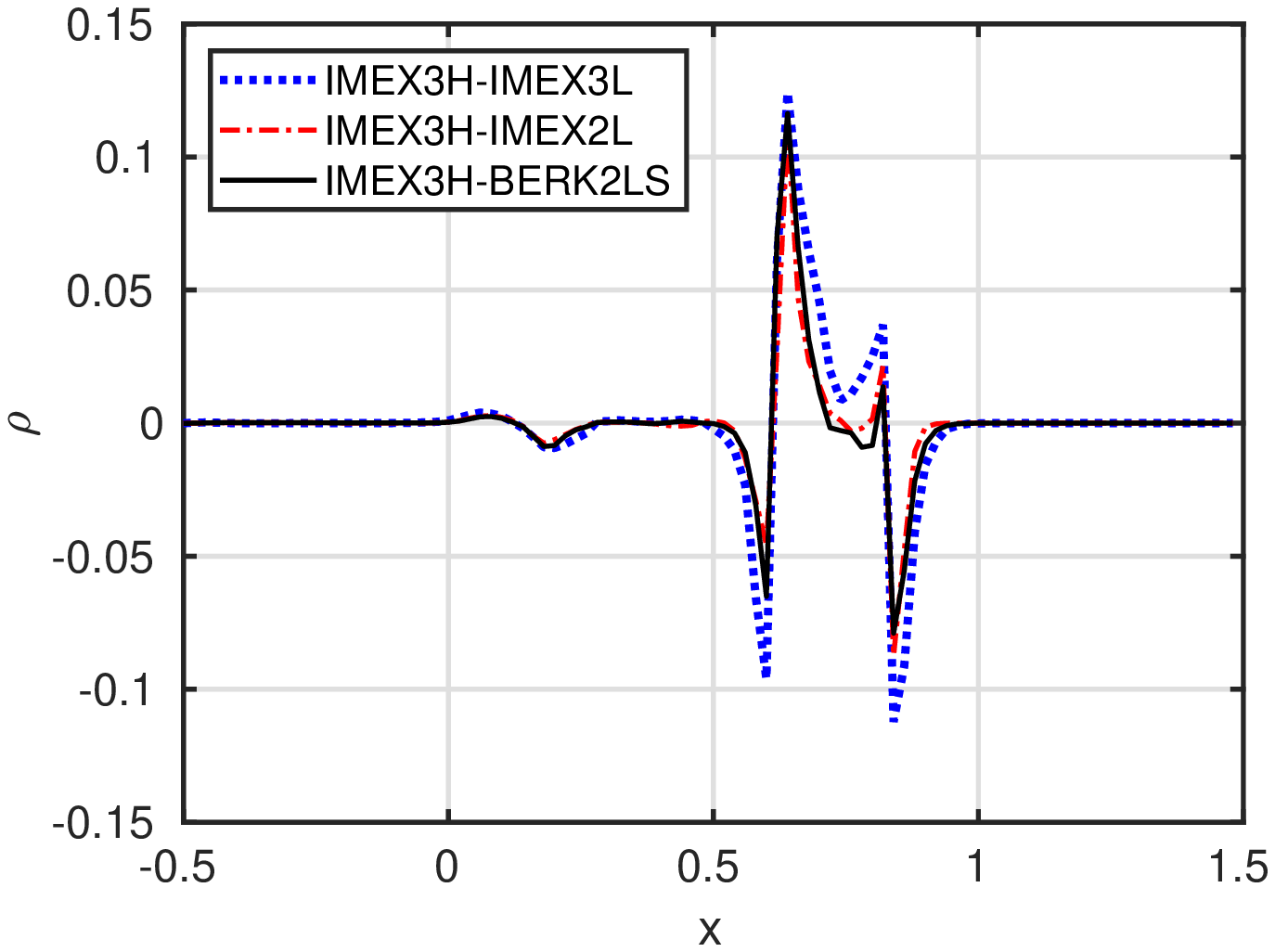}}
  \subfloat[][$u$]{\includegraphics[width=.26\textwidth]{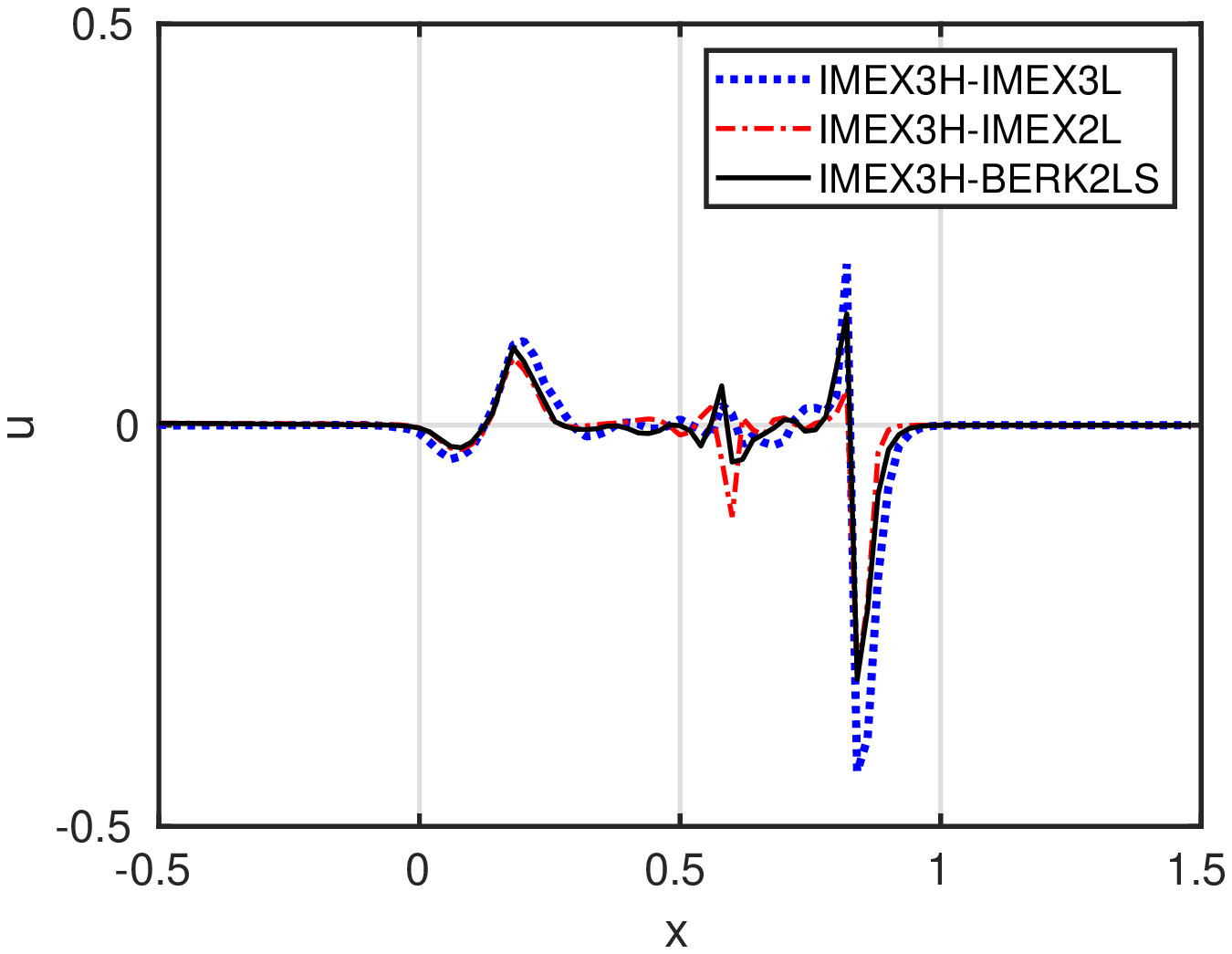}}
  \subfloat[][$\theta$]{\includegraphics[width=.26\textwidth]{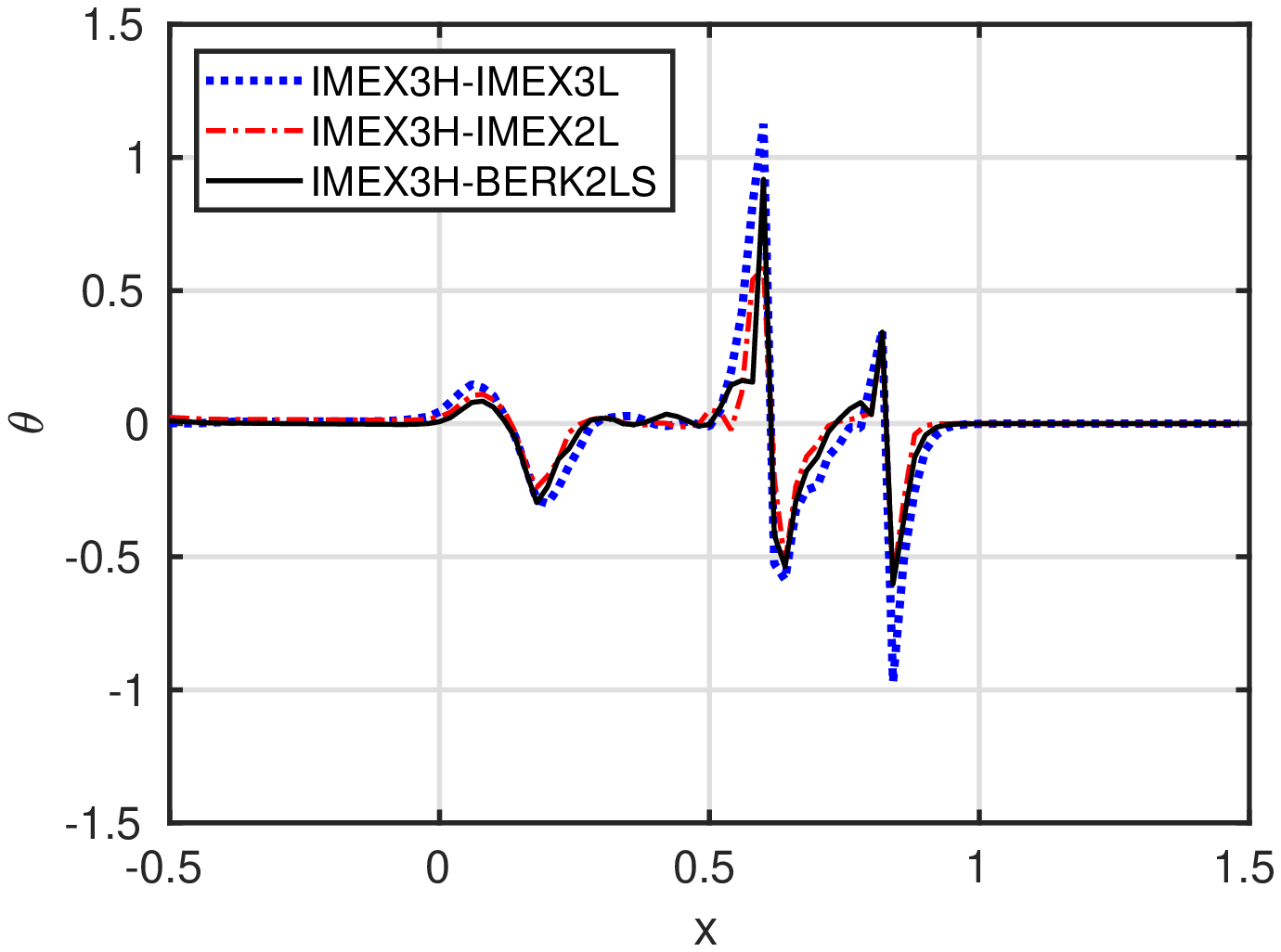}}\par
  \caption{Comparison of numerical solutions for the Lax problem at t=0.1 with various $\epsilon$. Row 1: Euler and IMEX3H at different values of $\epsilon$. Row 2: For $\epsilon =1$, and for each unknown ($\rho$, $u$, $\theta$), we plot the differences; IMEX3H-IMEX3L, IMEX3H–IMEX2L and IMEX3H–BERK2L. Row 3: Same as Row 1, but for $\epsilon=10^{-2}$, Row 4: Same as Row 1, but for $\epsilon=10^{-6}$. The computation time for each method is shown in the first column. Row 5: Comparisons of low-resolution solution for $\epsilon=0.01$. (Parameters for each method - Euler: $N_x=3000$, $C=0.1$;\quad IMEX3H: $N_x=N_v=1000$, $C=0.14$;\quad IMEX3L: $N_x=N_v=100$, $C=0.14$;\quad IMEX2L: $N_x=N_v=100$, $C=0.2$;\quad BERK2L: $N_x=N_v=100$, $C=0.2$)}
  \label{fig:comparison_Lax_diff_plot}
\end{figure}
\begin{figure}[ht!]
  \centering
  \raisebox{35pt}{\parbox[b]{.16\textwidth}{IMEX3H}}%
  \subfloat[][$\rho$]{\includegraphics[width=.27\textwidth]{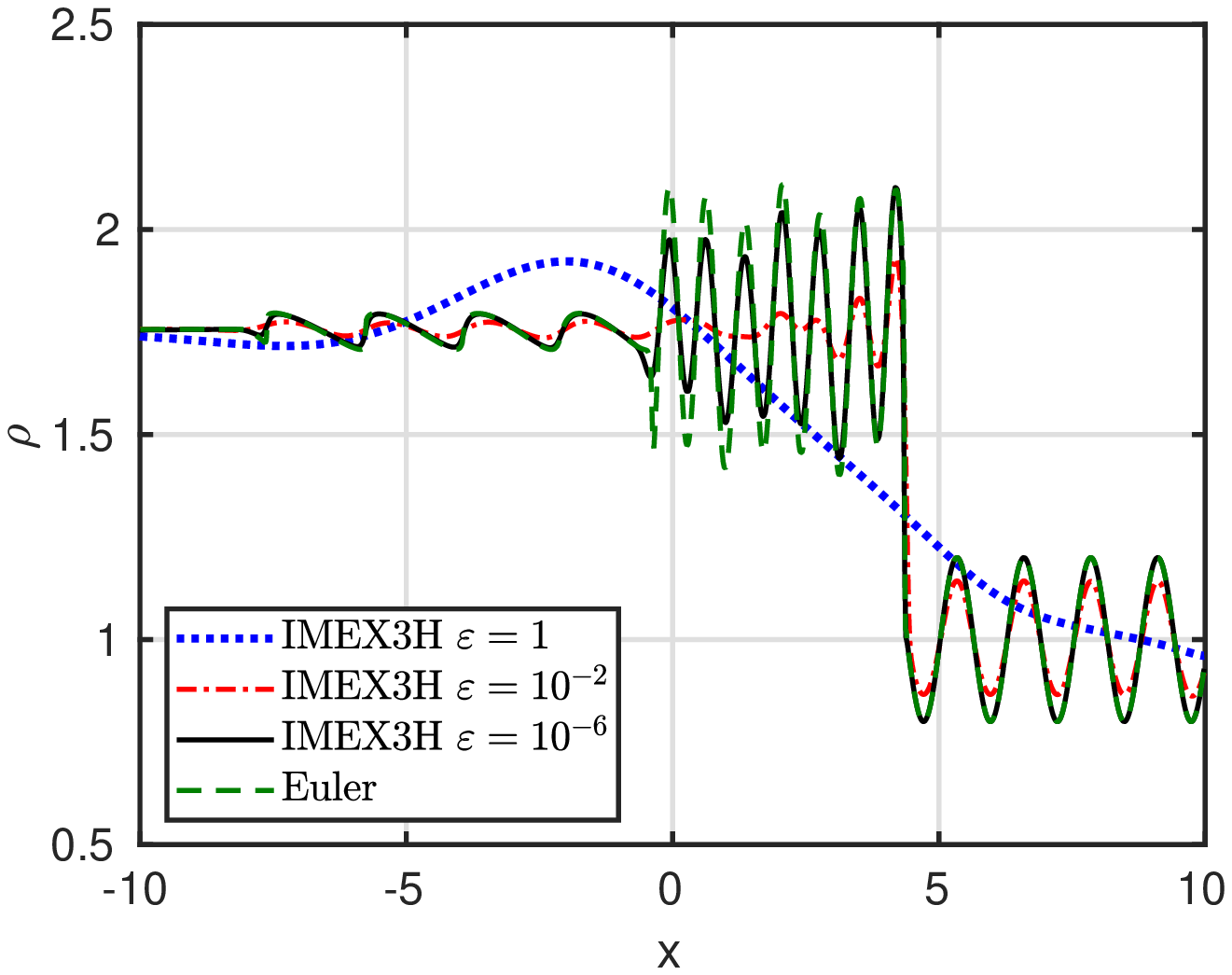}}
  \subfloat[][$u$]{\includegraphics[width=.27\textwidth]{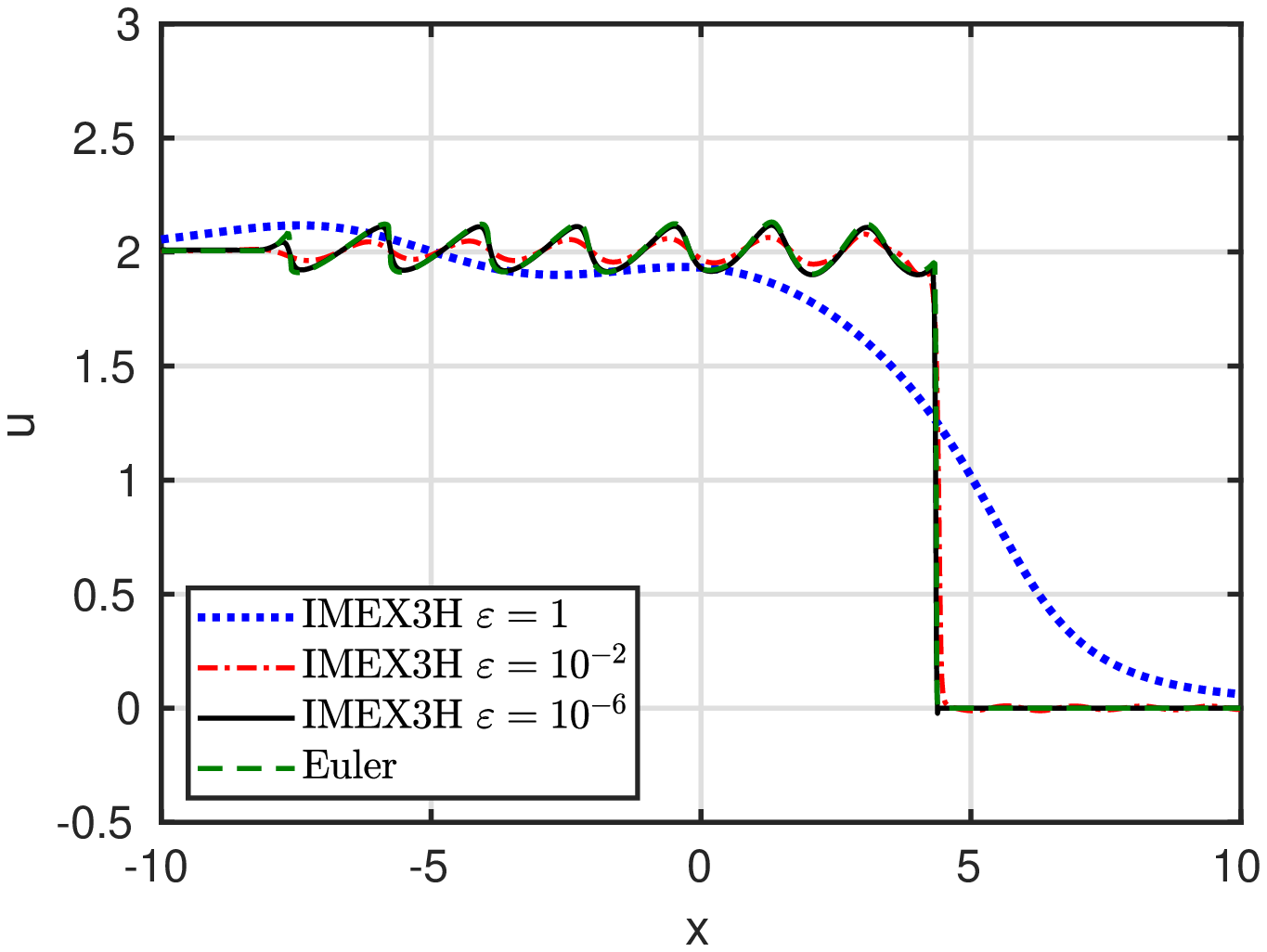}}
  \subfloat[][$\theta$]{\includegraphics[width=.27\textwidth]{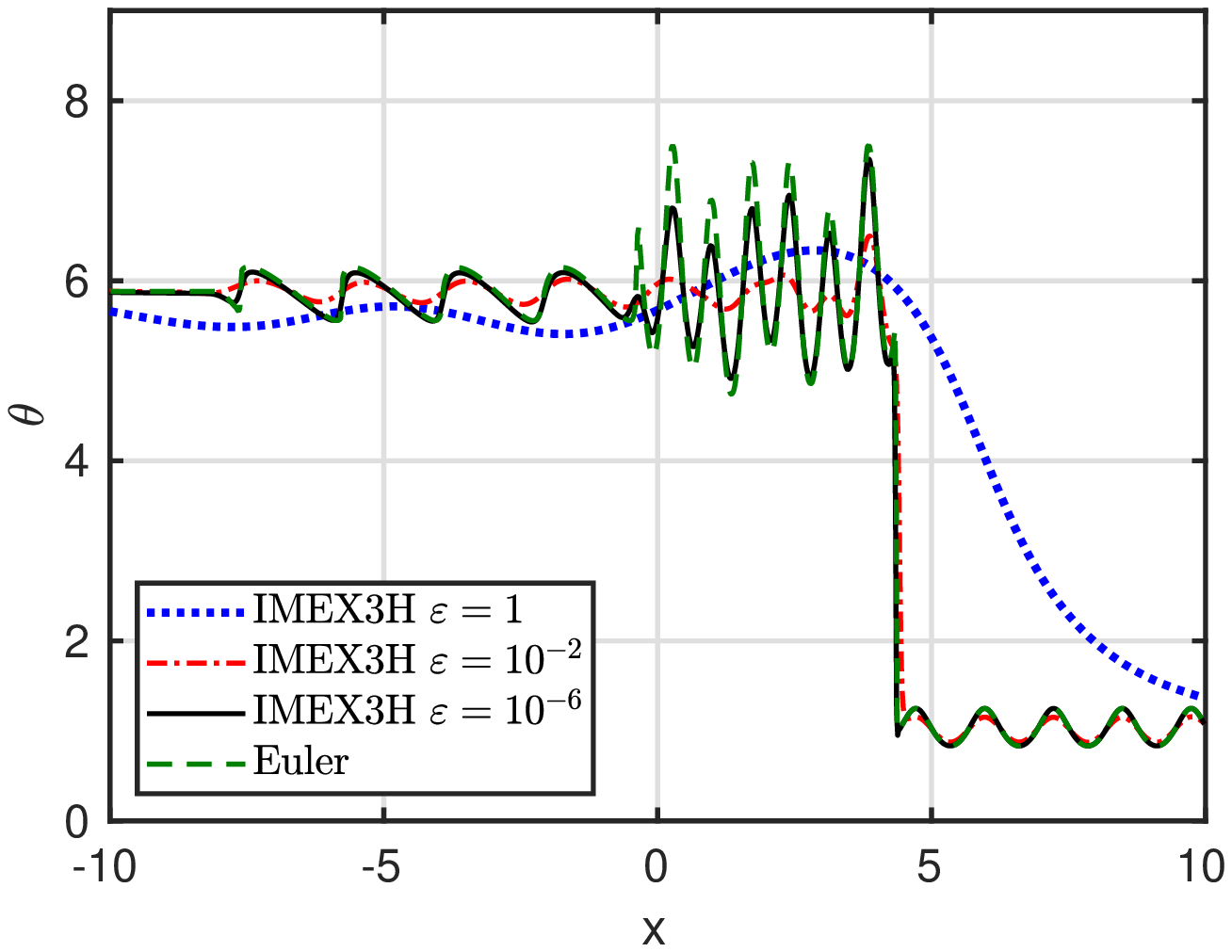}}\par  
  \raisebox{35pt}{\parbox[b]{.16\textwidth}{\tiny$\epsilon=1$\\IMEX3L (C=0.14): 292.33s\\IMEX2L (C=0.2): 120.93s\\BERK2L (C=0.2): 29.77s}}%
  \subfloat[][$\rho$]{\includegraphics[width=.27\textwidth]{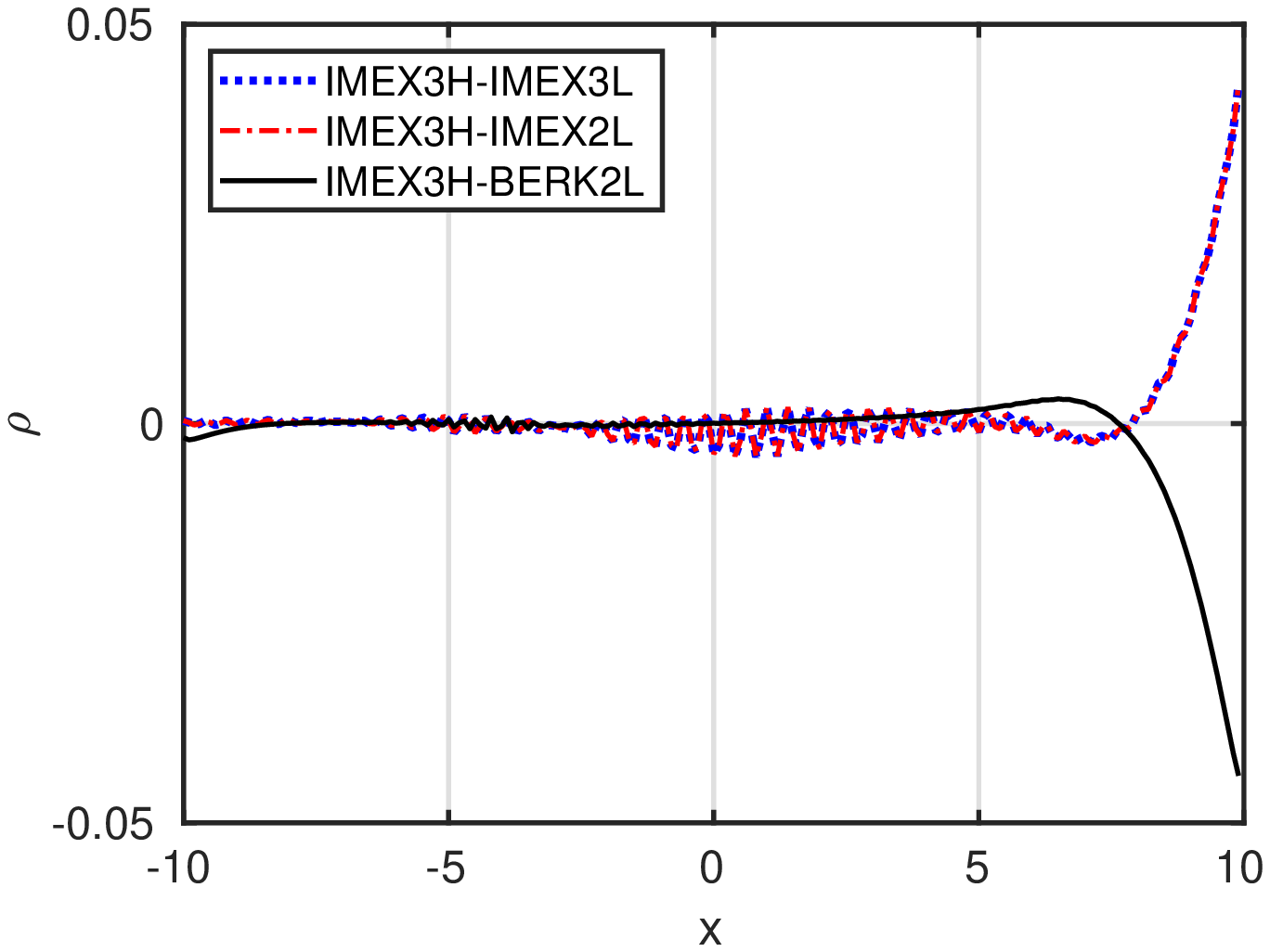}}
  \subfloat[][$u$]{\includegraphics[width=.27\textwidth]{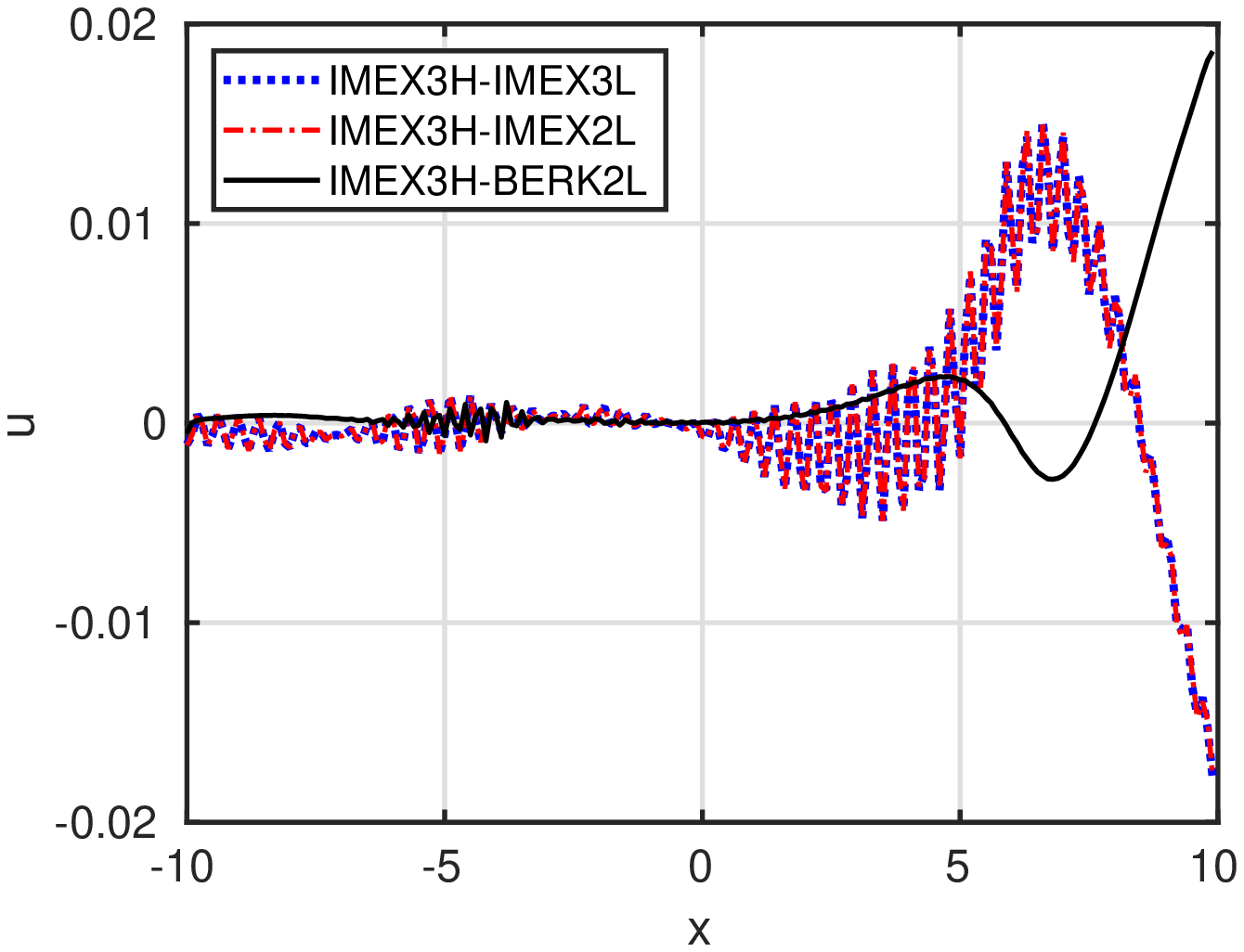}}
  \subfloat[][$\theta$]{\includegraphics[width=.27\textwidth]{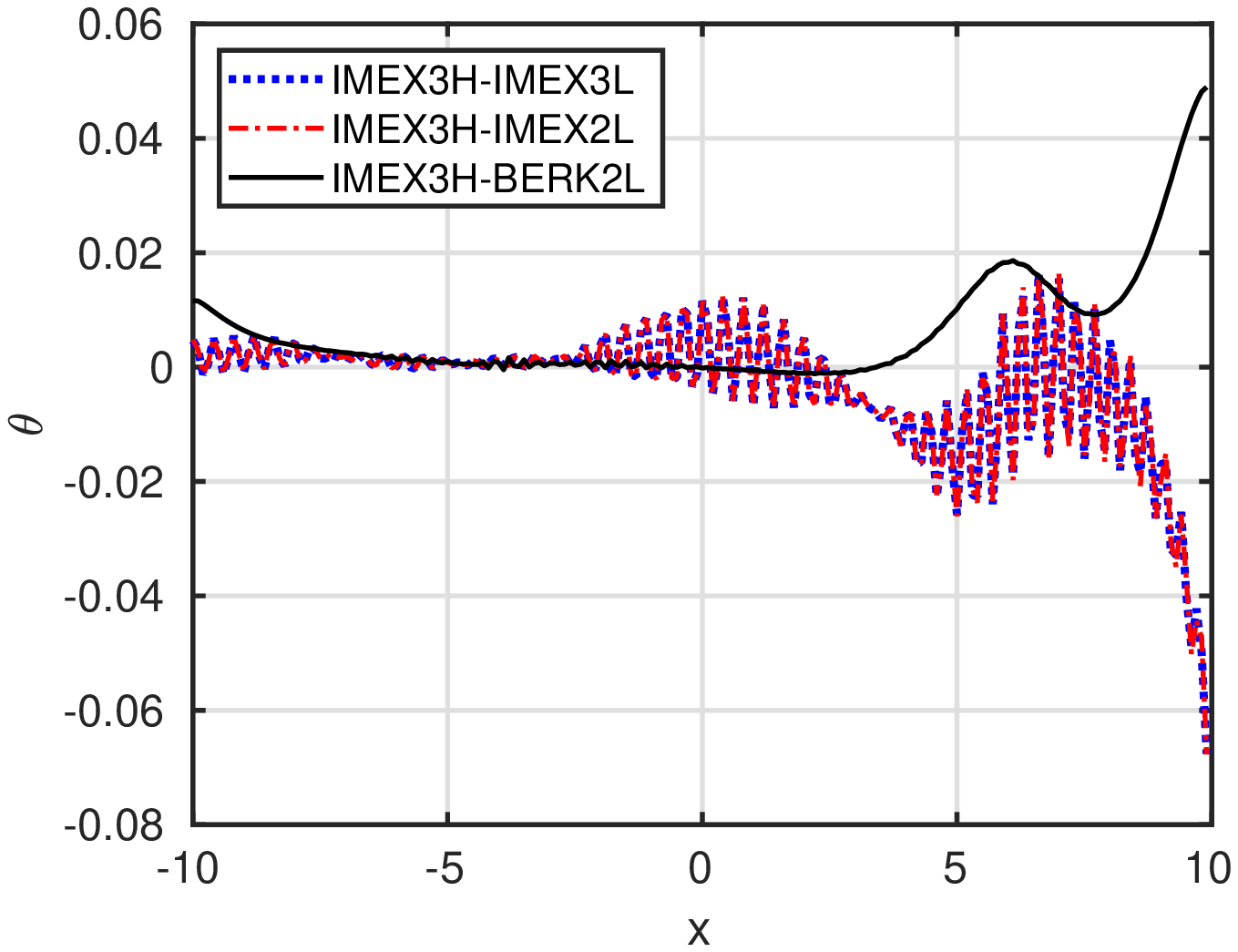}}\par
  \raisebox{35pt}{\parbox[b]{.16\textwidth}{\tiny$\epsilon=10^{-2}$\\IMEX3L (C=0.14): 287.11s\\IMEX2L (C=0.2): 121.02s\\BERK2L (C=0.2): 29.50s}}%
  \subfloat[][$\rho$]{\includegraphics[width=.27\textwidth]{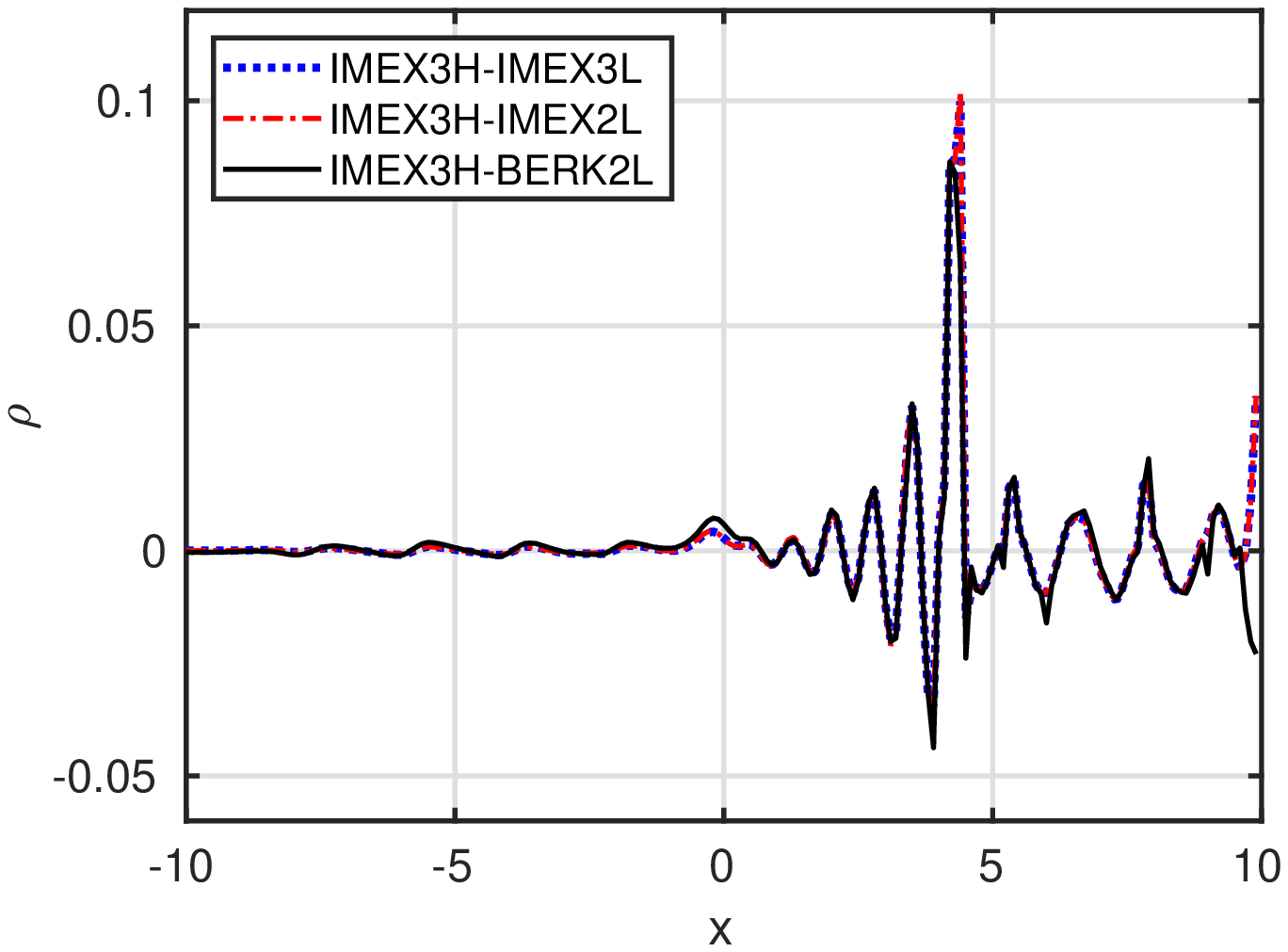}}
  \subfloat[][$u$]{\includegraphics[width=.27\textwidth]{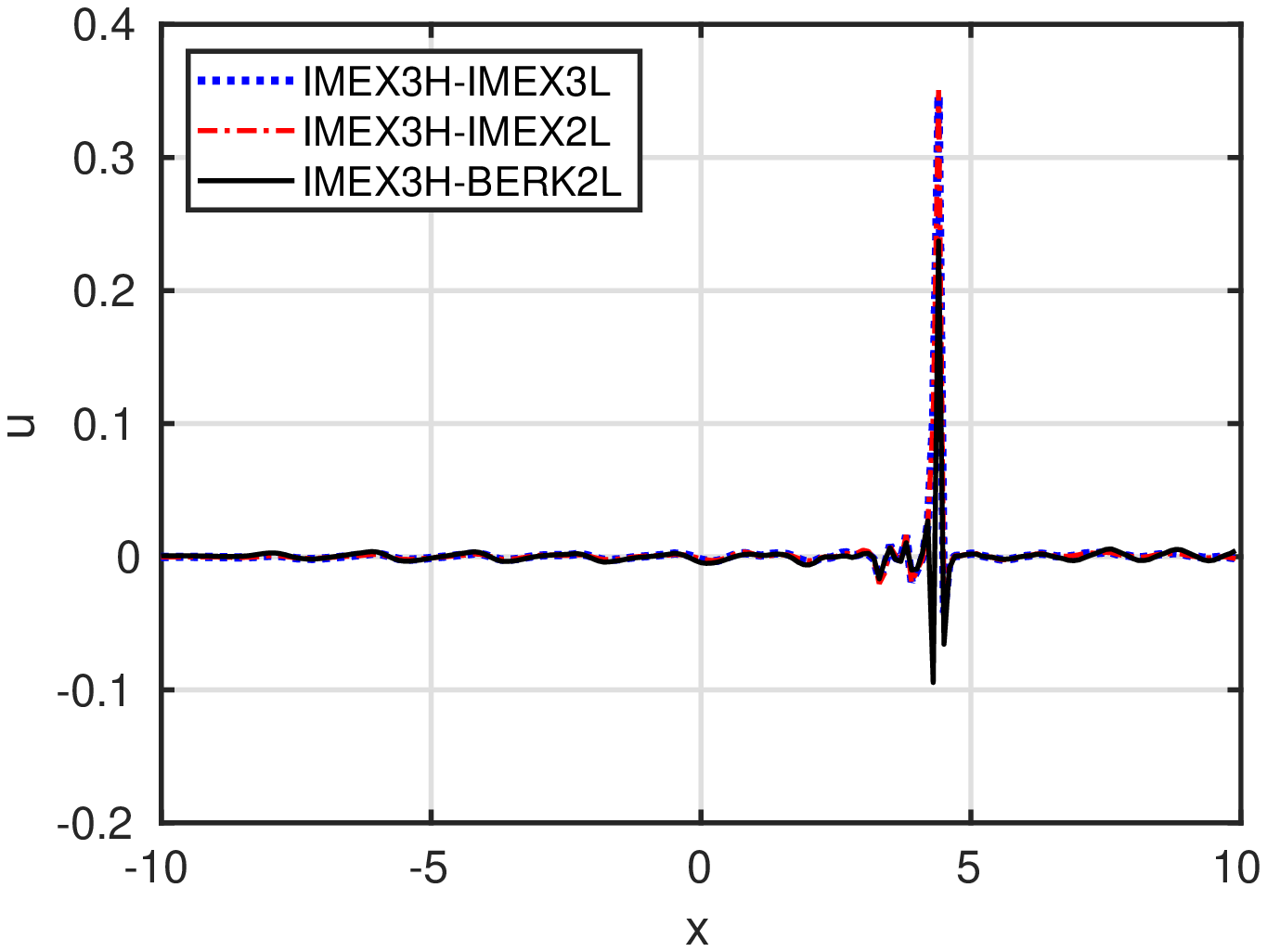}}
  \subfloat[][$\theta$]{\includegraphics[width=.27\textwidth]{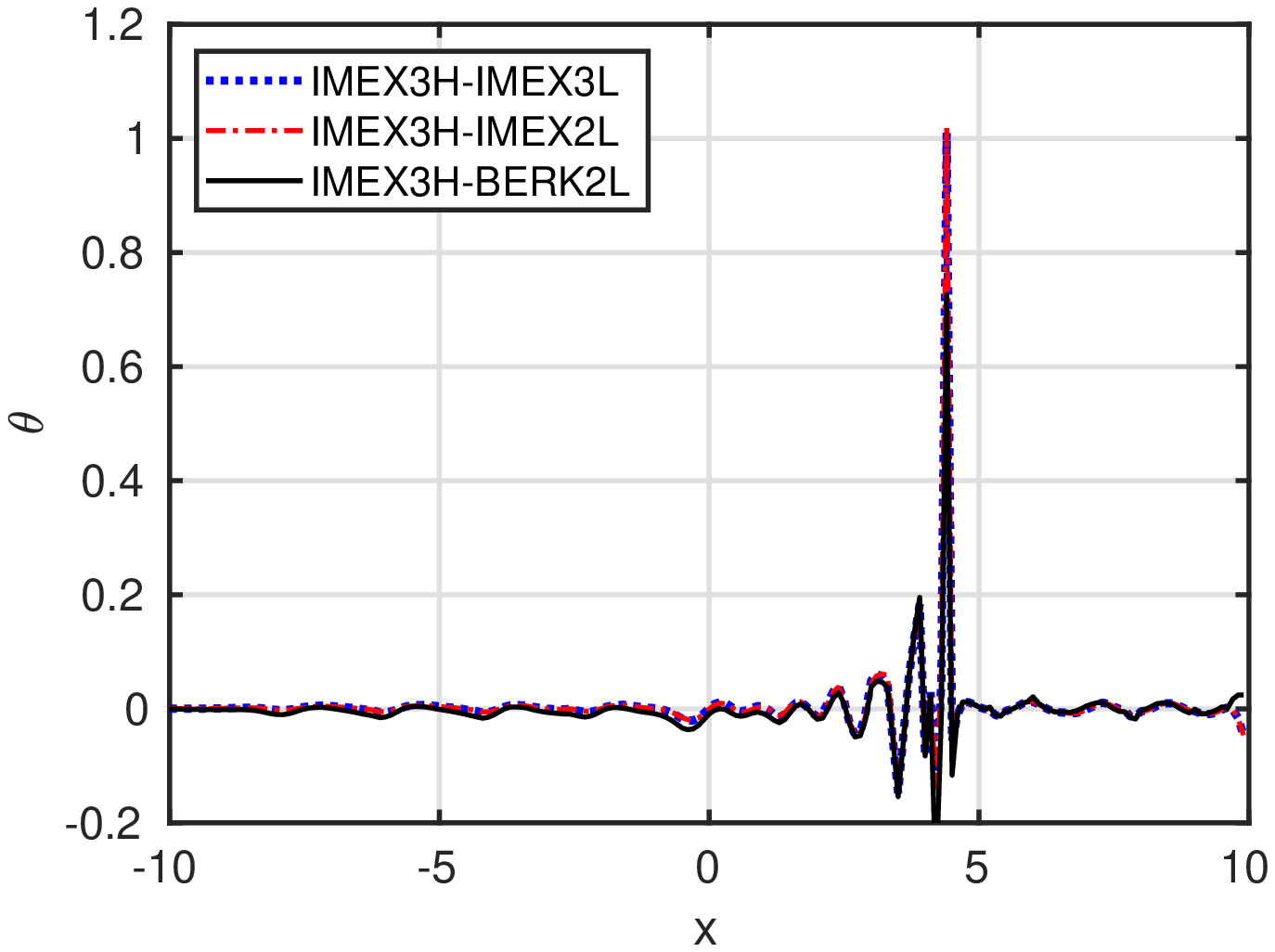}}\par
  \raisebox{35pt}{\parbox[b]{.16\textwidth}{\tiny$\epsilon=10^{-2}$\\IMEX3L (C=0.14): 287.11s\\IMEX2L (C=0.2): 121.02s\\BERK2LS (C=0.1): 61.20s}}%
  \subfloat[][$\rho$]{\includegraphics[width=.27\textwidth]{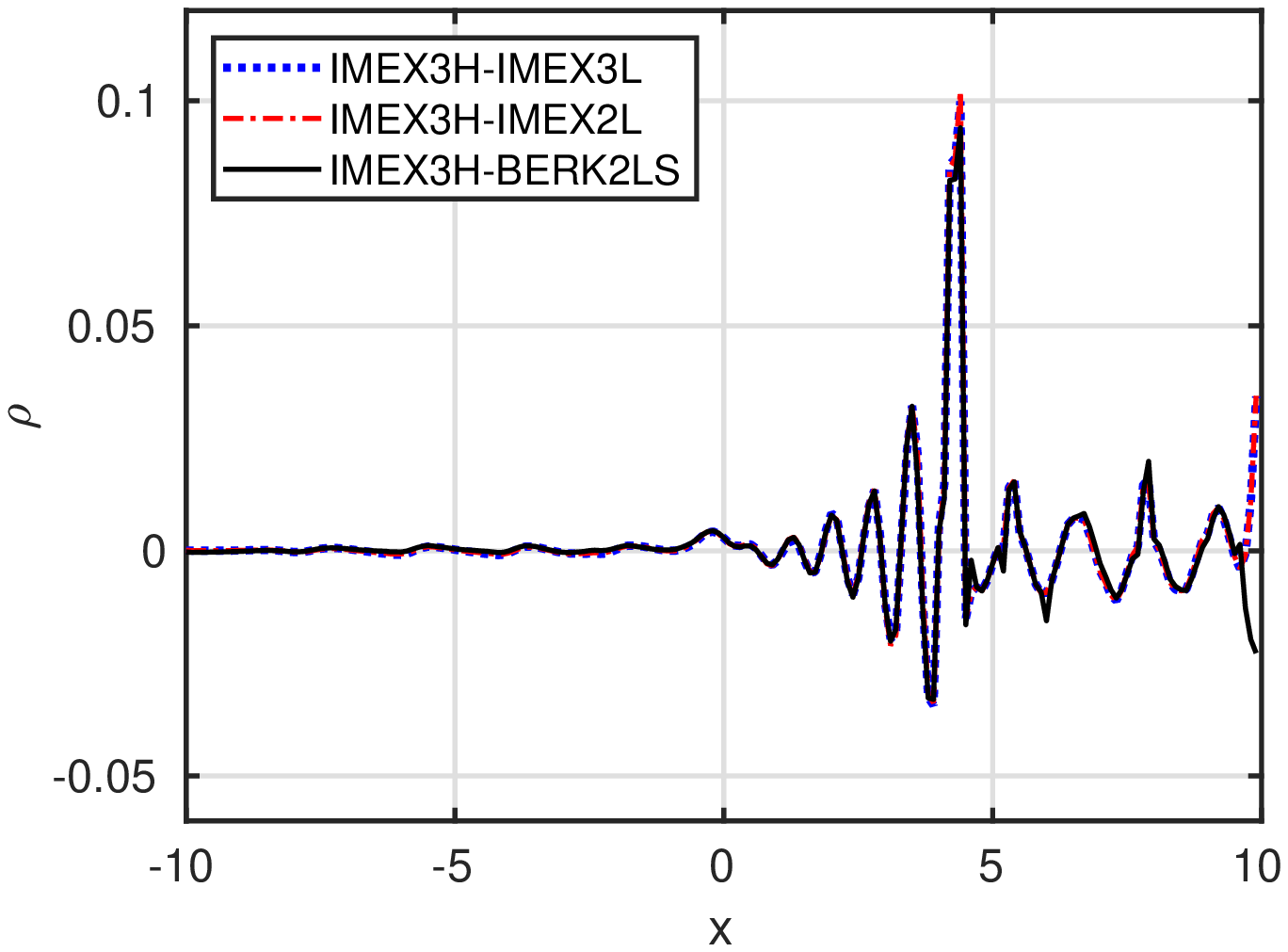}}
  \subfloat[][$u$]{\includegraphics[width=.27\textwidth]{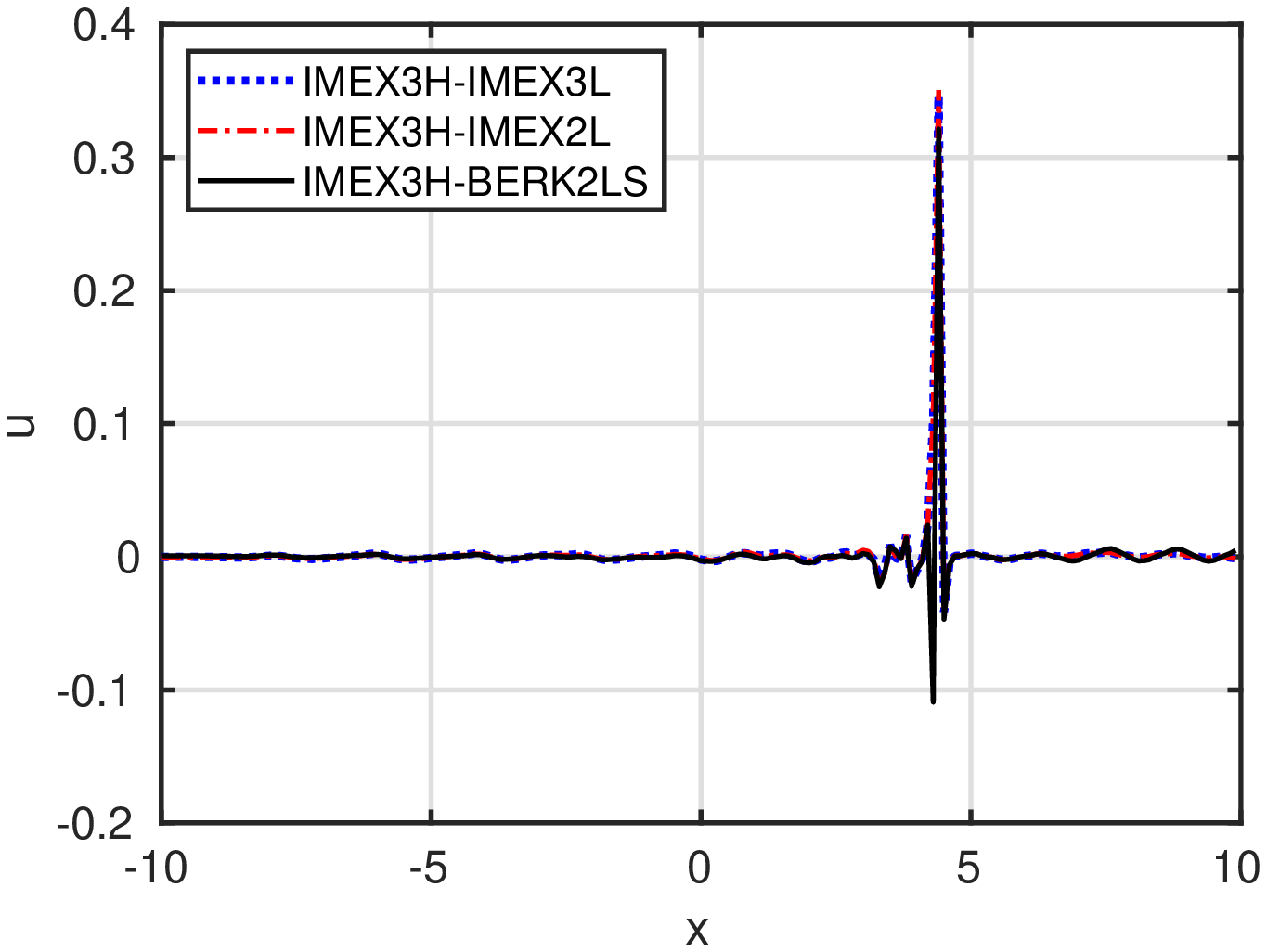}}
  \subfloat[][$\theta$]{\includegraphics[width=.27\textwidth]{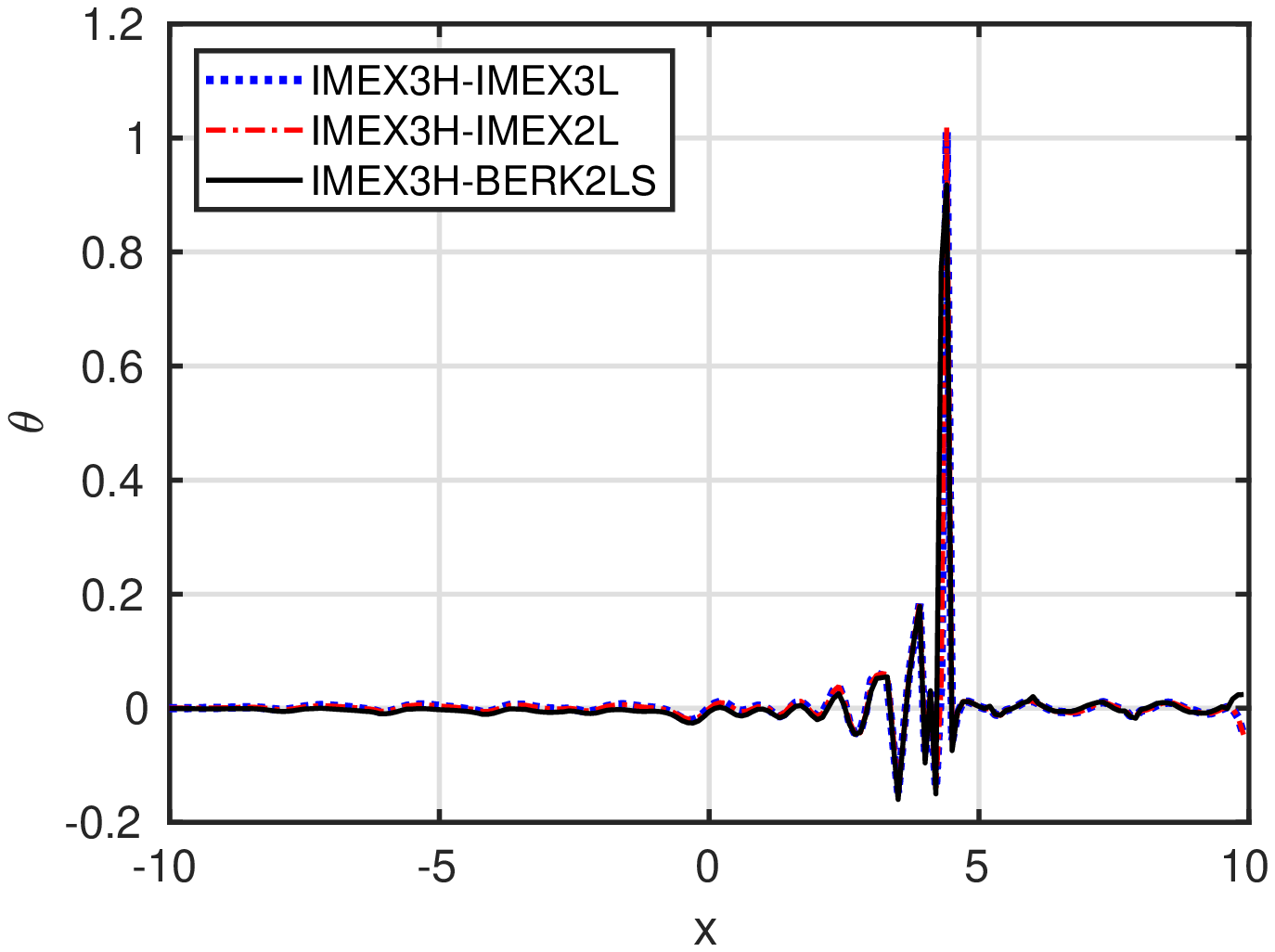}}\par  
  \raisebox{35pt}{\parbox[b]{.16\textwidth}{\tiny$\epsilon=10^{-6}$\\IMEX3L (C=0.14): 291.56s\\IMEX2L (C=0.2): 120.36s\\BERK2L (C=0.2): 29.92s}}%
  \subfloat[][$\rho$]{\includegraphics[width=.27\textwidth]{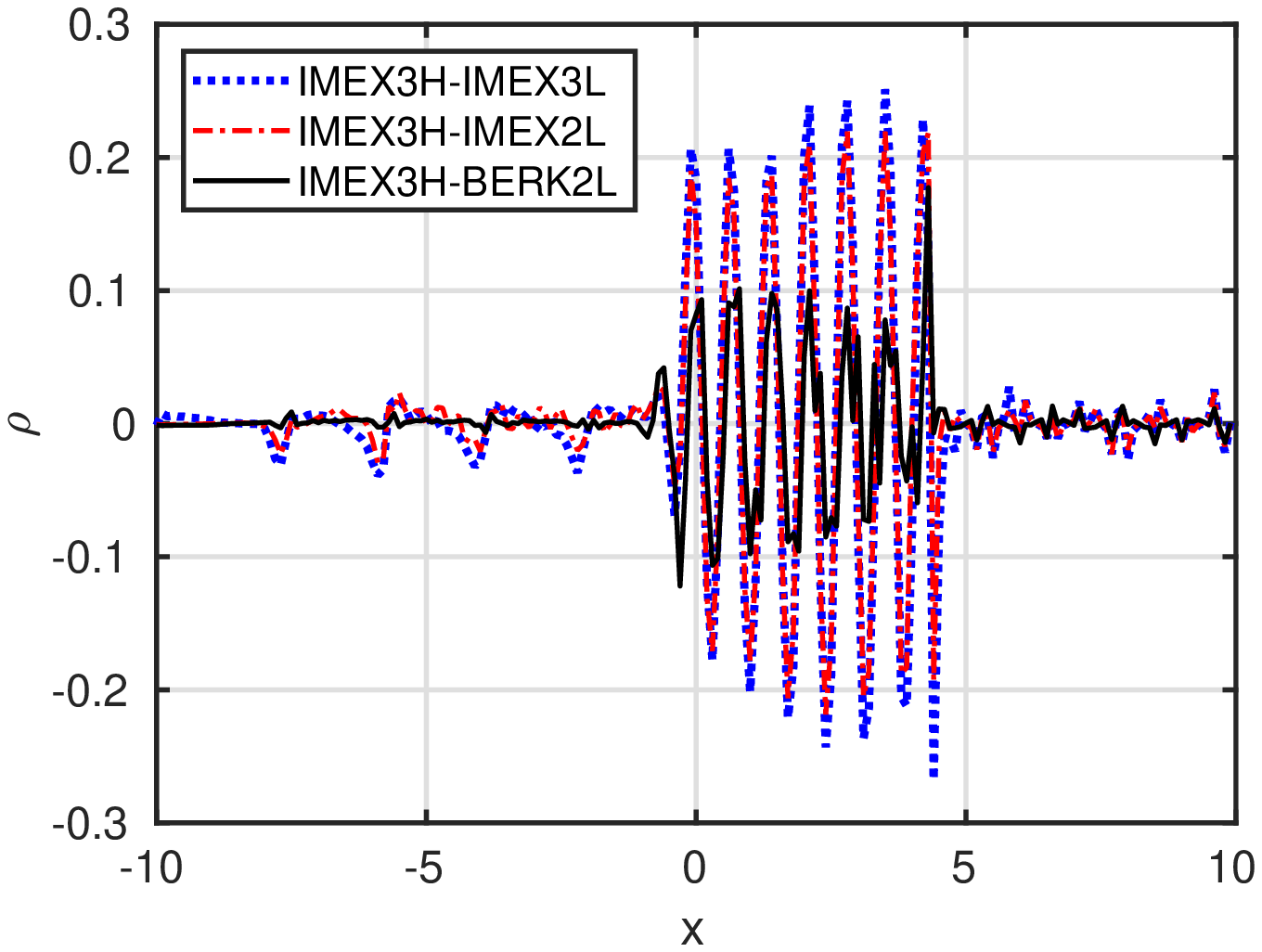}}
  \subfloat[][$u$]{\includegraphics[width=.27\textwidth]{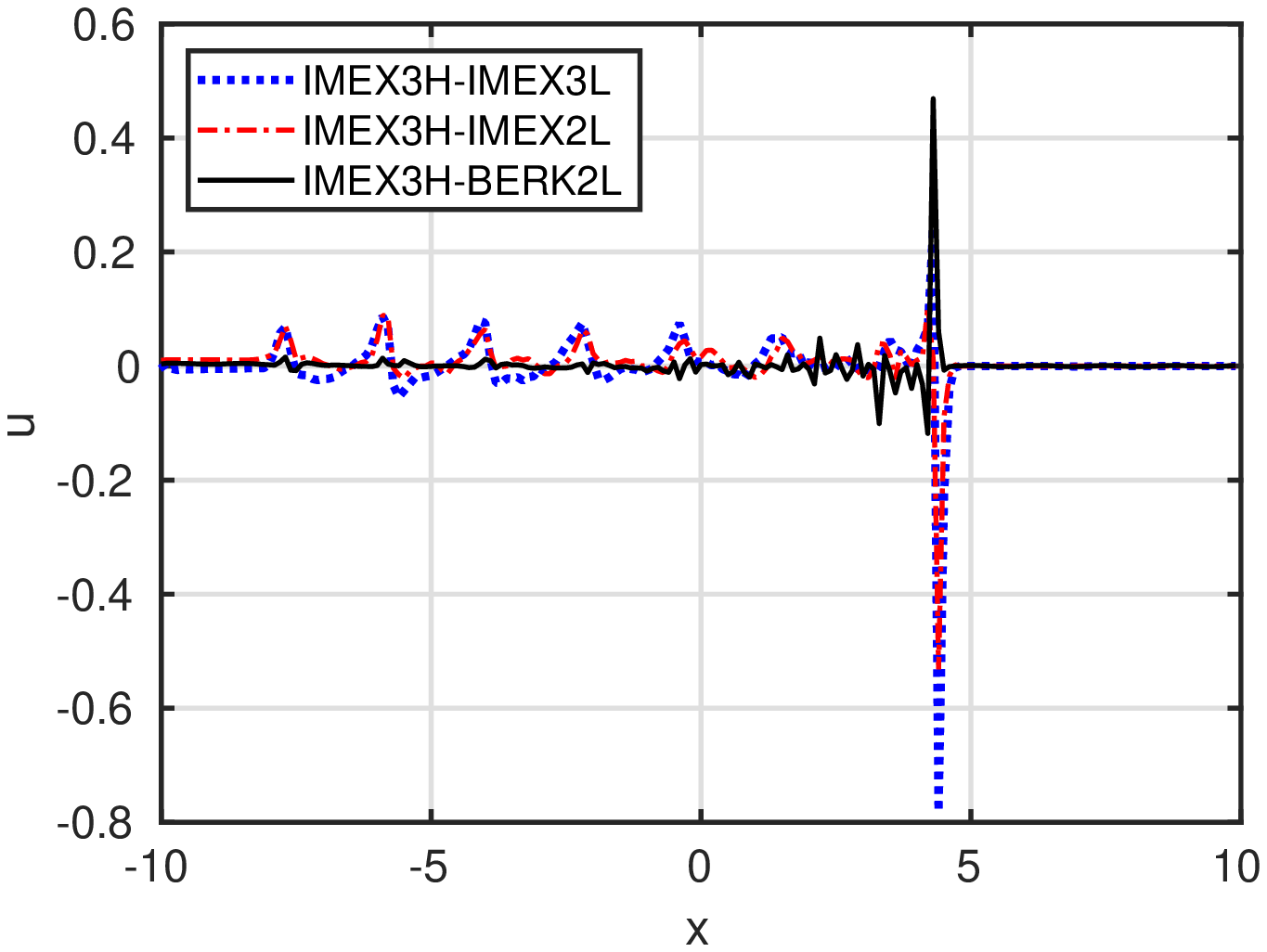}}
  \subfloat[][$\theta$]{\includegraphics[width=.27\textwidth]{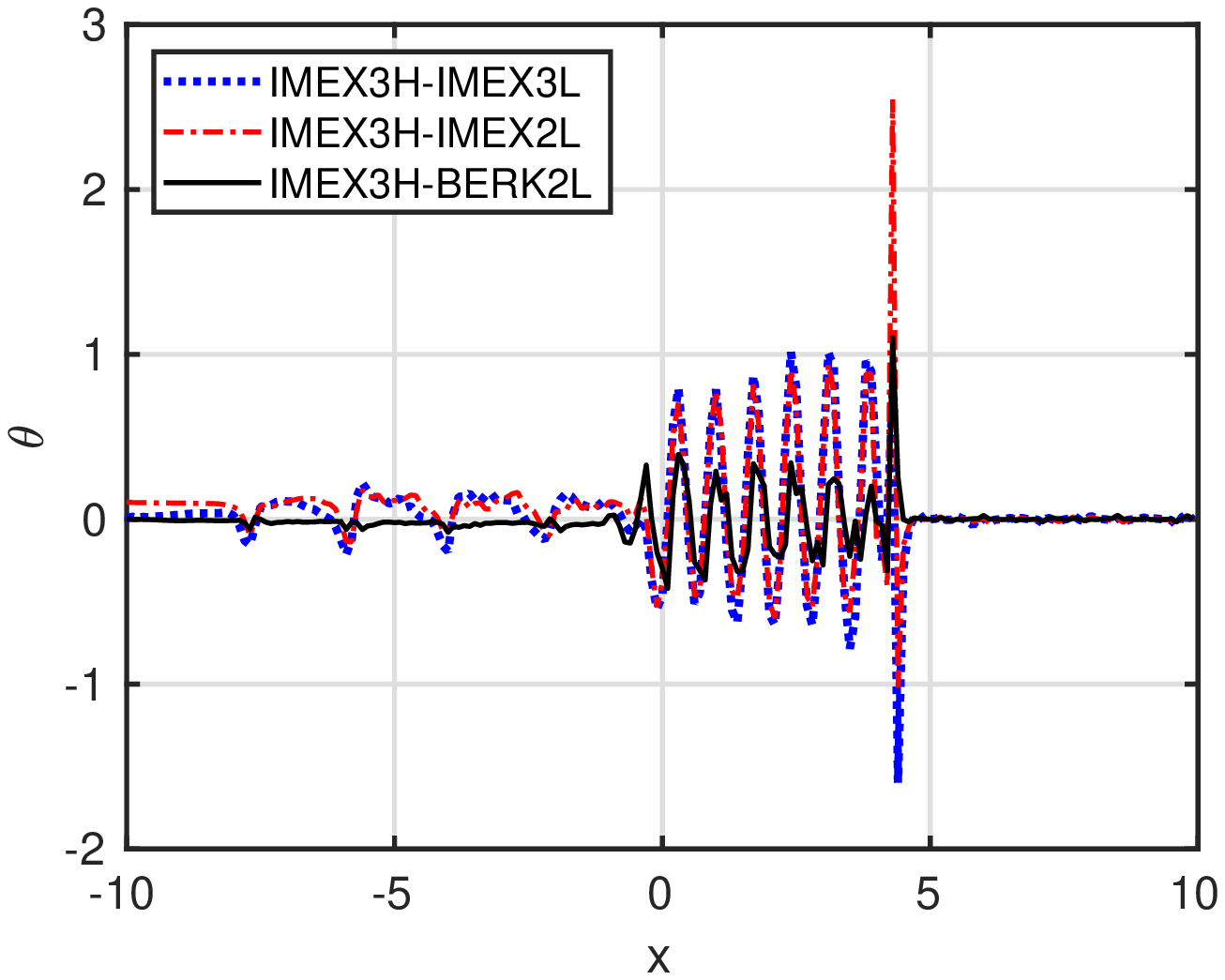}}\par
\caption{Comparison of numerical solutions for the Shu-Osher problem at t=1.8 with various $\epsilon$. Row 1: Euler and IMEX3H at different values of $\epsilon$. Row 2: For $\epsilon =1$, and for each unknown ($\rho$, $u$, $\theta$), we plot the differences; IMEX3H-IMEX3L, IMEX3H–IMEX2L and IMEX3H–BERK2L. Row 3: Same as Row 1, but for $\epsilon=10^{-2}$, Row 4: Same as Row 1, but for $\epsilon=10^{-6}$. The computation time for each method is shown in the first column. Row 5: Comparisons of low-resolution solution for $\epsilon=0.01$. (Parameters for each method - Euler: $N_x=3000$, $C=0.1$;\quad IMEX3H: $N_x=N_v=1000$, $C=0.14$;\quad IMEX3L: $N_x=N_v=200$, $C=0.14$;\quad IMEX2L: $N_x=N_v=200$, $C=0.2$;\quad BERK2L: $N_x=N_v=200$, $C=0.2$)}
  \label{fig:comparison_ShuOsher_diff_plot}
\end{figure}
\begin{figure}[ht!]
  \centering
  \raisebox{35pt}{\parbox[b]{.15\textwidth}{IMEX3H}}%
  \subfloat[][$\rho$]{\includegraphics[width=0.8\textwidth]{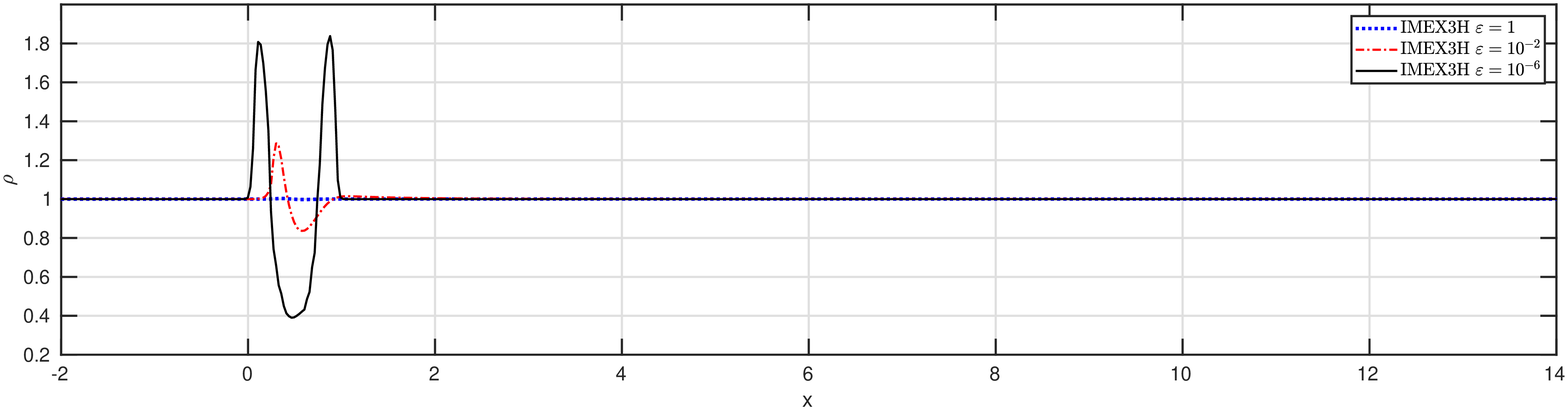}} \\
  \raisebox{35pt}{\parbox[b]{.15\textwidth}{\tiny$\epsilon=1$\\IMEX3L (C=0.1): 697.96s\\IMEX2L (C=0.1): 399.83s\\BERK2L (C=0.1): 3.09s}}\subfloat[][$\rho$]{\includegraphics[width=0.8\textwidth]{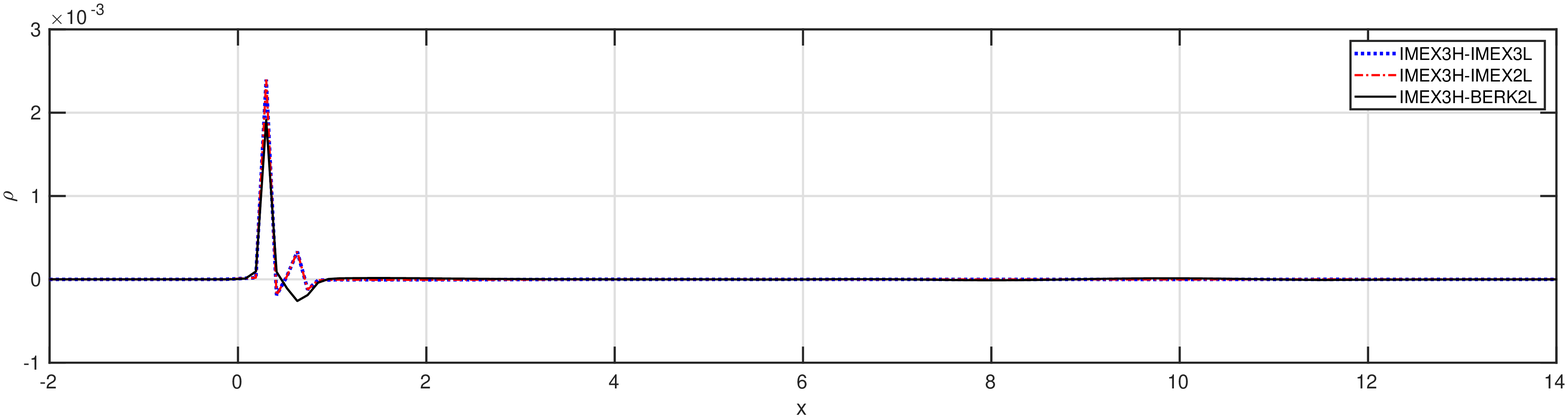}} \\
  \raisebox{35pt}{\parbox[b]{.15\textwidth}{\tiny$\epsilon=10^{-2}$\\IMEX3L (C=0.1): 699.51s\\IMEX2L (C=0.1): 404.01s\\BERK2L (C=0.1): 6.09s}}\subfloat[][$\rho$]{\includegraphics[width=0.8\textwidth]{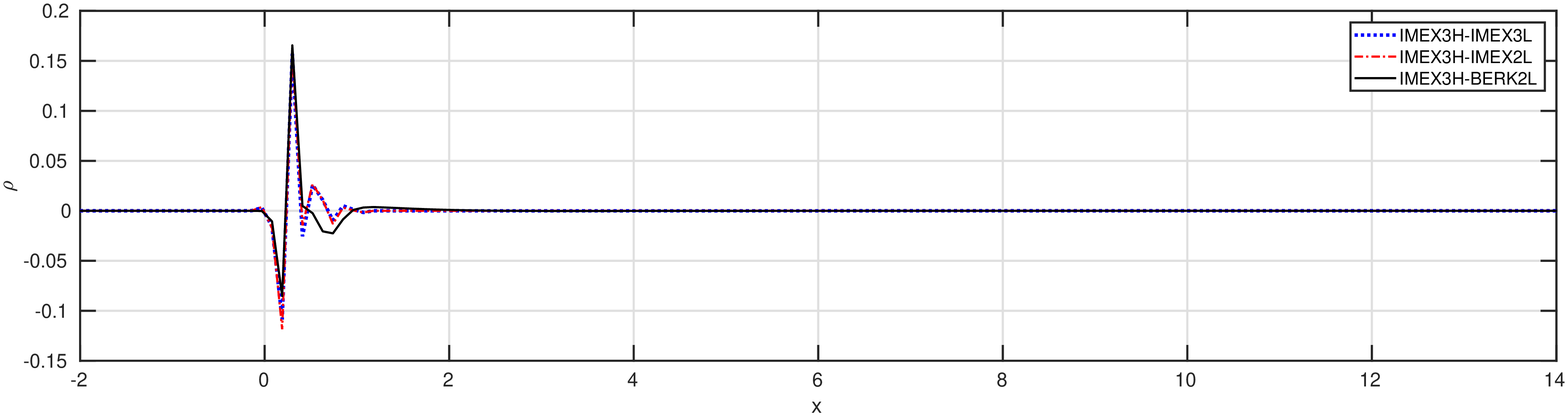}} \\
  \raisebox{35pt}{\parbox[b]{.15\textwidth}{\tiny$\epsilon=10^{-6}$\\IMEX3L (C=0.1): 700.40s\\IMEX2L (C=0.1): 405.88s\\BERK2L (C=0.1): 10.21s}}\subfloat[][$\rho$]{\includegraphics[width=0.8\textwidth]{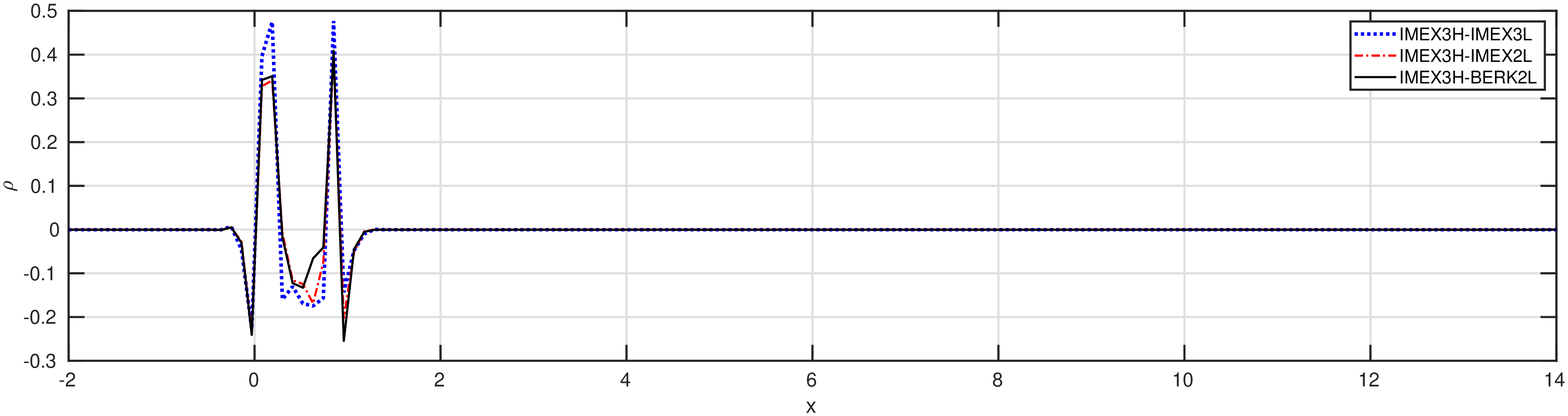}}
  \caption{Comparison of numerical solutions for the gas injection problem at t=0.1 with various $\epsilon$. Row 1: IMEX3H at different values of $\epsilon$. Row 2: For $\epsilon =1$, and for $\rho$, we plot the differences; IMEX3H-IMEX3L, IMEX3H–IMEX2L and IMEX3H–BERK2L. Row 3: Same as Row 1, but for $\epsilon=10^{-2}$, Row 4: Same as Row 1, but for $\epsilon=10^{-6}$. The computation time for each method is shown in the first column. (Parameters for each method - IMEX3H: $N_x=800$, $N_v=1100$, $C=0.1$;\quad IMEX3L: $N_x=200$, $N_v=1000$, $C=0.1$;\quad IMEX2L: $N_x=200$, $N_v=1000$, $C=0.1$;\quad BERK2L: $N_x=200$, $N_v=1000$, $C=0.1$).}
  \label{fig:comparison_gas_injection_diff_wide_rho_plot}
\end{figure}
\begin{figure}[ht!]
  \centering
  \raisebox{35pt}{\parbox[b]{.15\textwidth}{IMEX3H}}%
  \subfloat[][$u$]{\includegraphics[width=0.8\textwidth]{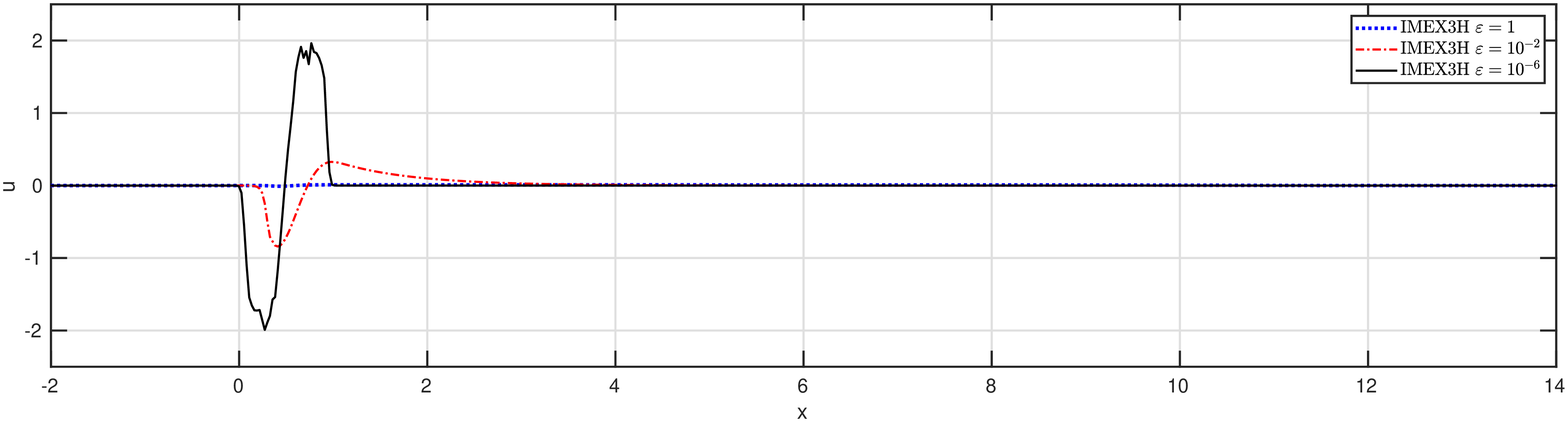}} \\
  \raisebox{35pt}{\parbox[b]{.15\textwidth}{\tiny$\epsilon=1$\\IMEX3L (C=0.1): 697.96s\\IMEX2L (C=0.1): 399.83s\\BERK2L (C=0.1): 3.09s}}\subfloat[][$u$]{\includegraphics[width=0.8\textwidth]{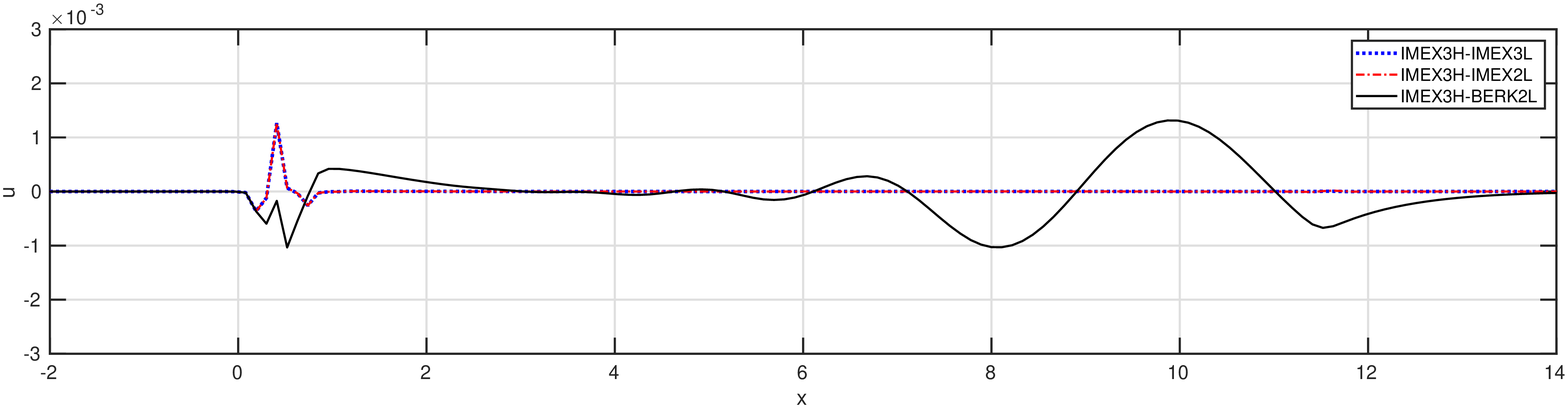}} \\
  \raisebox{35pt}{\parbox[b]{.15\textwidth}{\tiny$\epsilon=10^{-2}$\\IMEX3L (C=0.1): 699.51s\\IMEX2L (C=0.1): 404.01s\\BERK2L (C=0.1): 6.09s}}\subfloat[][$u$]{\includegraphics[width=0.8\textwidth]{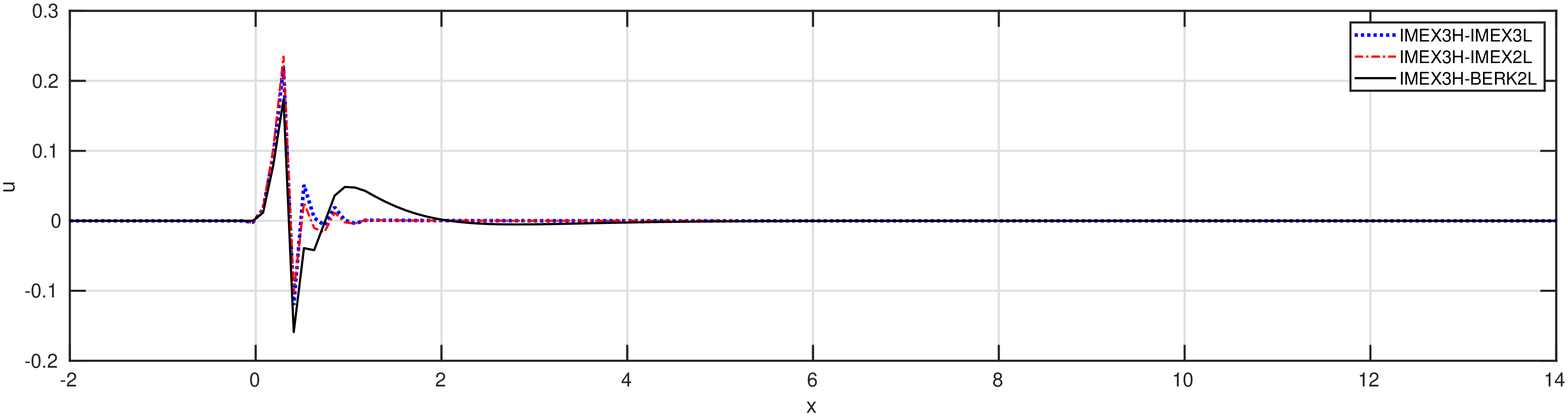}} \\
  \raisebox{35pt}{\parbox[b]{.15\textwidth}{\tiny$\epsilon=10^{-6}$\\IMEX3L (C=0.1): 700.40s\\IMEX2L (C=0.1): 405.88s\\BERK2L (C=0.1): 10.21s}}\subfloat[][$u$]{\includegraphics[width=0.8\textwidth]{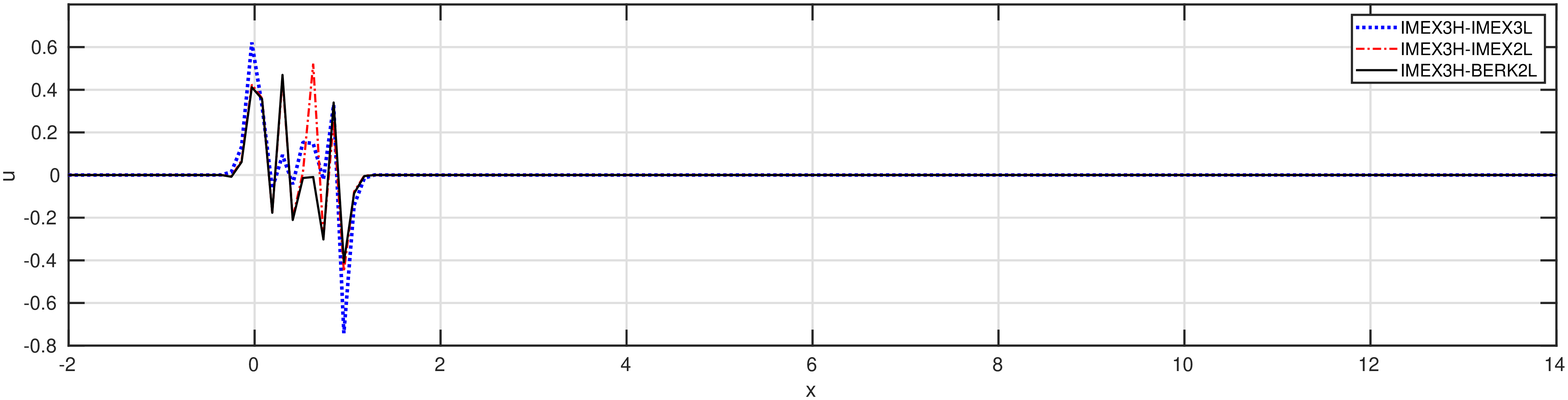}}
  \caption{Comparison of numerical solutions for the gas injection problem at t=0.1 with various $\epsilon$. Row 1: IMEX3H at different values of $\epsilon$. Row 2: For $\epsilon =1$, and for $u$, we plot the differences; IMEX3H-IMEX3L, IMEX3H–IMEX2L and IMEX3H–BERK2L. Row 3: Same as Row 1, but for $\epsilon=10^{-2}$, Row 4: Same as Row 1, but for $\epsilon=10^{-6}$. The computation time for each method is shown in the first column. (Parameters for each method - IMEX3H: $N_x=800$, $N_v=1100$, $C=0.1$;\quad IMEX3L: $N_x=200$, $N_v=1000$, $C=0.1$;\quad IMEX2L: $N_x=200$, $N_v=1000$, $C=0.1$;\quad BERK2L: $N_x=200$, $N_v=1000$, $C=0.1$).}
  \label{fig:comparison_gas_injection_diff_wide_u_plot}
\end{figure}
\begin{figure}[ht!]
  \centering
  \raisebox{35pt}{\parbox[b]{.15\textwidth}{IMEX3H}}%
  \subfloat[][$\theta$]{\includegraphics[width=0.8\textwidth]{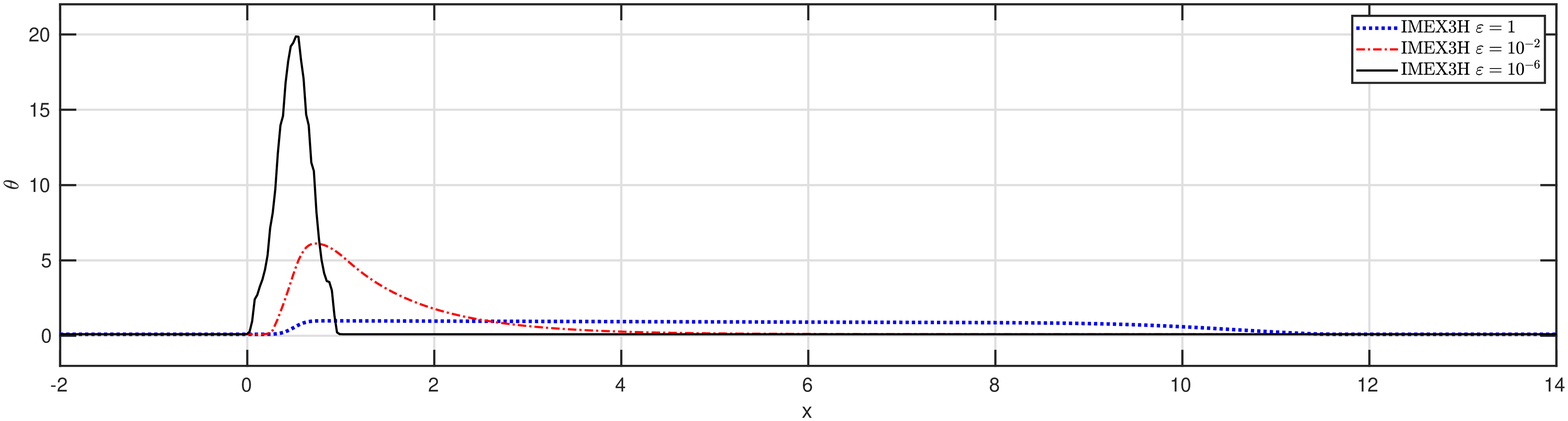}} \\
  \raisebox{35pt}{\parbox[b]{.15\textwidth}{\tiny$\epsilon=1$\\IMEX3L (C=0.1): 697.96s\\IMEX2L (C=0.1): 399.83s\\BERK2L (C=0.1): 3.09s}}\subfloat[][$\theta$]{\includegraphics[width=0.8\textwidth]{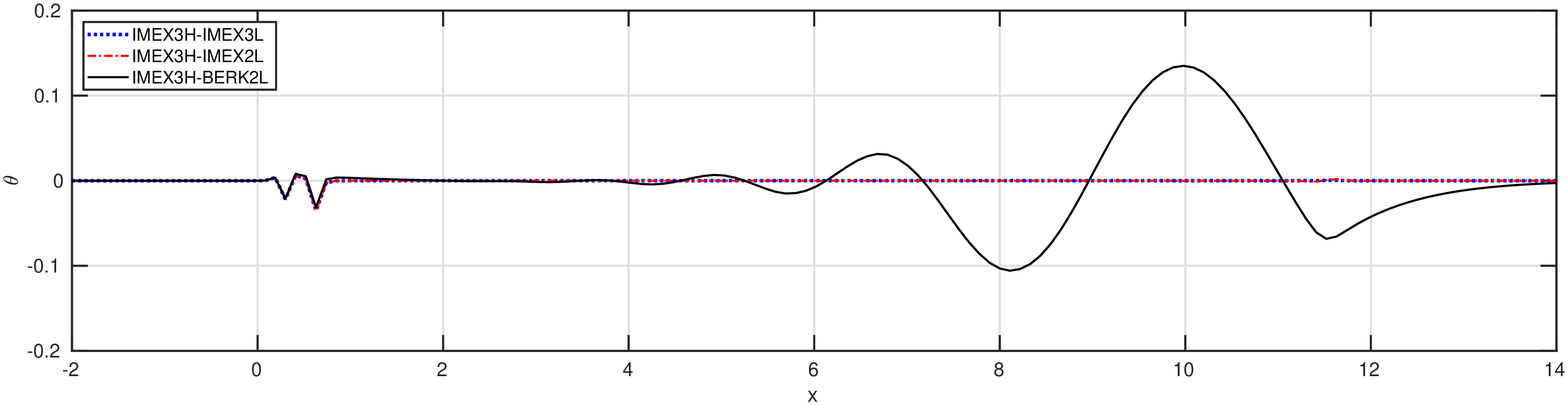}} \\
  \raisebox{35pt}{\parbox[b]{.15\textwidth}{\tiny$\epsilon=10^{-2}$\\IMEX3L (C=0.1): 699.51s\\IMEX2L (C=0.1): 404.01s\\BERK2L (C=0.1): 6.09s}}\subfloat[][$\theta$]{\includegraphics[width=0.8\textwidth]{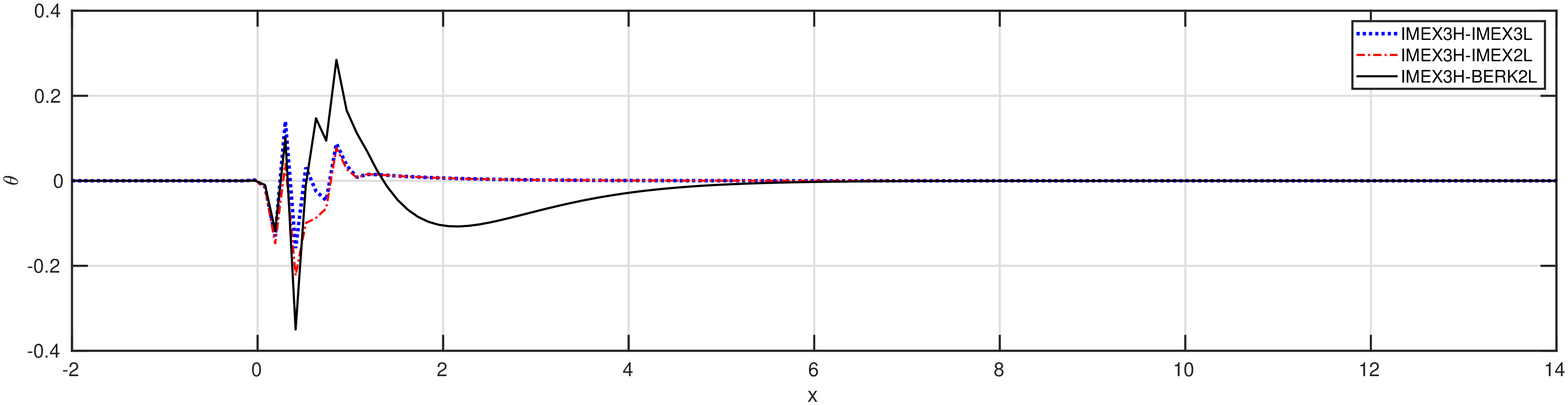}} \\
  \raisebox{35pt}{\parbox[b]{.15\textwidth}{\tiny$\epsilon=10^{-6}$\\IMEX3L (C=0.1): 700.40s\\IMEX2L (C=0.1): 405.88s\\BERK2L (C=0.1): 10.21s}}\subfloat[][$\theta$]{\includegraphics[width=0.8\textwidth]{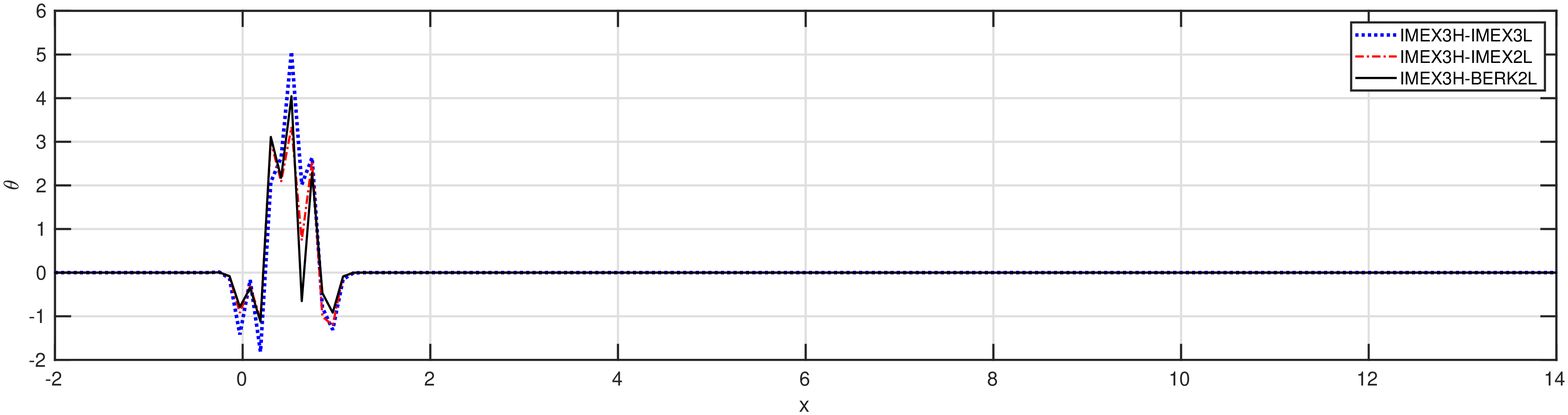}}
  \caption{Comparison of numerical solutions for the gas injection problem at t=0.1 with various $\epsilon$. Row 1: IMEX3H at different values of $\epsilon$. Row 2: For $\epsilon =1$, and for $\theta$, we plot the differences; IMEX3H-IMEX3L, IMEX3H–IMEX2L and IMEX3H–BERK2L. Row 3: Same as Row 1, but for $\epsilon=10^{-2}$, Row 4: Same as Row 1, but for $\epsilon=10^{-6}$. The computation time for each method is shown in the first column. (Parameters for each method - IMEX3H: $N_x=800$, $N_v=1100$, $C=0.1$;\quad IMEX3L: $N_x=200$, $N_v=1000$, $C=0.1$;\quad IMEX2L: $N_x=200$, $N_v=1000$, $C=0.1$;\quad BERK2L: $N_x=200$, $N_v=1000$, $C=0.1$).}
  \label{fig:comparison_gas_injection_diff_wide_T_plot}
\end{figure}
\begin{figure}[ht!]
  \centering
  \raisebox{35pt}{\parbox[b]{.15\textwidth}{IMEX3H}}%
  \subfloat[][$\rho$]{\includegraphics[width=0.8\textwidth]{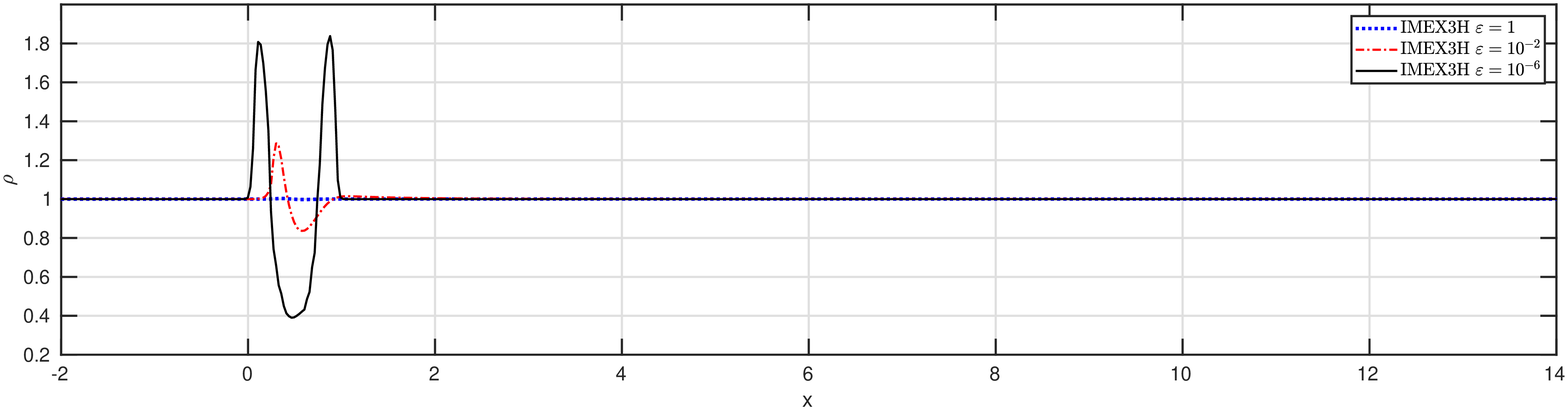}} \\
  \raisebox{35pt}{\parbox[b]{.17\textwidth}{\tiny$\epsilon=1$\\IMEX3L (C=0.1): 697.96s\\IMEX2L (C=0.1): 399.83s\\BERK2LS (C=0.025): 10.33s}}\subfloat[][$\rho$]{\includegraphics[width=0.8\textwidth]{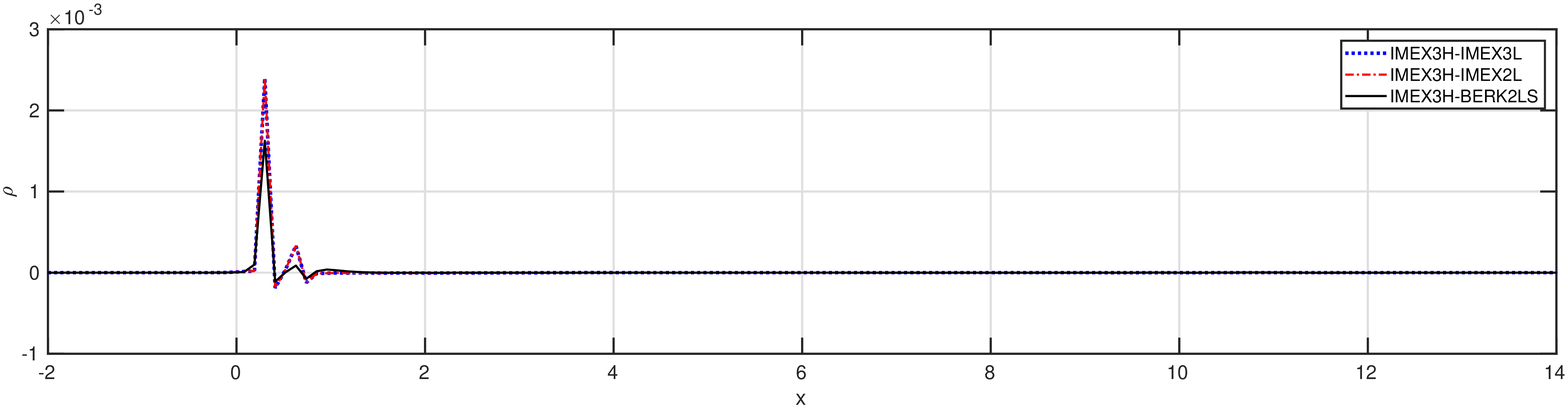}} \\
  \raisebox{35pt}{\parbox[b]{.17\textwidth}{\tiny$\epsilon=10^{-2}$\\IMEX3L (C=0.1): 699.51s\\IMEX2L (C=0.1): 404.01s\\BERK2LS (C=0.025): 18.84s}}\subfloat[][$\rho$]{\includegraphics[width=0.8\textwidth]{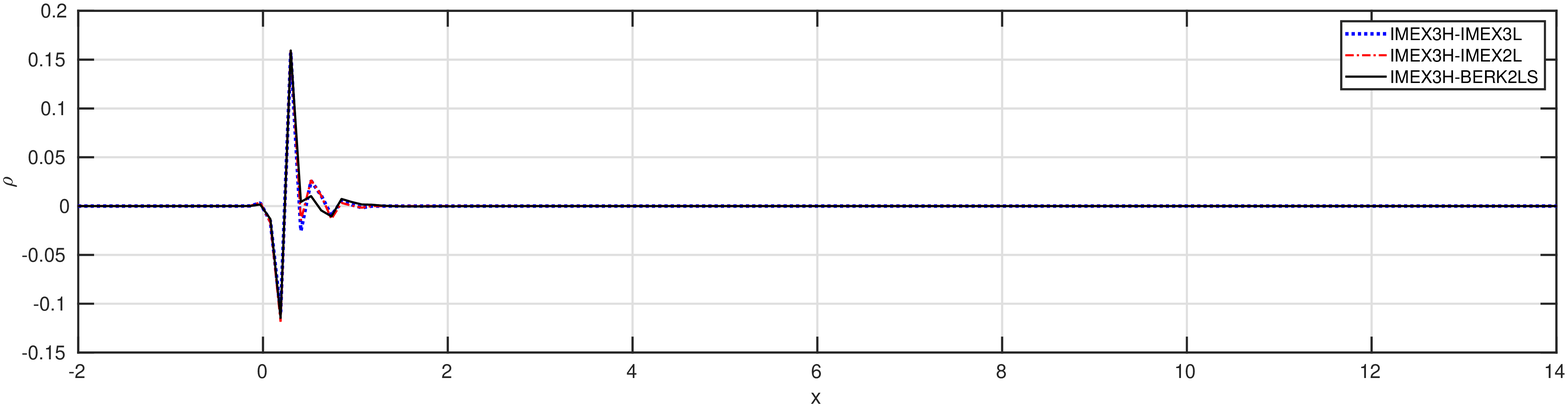}} \\
  \raisebox{35pt}{\parbox[b]{.17\textwidth}{\tiny$\epsilon=10^{-6}$\\IMEX3L (C=0.1): 700.40s\\IMEX2L (C=0.1): 405.88s\\BERK2LS (C=0.025): 33.85s}}\subfloat[][$\rho$]{\includegraphics[width=0.8\textwidth]{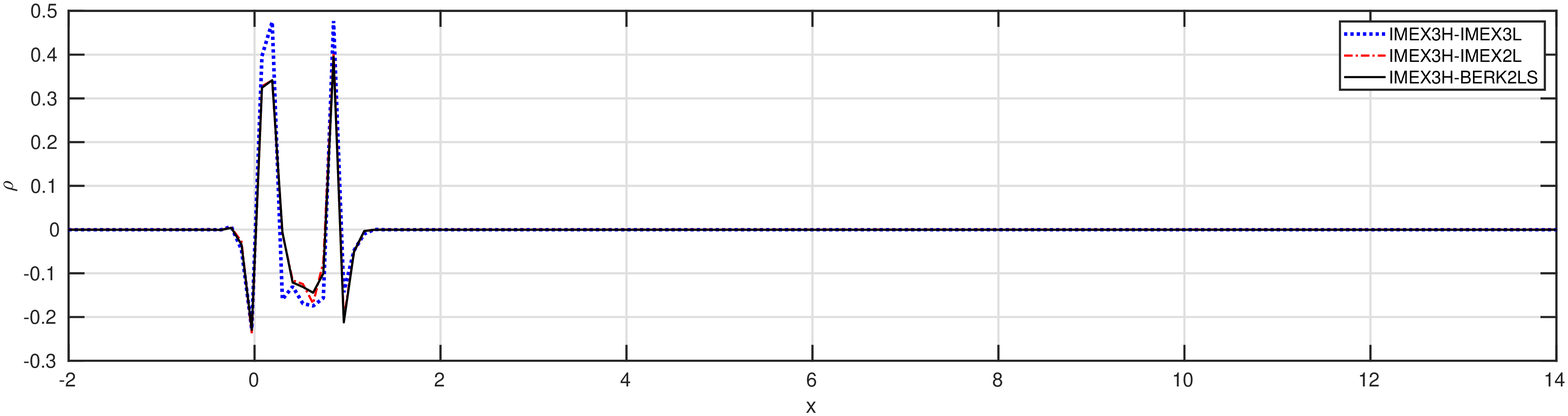}}
  \caption{Comparison of numerical solutions for the gas injection problem at t=0.1 with various $\epsilon$. Row 1: IMEX3H at different values of $\epsilon$. Row 2: For $\epsilon =1$, and for $\rho$, we plot the differences; IMEX3H-IMEX3L, IMEX3H–IMEX2L and IMEX3H–BERK2LS. Row 3: Same as Row 1, but for $\epsilon=10^{-2}$, Row 4: Same as Row 1, but for $\epsilon=10^{-6}$. The computation time for each method is shown in the first column. (Parameters for each method - IMEX3H: $N_x=800$, $N_v=1100$, $C=0.1$;\quad IMEX3L: $N_x=200$, $N_v=1000$, $C=0.1$;\quad IMEX2L: $N_x=200$, $N_v=1000$, $C=0.1$;\quad BERK2LS: $N_x=200$, $N_v=1000$, $C=0.025$).}
  \label{fig:comparison_gas_injection_diff_S_wide_rho_plot}
\end{figure}
\begin{figure}[ht!]
  \centering
  \raisebox{35pt}{\parbox[b]{.15\textwidth}{IMEX3H}}%
  \subfloat[][$u$]{\includegraphics[width=0.8\textwidth]{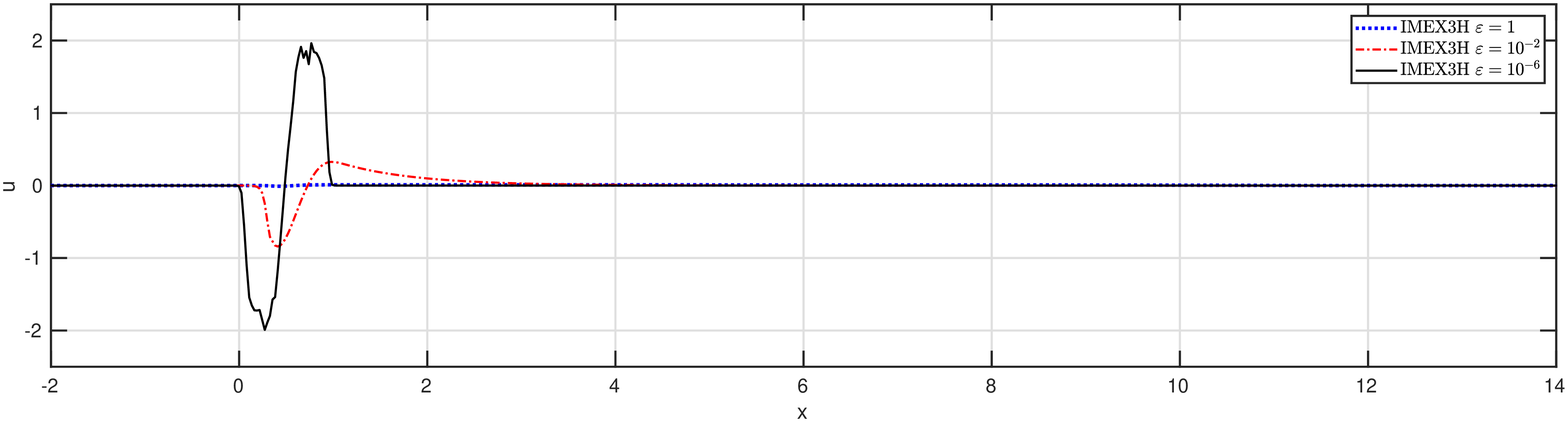}} \\
  \raisebox{35pt}{\parbox[b]{.17\textwidth}{\tiny$\epsilon=1$\\IMEX3L (C=0.1): 697.96s\\IMEX2L (C=0.1): 399.83s\\BERK2LS (C=0.025): 10.33s}}\subfloat[][$u$]{\includegraphics[width=0.8\textwidth]{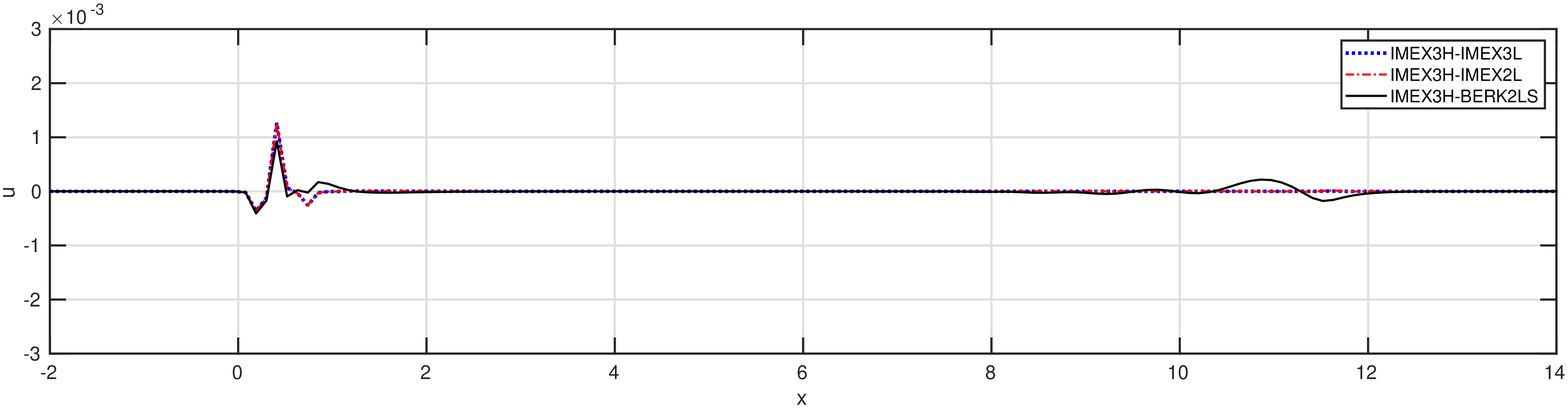}} \\
  \raisebox{35pt}{\parbox[b]{.17\textwidth}{\tiny$\epsilon=10^{-2}$\\IMEX3L (C=0.1): 699.51s\\IMEX2L (C=0.1): 404.01s\\BERK2LS (C=0.025): 18.84s}}\subfloat[][$u$]{\includegraphics[width=0.8\textwidth]{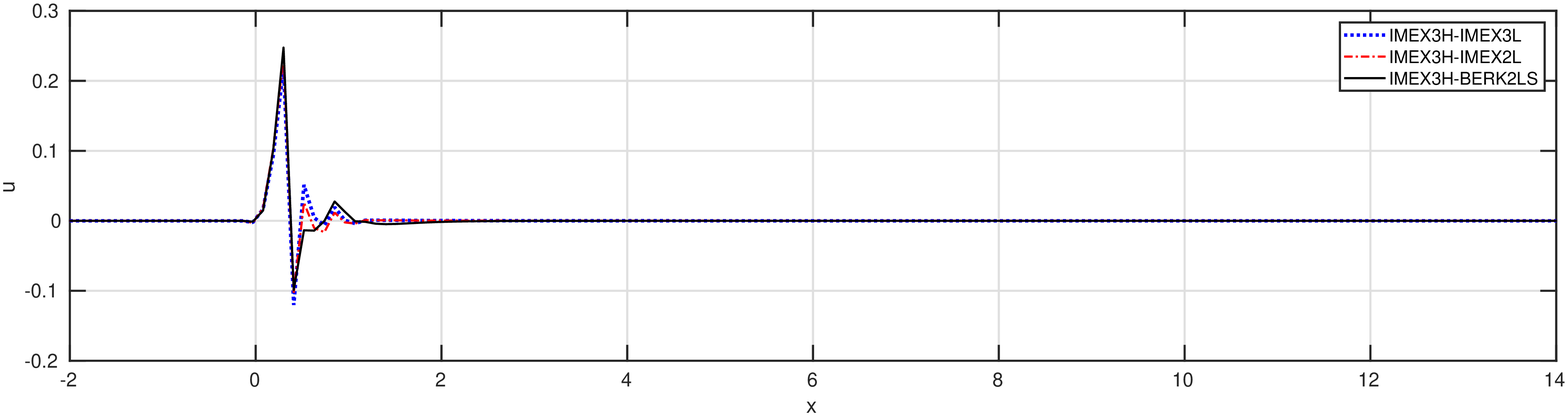}} \\
  \raisebox{35pt}{\parbox[b]{.17\textwidth}{\tiny$\epsilon=10^{-6}$\\IMEX3L (C=0.1): 700.40s\\IMEX2L (C=0.1): 405.88s\\BERK2LS (C=0.025): 33.85s}}\subfloat[][$u$]{\includegraphics[width=0.8\textwidth]{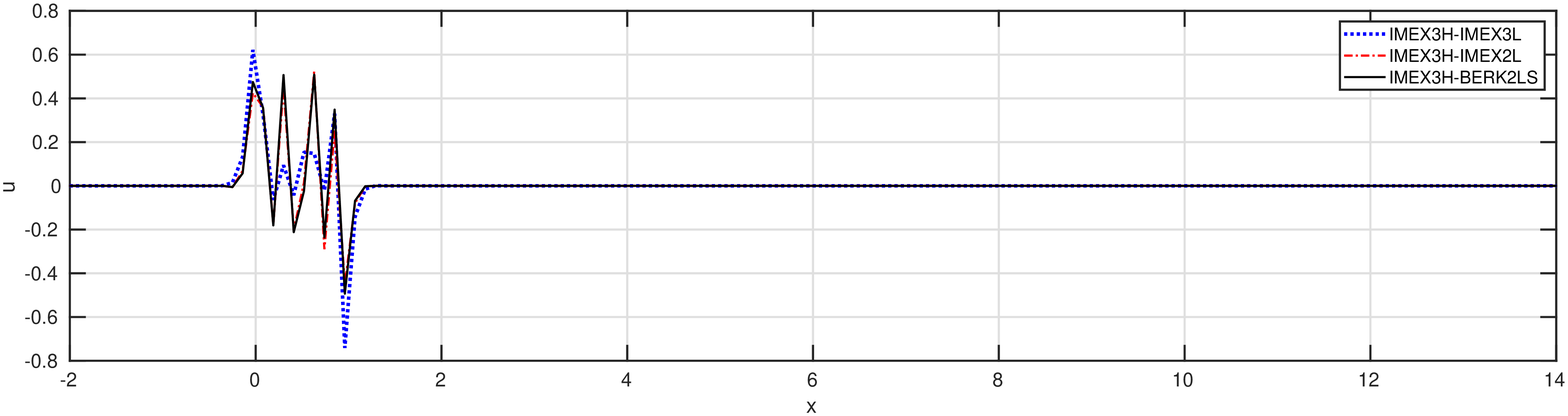}}
  \caption{Comparison of numerical solutions for the gas injection problem at t=0.1 with various $\epsilon$. Row 1: IMEX3H at different values of $\epsilon$. Row 2: For $\epsilon =1$, and for $u$, we plot the differences; IMEX3H-IMEX3L, IMEX3H–IMEX2L and IMEX3H–BERK2LS. Row 3: Same as Row 1, but for $\epsilon=10^{-2}$, Row 4: Same as Row 1, but for $\epsilon=10^{-6}$. The computation time for each method is shown in the first column. (Parameters for each method - IMEX3H: $N_x=800$, $N_v=1100$, $C=0.1$;\quad IMEX3L: $N_x=200$, $N_v=1000$, $C=0.1$;\quad IMEX2L: $N_x=200$, $N_v=1000$, $C=0.1$;\quad BERK2LS: $N_x=200$, $N_v=1000$, $C=0.025$).}
  \label{fig:comparison_gas_injection_diff_S_wide_u_plot}
\end{figure}
\begin{figure}[ht!]
  \centering
  \raisebox{35pt}{\parbox[b]{.15\textwidth}{IMEX3H}}%
  \subfloat[][$\theta$]{\includegraphics[width=0.8\textwidth]{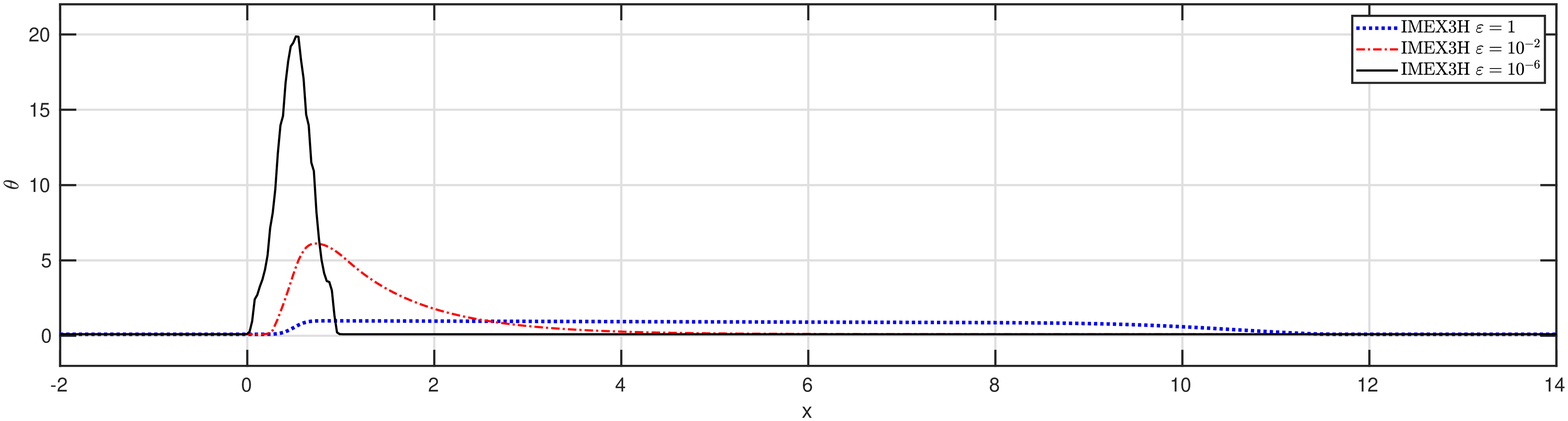}} \\
  \raisebox{35pt}{\parbox[b]{.17\textwidth}{\tiny$\epsilon=1$\\IMEX3L (C=0.1): 697.96s\\IMEX2L (C=0.1): 399.83s\\BERK2LS (C=0.025): 10.33s}}\subfloat[][$\theta$]{\includegraphics[width=0.8\textwidth]{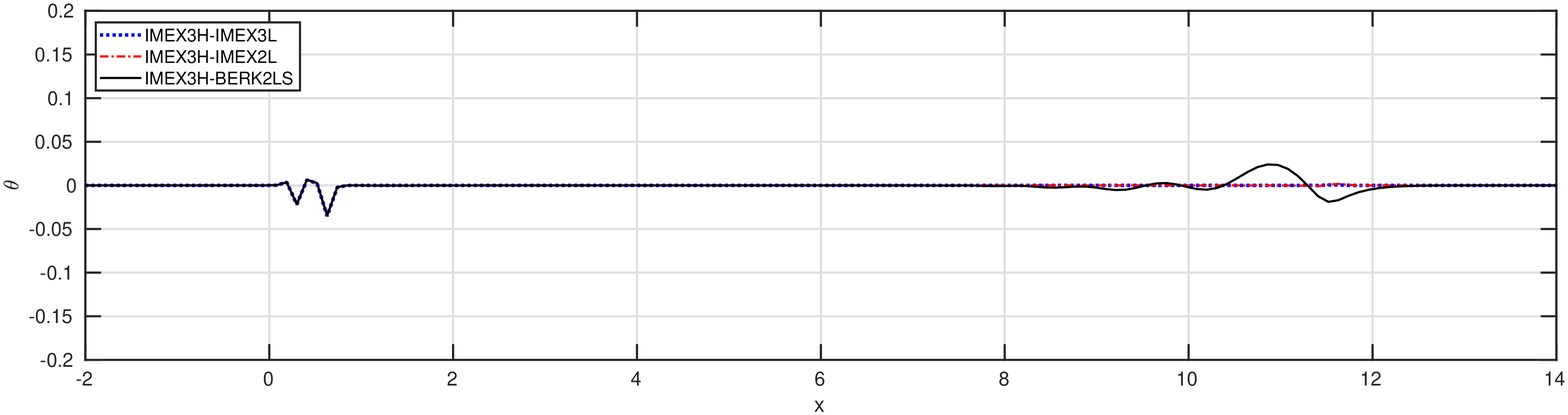}} \\
  \raisebox{35pt}{\parbox[b]{.17\textwidth}{\tiny$\epsilon=10^{-2}$\\IMEX3L (C=0.1): 699.51s\\IMEX2L (C=0.1): 404.01s\\BERK2LS (C=0.025): 18.84s}}\subfloat[][$\theta$]{\includegraphics[width=0.8\textwidth]{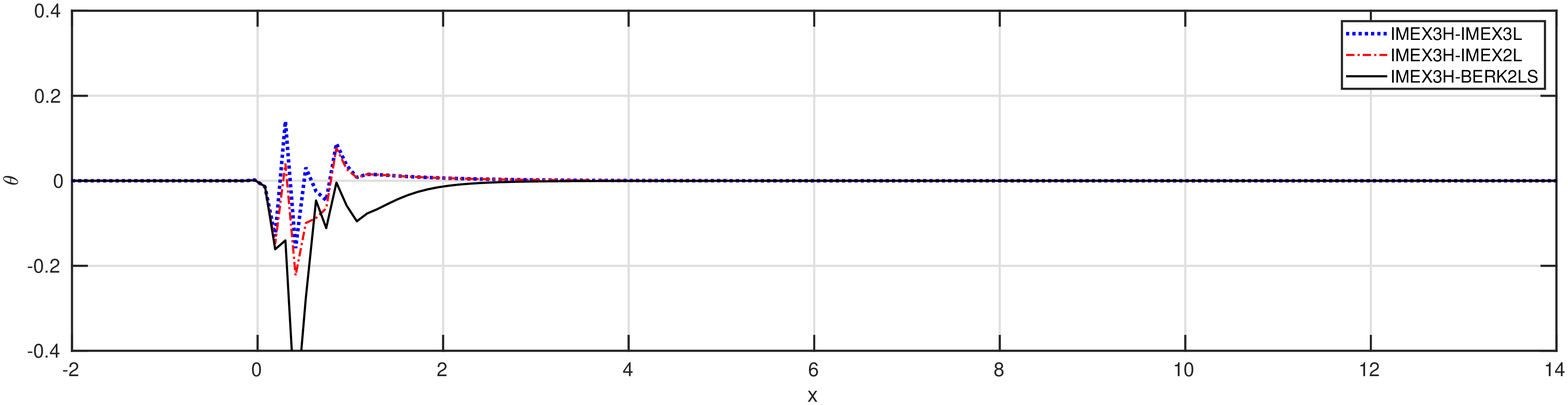}} \\
  \raisebox{35pt}{\parbox[b]{.17\textwidth}{\tiny$\epsilon=10^{-6}$\\IMEX3L (C=0.1): 700.40s\\IMEX2L (C=0.1): 405.88s\\BERK2LS (C=0.025): 33.85s}}\subfloat[][$\theta$]{\includegraphics[width=0.8\textwidth]{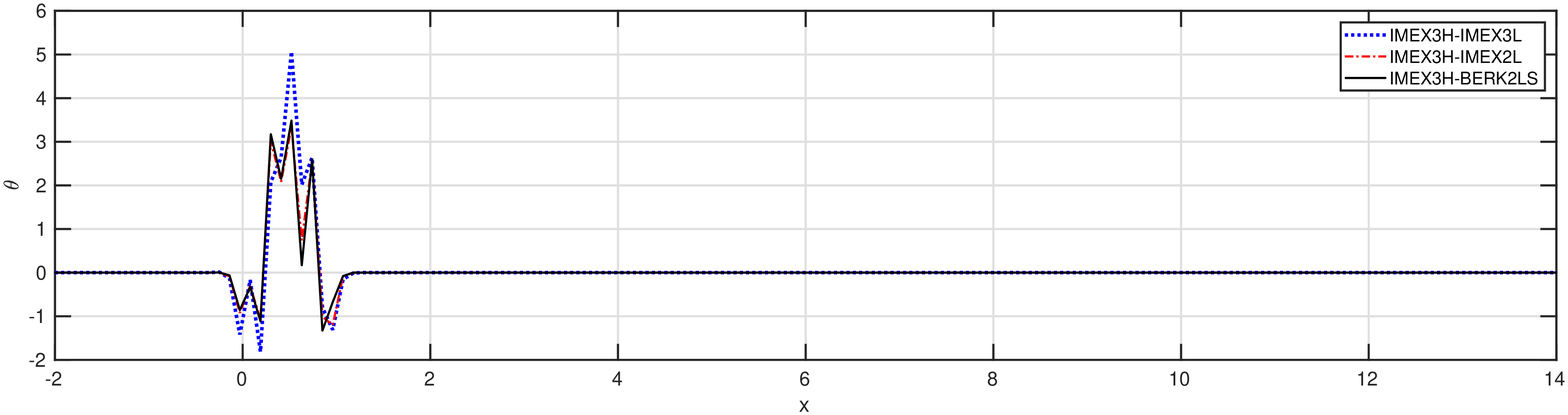}}
  \caption{Comparison of numerical solutions for the gas injection problem at t=0.1 with various $\epsilon$. Row 1: IMEX3H at different values of $\epsilon$. Row 2: For $\epsilon =1$, and for $\theta$, we plot the differences; IMEX3H-IMEX3L, IMEX3H–IMEX2L and IMEX3H–BERK2LS. Row 3: Same as Row 1, but for $\epsilon=10^{-2}$, Row 4: Same as Row 1, but for $\epsilon=10^{-6}$. The computation time for each method is shown in the first column. (Parameters for each method - IMEX3H: $N_x=800$, $N_v=1100$, $C=0.1$;\quad IMEX3L: $N_x=200$, $N_v=1000$, $C=0.1$;\quad IMEX2L: $N_x=200$, $N_v=1000$, $C=0.1$;\quad BERK2LS: $N_x=200$, $N_v=1000$, $C=0.025$).}
  \label{fig:comparison_gas_injection_diff_S_wide_T_plot}
\end{figure}